\documentclass[a4paper,11pt]{article}

\usepackage{fancyhdr}
\usepackage[centertags]{amsmath}
\usepackage{amsfonts}
\usepackage{graphics}
\usepackage{graphicx}
\usepackage{booktabs}
\usepackage{hyperref}
\usepackage{amssymb}
\usepackage{amsthm}
\usepackage{newlfont}
\usepackage{url}
\usepackage{color}
\usepackage{epsfig}
\usepackage{ifthen}
\usepackage{ifpdf}
\usepackage{caption}
\usepackage{amsmath}
\usepackage[ruled]{algorithm}
\usepackage{algorithmic}
\usepackage{multirow,multicol}
\usepackage{subcaption}

\DeclareMathOperator{\prox}{prox}

\def\st{\mbox{s.t.}}
\newcommand{\N}{\mathbb{N}}
\newcommand{\R}{\mathbb{R}}
\newcommand{\C}{\mathbb{C}}
\def\bb{\mathbf b}

\def\bd{\mathbf d}
\def\be{\mathbf e}
\def\bu{\mathbf u}
\def\bv{\mathbf v}
\def\bw{\mathbf w}
\def\bx{\mathbf x}
\def\by{\mathbf y}
\def\bz{\mathbf z}
\def\b0{\mathbf 0}

\def\bgamma{\boldsymbol \gamma}
\def\bmu{\boldsymbol \mu}
\def\bxi{\boldsymbol \xi}

\def\mE{\mathcal{E}}
\def\mF{\mathcal{F}}
\def\mL{\mathcal{L}}

\newtheorem{theorem}{Remark}
\newtheorem{Assumption}{Assumption}[theorem]
\newtheorem{Remark}{Remark}[theorem]

\ifpdf
\DeclareGraphicsExtensions{.pdf,.png,.jpg}
\else
\DeclareGraphicsExtensions{.eps}
\fi
\graphicspath{{figures/}}

\title{Directional TGV-based image restoration under Poisson noise\footnote{This research was partially supported by the Istituto Nazionale di Alta Matematica, Gruppo Nazionale per il Calcolo Scientifico (INdAM-GNCS). D.~di~Serafino and M.~Viola were also funded by the V:ALERE Program of the University of Campania ``L. Vanvitelli''.}}
\author{Daniela di Serafino$^{1}$, Germana Landi$^{2}$, Marco Viola$^{3}$\\[10pt]
    {\footnotesize$^{1}$ Department of Mathematics and Applications ``R. Caccioppoli'',}\\
    {\footnotesize University of Naples Federico II, Naples, Italy; \texttt{daniela.diserafino@unina.it}}\\
    {\footnotesize$^{2}$ Department of Mathematics, University of Bologna,}\\
    {\footnotesize Bologna, Italy; \texttt{germana.landi@unibo.it}}\\
    {\footnotesize$^{3}$ Department of Mathematics and Physics,}\\
    {\footnotesize University of Campania ``L. Vanvitelli'', Caserta, Italy; \texttt{marco.viola@unicampania.it}}
}

\date{\textsc{version 2} -- June 11, 2021}
 
\begin{document}
    \maketitle
    \begin{abstract}
        We are interested in the restoration of noisy and blurry images where the texture
        mainly follows a single direction (i.e., directional images). Problems of this type arise, for example,
        in microscopy or computed tomography for carbon or glass fibres. {In order to deal with these problems,
        the Directional Total Generalized Variation (DTGV) was developed by Kongskov et al.~in 2017 and 2019,
        in the case of impulse and Gaussian noise.}
        In this article we focus on images corrupted by Poisson noise,
        extending the DTGV regularization to image restoration models where the data fitting term is the generalized
        Kullback--Leibler divergence. We also propose a technique for the identification of the main texture direction,
        which improves upon the techniques used in the aforementioned work about DTGV. We solve the problem by an
        ADMM algorithm with proven convergence and subproblems that can be solved exactly at a low computational cost.
        Numerical results on both phantom and real images demonstrate the effectiveness of our approach.\\[10pt]
        {\textbf{Keywords:} directional image restoration; Poisson noise; DTGV regularization; ADMM method}
    \end{abstract}

\section{Introduction}

Poisson noise appears in processes where digital images are obtained by the count of particles
(generally photons). This is the case of X-ray computed tomography,
positron emission tomography, confocal and fluorescence microscopy and optical/infrared
astronomical imaging, to name just a few applications (see, e.g., \cite{bertero:2009} and the references
therein). In this case, the object to be restored can be represented as a vector $\bu \in \R^n$ 
and the data can be assumed to be a vector {$\bb \in \N_0^n$},
whose entries $b_j$ are sampled from {$n$} independent Poisson random variables $B_j$ with probability 

$$
P(B_j = b_j) = \frac{e^{-(A \bu + \bgamma)_j}(A\bu + \bgamma)_ {j}^{b_j}}{b_j!}.
$$

The matrix {$A = (a_{ij}) \in \R^{n \times n}$} models the observation mechanism of the imaging system
and the following standard assumptions are made:

\begin{equation}
    \label{eq:matrixa}
    a_{ij} \geq 0 \mbox{ for all }  i, j, \qquad \sum_{i=1}^{n} a_{ij} = 1 \mbox{ for all } j.
\end{equation}

The vector $\bgamma \in \R^{n}$,  with $\bgamma > 0$, models the background radiation detected by \mbox{the sensors.}

By applying a maximum-likelihood approach~\cite{bertero:2009,shepp:1982}, we can estimate $\bu$ by minimizing
the Kullback--Leibler (KL) divergence of $A \bu+\bgamma$ from $\bb$:

\begin{equation} \label{eq:kl}
    D_{KL}(A \bu + \bgamma, \bb) = \sum_{i=1}^{n} \left( b_i \ln \frac{b_i}{[A \bu+\bgamma]_i} + [A \bu+\bgamma]_i - b_i \right),
\end{equation}

\noindent where we set

$$
b_i \ln \frac{b_i}{[A \bu+\bgamma]_i} = 0 \quad \mbox{if } \; b_i = 0.
$$

Regularization is usually introduced in~\eqref{eq:kl} to deal with the ill-conditioning of this problem.
The Total Variation (TV) regularization~\cite{rudin:1992} has been widely used in this context, because
it preserves edges and is able to smooth flat areas of the image. However, since it may produce staircase artifacts,
other TV-based regularizers have been proposed. For example, the Total Generalized Variation (TGV) has been
proposed and applied in \cite{bredies:2010, bredies:2014, gao:2018, diserafino:2021conf} to overcome the staircasing effect while keeping
the ability of identifying edges. On the other hand, to improve the quality of restoration for directional images,
the Directional TV (DTV) regularization has been considered in \cite{bayram:2012}, in the discrete setting.
In \cite{kongskov:2017, kongskov:2019}, a regularizer combining DTV and TGV, named Directional
TGV (DTGV), has been successfully applied to directional images affected by impulse and Gaussian noise.

{
    Given an image $\bu\in\R^n$, the discrete second-order Directional TGV of $\bu$ is defined as
    \begin{equation}\label{eqn:dtgv2}
        \mathrm{DTGV}^2(\bu) = \min\limits_{\bw\in\R^{2n}} \alpha_0 \left\|\widetilde{\nabla} \bu - \bw\right\|_{2,1|\R^{2n}} + \alpha_1 \left\|\widetilde{\mE}\bw\right\|_{2,1|\R^{4n}},
    \end{equation}
    where $\bw\in\R^{2n}$, $\widetilde{\nabla}\in\R^{2n\times n}$
    and $\widetilde{\mE}\in\R^{4n\times2n}$ are the {discrete directional gradient operator} and the {directional symmetrized derivative}, respectively, and $\alpha_0,\alpha_1\in(0,\,+\infty)$. For any vector $\bv\in\R^{2n}$ we set 
    
    \begin{equation}\label{eq:def_21norm_z2}
        \|\bv\|_{2,1|\R^{2n}} = \sum_{j=1}^n \sqrt{v_j^2+v_{n+j}^2} ,
    \end{equation}
    
    
    \noindent
    and for any vector $\by \in \R^{4n}$ we set
    
    \begin{equation}\label{eq:def_21norm_z3}
        \| \by \|_{2,1 | \R^{4n}} = \sum_{j=1}^n \sqrt{y_j^2 + y_{n+j}^2 + y_{2n+j}^2 + y_{3n+j}^2}.
    \end{equation}
    
    
    Given an angle $\theta\in[-\pi,\pi]$ and a scaling parameter $a>0$, we have that the discrete directional gradient operator has the form
    
    \begin{equation*}
        \widetilde{\nabla} = \left[\begin{array}{c}
            D_\theta \\ D_{\theta^\bot} \end{array}\right] =
        \left[\begin{array}{c} \cos(\theta) D_H +  \sin(\theta) D_V \\ a\left(-\sin(\theta) D_H +  \cos(\theta) D_V\right)
        \end{array}\right],
    \end{equation*}
    
    \noindent
    where $D_\theta, D_{\theta^\bot}\in\R^{n\times n}$ represent the forward finite-difference operators
    along the directions determined by $\theta$ and $\theta^\bot=\theta+\frac{\pi}{2}$, respectively,
    and $D_H, D_V \in\R^{n\times n}$ represent the forward finite-difference operators along the
    horizontal and the vertical direction, respectively. Moreover, the directional symmetrized derivative is defined in block-wise \mbox{form as}
    \begin{equation*}
        \widetilde{\mE} = \left[\begin{array}{cc}
            D_\theta & 0 \\
            \frac{1}{2}D_{\theta^\bot} & \frac{1}{2}D_{\theta}\\
            \frac{1}{2}D_{\theta^\bot} & \frac{1}{2}D_{\theta}\\
            0 & D_{\theta^\bot}
        \end{array}\right].
    \end{equation*}
    
    It is worth noting that, by fixing $\theta=0$ and $a=1$, we have $D_\theta = D_H$ and $D_{\theta^\bot} = D_V$, and the operators
    $\widetilde{\nabla}$ and $\widetilde{\mE}$ define the TGV$^2$ regularization~\cite{bredies:2010}.
    
    We observe that the definition of both the matrix $A$ and the finite difference operators $D_H$ and $D_V$ depend on the choice of boundary conditions. We make the \mbox{following assumption.}
    \begin{Assumption}\label{ass:bccb_structure}
        We assume that periodic boundary conditions are considered for $A$, $D_H$ and $D_V$.
        Therefore, those matrices are Block Circulant with Circulant Blocks (BCCB).
    \end{Assumption}
}

In this work we focus on directional images affected by Poisson noise, with the aim of assessing the behaviour of DTGV in this case.
Besides extending the use of DTGV to Poisson noise, we introduce a novel technique for estimating the main direction
of the image, which appears to be more efficient than the techniques applied in \cite{kongskov:2017, kongskov:2019}.
We solve the resulting optimization problem by using a customized version of the Alternating Direction Method of Multipliers (ADMM).
We note that all the ADMM subproblems can be solved exactly at a low cost, thanks also to the use of FFTs, and that the method
has proven convergence. Finally, we show the effectiveness of our approach on a set of test images, corrupted
by out-of-focus and Gaussian blurs and noise with different signal-to-noise ratios. In particular, the KL-DTGV model of our problem
is described in Section~\ref{sec:kl-dtgv} and the technique for estimating the main direction is presented in Section~\ref{sec:direction}.
A detailed description of the ADMM version used for the minimization is given in Section~\ref{sec:admm} and the results
of the numerical experiments are discussed in Section~\ref{sec:results}. Conclusions are given in Section~\ref{sec:conclusions}.

Throughout this work we denote matrices with uppercase lightface letters, vectors with lowercase boldface
letters and scalars with lowercase lightface letters. All the vectors are column vectors. Given a vector $\bv$, we use $v_i$ or $(\bv)_i$
to denote its $i$-th entry. We use $\R_+$ to indicate the set of real nonnegative numbers and $\| \cdot \|$ to indicate the two-norm.
For brevity, given any vectors $\bv$ and $\bw$ we use the notation $( \bv, \bw )$ instead of $[ \bv^\top \; \bw^\top ]^\top$.
Likewise, given any scalars $v$ and $w$, we use $(v, \, w)$ to indicate the vector $[ v \;\; w ]^\top$. We also use the notation
$([\bv]_1, [\bv]_2)$ to highlight the subvectors $[\bv]_1$ and $[\bv]_2$ forming the vector $\bv$.
Finally, by writing $\bv > 0$ we mean that all the entries of $\bv$ are nonnegative and at least one of them is positive.


\section{The KL-DTGV$^2$ Model\label{sec:kl-dtgv}}

We briefly describe the KL-DTGV$^2$ model for the restoration of directional images corrupted by Poisson noise.
Let {$\bb\in\R^n$} be the observed image. We want to recover the original image by minimizing a combination
of the KL divergence~\eqref{eq:kl} and the DTGV$^2$ regularizer~\eqref{eqn:dtgv2}, i.e., by solving the optimization problem

\begin{equation}\label{eq:kldtgv}
    \begin{array}{ll}
        \min\limits_{\bu,\bw} & \displaystyle \lambda\,D_{KL}(A\bu+\bgamma,\bb)  +
        \alpha_0 \left\|\widetilde{\nabla} \bu - \bw\right\|_{2,1|\R^{2n}} + \alpha_1 \left\|\widetilde{\mE}\bw\right\|_{2,1|\R^{4n}} \\
        \st       & \bu \geq 0,
    \end{array}
\end{equation}

\noindent
where $\bu\in\R^n$, {$A\in\R^{n\times n}$}, $\bgamma,\bb\in{\R^{n}}$, $\bw\in\R^{2n}$, and $\widetilde{\nabla}\in\R^{2n\times n}$
and $\widetilde{\mE}\in\R^{4n\times2n}$ are the linear operators defining the DTGV$^2$ regularization. The parameters $\lambda\in(0,\,+\infty)$ and $\alpha_0,\alpha_1\in(0,\,1)$ determine the balance between the KL data fidelity term and the two components of the regularization term.

We note that problem~\eqref{eq:kldtgv} is a nonsmooth convex optimization problem because of the properties
of the KL divergence (see, e.g., \cite{diserafino:2020amc}) and the DTGV operator (see, e.g., \cite{kongskov:2019}).


\section{Efficient Estimation of the Image Direction\label{sec:direction}}

An essential ingredient in the DTGV regularization is the estimation of the angle $\theta$ representing
the image texture direction. In~\cite{kongskov:2019}, an estimation algorithm based on the one in \cite{setzer:2008}
is proposed, whose basic idea is to compute a pixelwise direction estimate and then $\theta$ as the average
of that estimate. In \cite{kongskov:2017}, which focuses on impulse noise removal, a more efficient and robust algorithm for
estimating the direction is presented, based on the Fourier transform. The main idea behind this algorithm is to
exploit the fact that two-dimensional Fourier basis functions can be seen as images with one-directional patterns.
However, despite being very efficient from a computational viewpoint, this technique does not appear to be fully
reliable in our tests on Poissonian images (see Section~\ref{sec:estimation_results}). Therefore, we propose
a different approach for estimating the direction,
based on classical tools of image processing: the Sobel filter~\cite{kanopoulos1988} and the Hough
transform~\cite{hough:1962, duda:1972}.

\begin{figure}[H]	    
\begin{center}
    \begin{tabular}[b]{cc}
        \includegraphics[width=5cm]{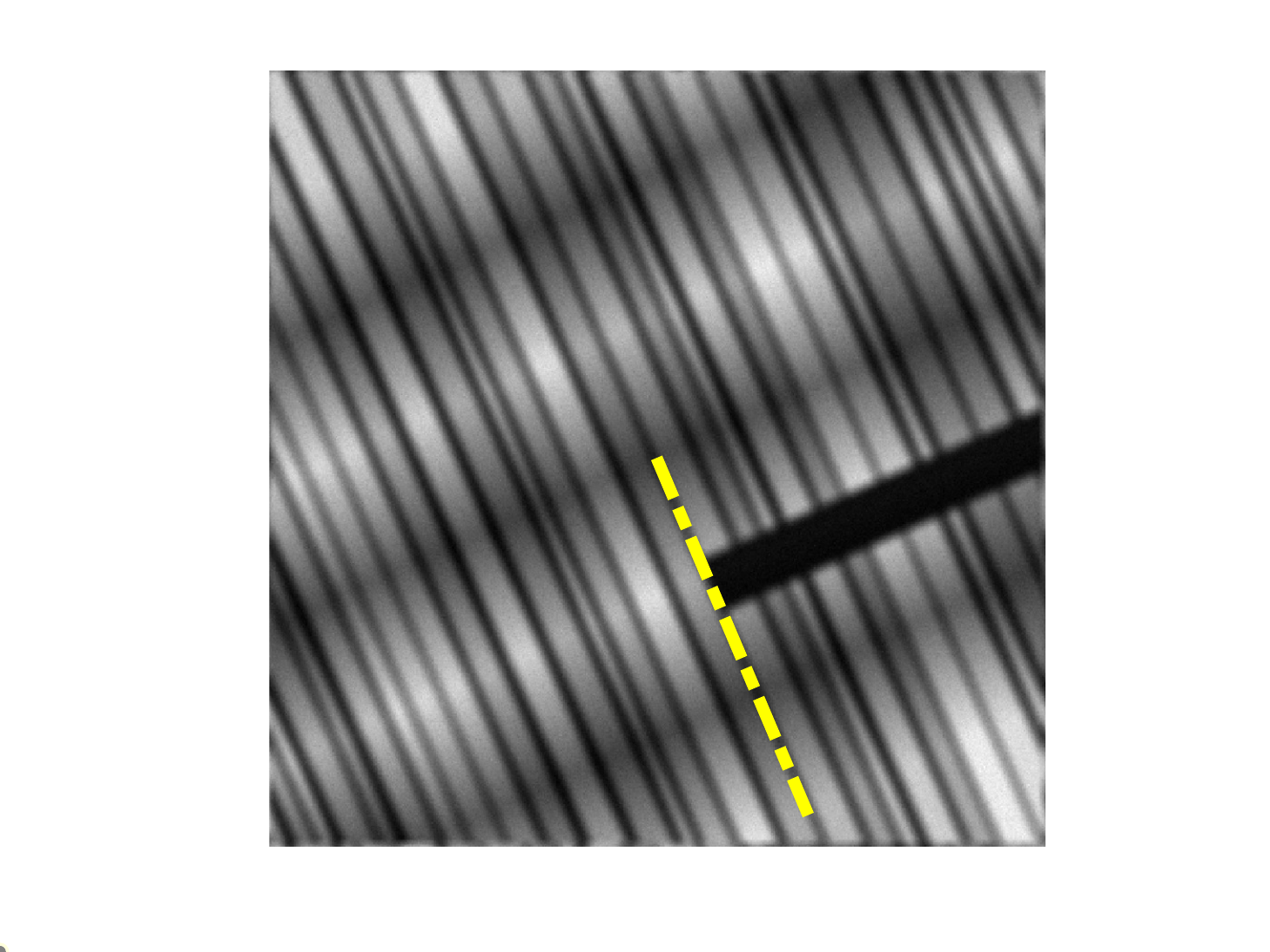} &
        \includegraphics[width=5cm]{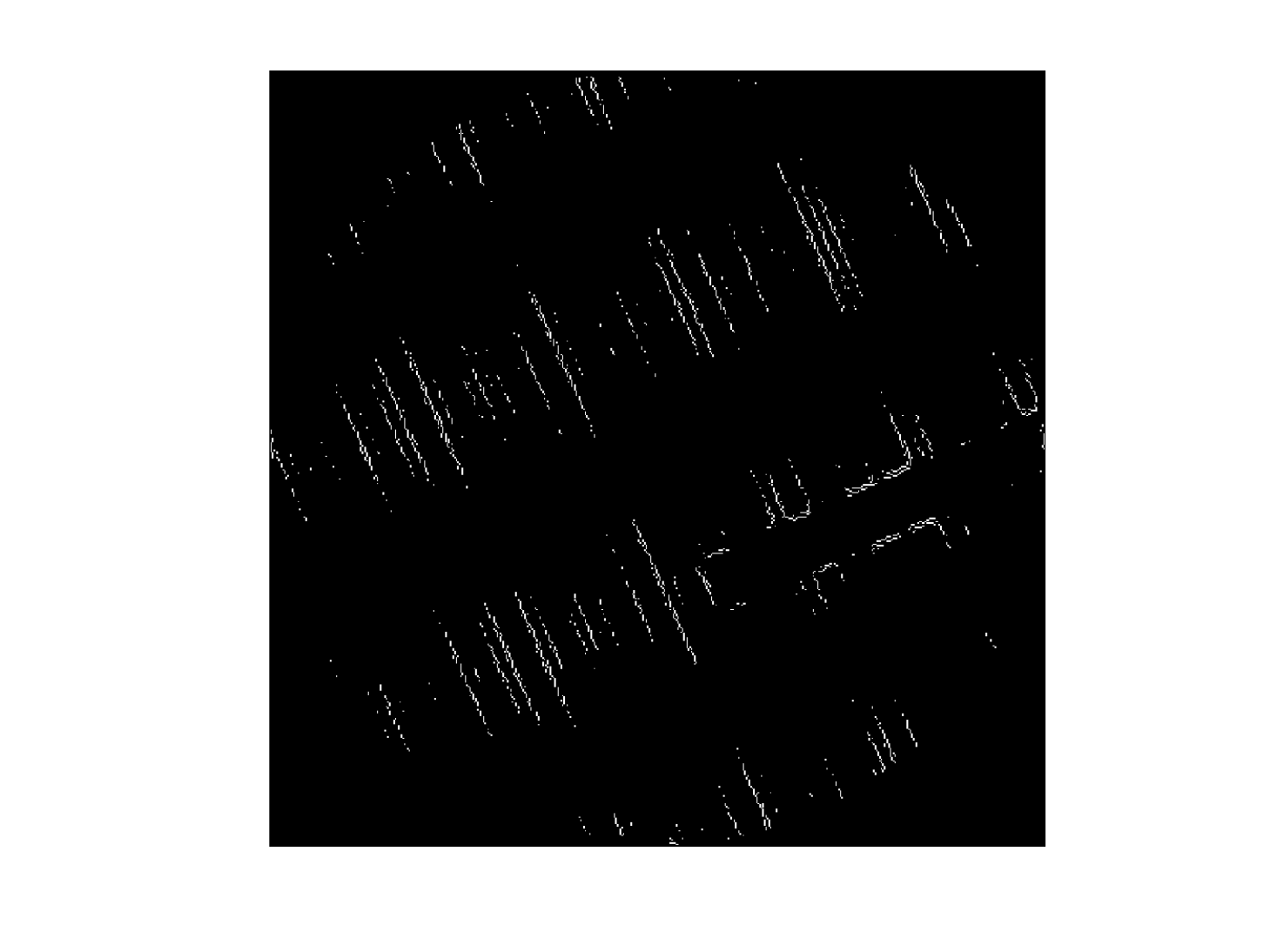} \\[-8pt]
        \footnotesize (\textbf{a}) image &
        \footnotesize (\textbf{b}) edge detection (Sobel filter + mask) \\[5pt]
        \includegraphics[width=5cm]{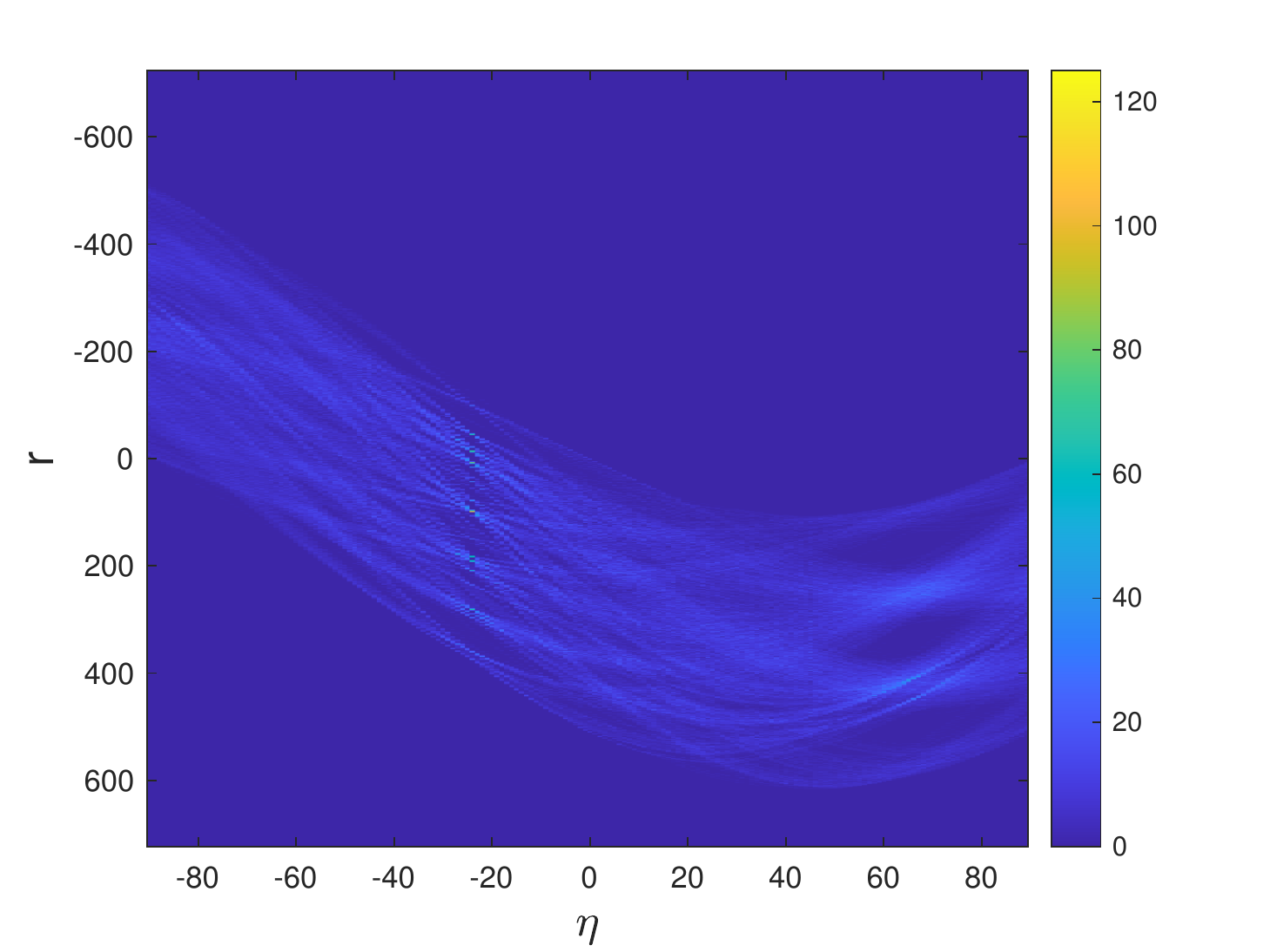} &
        \includegraphics[width=5cm]{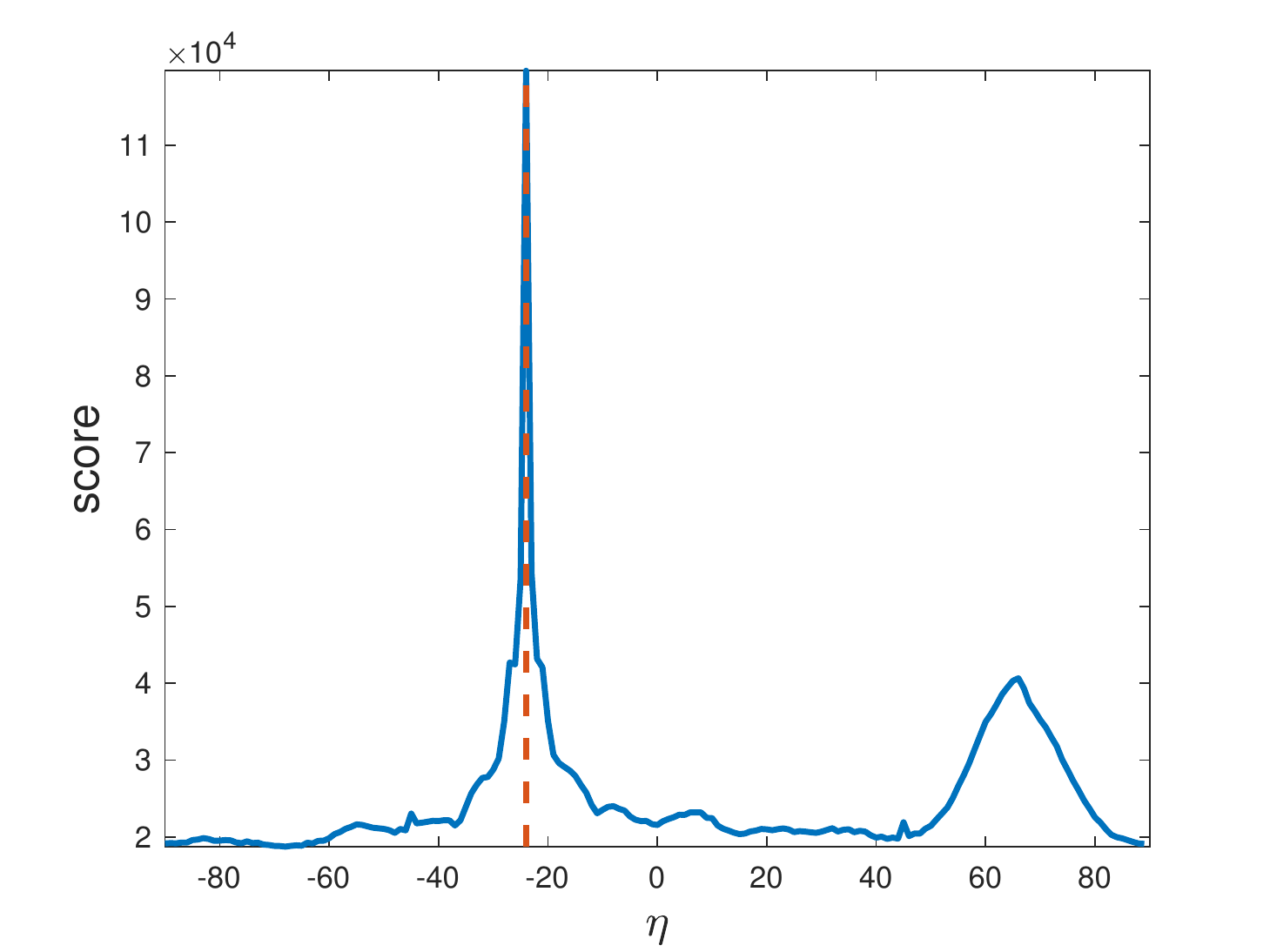} \\[-7pt]
        \footnotesize (\textbf{c}) Hough transform &
        \footnotesize (\textbf{d}) direction scoring
    \end{tabular}
    \caption{Workflow of Algorithm~\ref{alg:direction_estimation} on a random directional image.\label{fig:estimation_example}} 
\end{center}
\end{figure}

Our technique is based on the idea that if an image has a one-directional structure, i.e., its main pattern
consists of stripes, then the edges of the image mainly consist of lines going in the direction of the stripes.
The first stage of the proposed algorithm uses the Sobel filter to determine the edges of the noisy and blurry image.
Then, the Hough transform is applied to the edge image in order to detect the lines.
{The Hough transform is based on the idea that each straight line can be identified by a pair $(r,\,\eta)$ where $r$ is the distance of the line from the origin, and $\eta$ is the angle between the $x$ axis and the segment connecting the origin with its orthogonal projection on the line. The output of the transform is a matrix in which each entry is associated with a pair $(r,\,\eta)$, i.e., with a straight line in the image, and its value is the sum of the values in the pixels that are on the line. Hence, the elements with the highest value in the Hough transform indicate the lines that are most likely to be present in the input image.
    Because of its definition,} the Hough transform tends to overestimate diagonal lines in rectangular images (diagonal lines through the central part of the image
contain the largest number of pixels); therefore, before computing the transform we apply a mask to the edge image, considering only
the pixels inside the largest circle centered in the center of the image. After the Hough transform has been applied,
we compute the square of the two-norm of each column of the matrix resulting from the transform, to determine a score for each angle from $-90^o$ to $90^o$. {Intuitively, the score for each angle is related to the number of lines with that particular inclination which have been detected in the image.
    Finally, we set the direction estimate $\theta \in [-\pi, \pi]$ as
    \begin{equation*}
        \theta = \left\lbrace \begin{array}{ll}
            \displaystyle \frac{90-\eta_{\max}}{180}\pi, & \eta_{\max} \geq 0,\\[8pt]
            \displaystyle \frac{-90-\eta_{\max}}{180}\pi, & \eta_{\max} < 0.
        \end{array}\right.
    \end{equation*}
    where $\eta_{\max}$ is the value of $\eta$ corresponding to the maximum score.} A pseudocode for the estimation algorithm
is provided in Algorithm~\ref{alg:direction_estimation} and an example of the algorithm workflow is given
in Figure~\ref{fig:estimation_example}.

\begin{algorithm}[H]
    \caption{Direction estimation\label{alg:direction_estimation}.}
    {\small
        \begin{algorithmic}[1]
            \STATE Use the Sobel operator to obtain the image $\be$ of the edges of the noisy and blurry image $\bb$.
            \STATE Apply a disk mask to cut out some diagonal edges in $\be$, obtaining a new edge image $\widetilde \be$ (Figure~\ref{fig:estimation_example}b).
            \STATE Compute the Hough transform $h( \widetilde \be )$ (Figure~\ref{fig:estimation_example}c).
            \STATE { Set $\eta_{\max}$ as the value of $\eta$ corresponding to} the column of $h( \widetilde \be )$ with maximum 2-norm. (Figure~\ref{fig:estimation_example}d)
            \STATE { Set $\theta = \left\lbrace \begin{array}{ll}
                    \frac{90-\eta_{\max}}{180}\pi, & \eta_{\max}\geq 0,\\[2pt]
                    \frac{-90-\eta_{\max}}{180}\pi, & \eta_{\max}< 0.
                \end{array}\right.$} (yellow line in Figure~\ref{fig:estimation_example}a)
        \end{algorithmic}
    }
\end{algorithm}

\section{ADMM for Minimizing the KL-DTGV$^2$ Model\label{sec:admm}}

Although problem~\eqref{eq:kldtgv} is a bound-constrained convex optimization problem, the nondifferentiability
of the DTGV$^2$ regularizer does not allow its solution by classical optimization methods for smooth
problems, such as gradient methods (see \cite{bonettini:2009, birgin:2014, diserafino:2018amc} and the
references therein). However, the problem can be solved by methods based on splitting techniques, such as
\cite{goldstein:2009, boyd:2011admm, parik:2014, bonettini:2016, desimone:2020}.
Here we solve~\eqref{eq:kldtgv} by the Alternating Direction Method of Multipliers
(ADMM)~\cite{boyd:2011admm}. To this end, we first reformulate the problem as follows:
\begin{equation}\label{eq:kldtgv_reformulate}
    \begin{array}{ll}
        \min\limits_{\bu,\bw,\bz_1,\bz_2,\bz_3,\bz_4} & \displaystyle \lambda\,D_{KL}(\bz_1+\bgamma,\bb)  +  \alpha_0\,\|\bz_2\|_{2,1|\R^{2n}} + \alpha_1\,\|\bz_3\|_{2,1|\R^{4n}} + \chi_{\R_+^n}(\bz_4) \\
        \st     & \bz_1 = A\,\bu,\\
        & \bz_2 = \widetilde{\nabla}\bu - \bw,\\
        & \bz_3 = \widetilde{\mE}\bw,\\
        & \bz_4 = \bu,
    \end{array}
\end{equation}

\noindent
where {$\bz_1\in\R^n$}, $\bz_2\in\R^{2n}$, $\bz_3\in\R^{4n}$, $\bz_4\in\R^{n}$, and $\chi_{\R_+^n}(\bz_4)$
is the characteristic function of the nonnegative orthant in $\R^n$. {A similar splitting has been used in \cite{setzer:2010} for TV-based deblurring of Poissonian images.} By introducing the auxiliary variables
$\bx = ( \bu, \bw )$ and  $\bz = ( \bz_1, \bz_2, \bz_3, \bz_4 )$ we can
further reformulate the KL-DTGV$^2$ problem as
\begin{equation}\label{eq:kldtgv_reformulate_admm}
    \begin{array}{ll}
        \min\limits_{\bx,\bz} & \displaystyle F_1(\bx) + F_2(\bz) \\
        \st     & H\,\bx + G\,\bz = 0,
    \end{array}
\end{equation}

\noindent
where we set
\begin{equation}\label{eq:define_F1_F2}
    F_1(\bx)= 0,\quad F_2(\bz) = \lambda\,D_{KL}(\bz_1+\bgamma,\bb)  +  \alpha_0\,\|\bz_2\|_{2,1|\R^{2n}} + \alpha_1\,\|\bz_3\|_{2,1|\R^{4n}} + \chi_{\R_+^n}(\bz_4),
\end{equation}
and we define the matrices {$H\in\R^{8n\times 3n}$ and $G\in\R^{8n\times 8n}$} as
\begin{equation}\label{eq:define_H_G}
    H = \left[ \begin{array}{rr} A & 0 \\ \widetilde{\nabla} & -I_{2n} \\ 0 & \widetilde{\mE} \\ I_n & 0 \end{array} \right], \quad
    G = \left[ \begin{array}{rrrr} {-I_n} & 0 & 0 & 0 \\ 0 & -I_{2n} & 0 & 0 \\ 0 & 0 & -I_{4n} & 0 \\ 0 & 0 & 0 & -I_n\end{array} \right].
\end{equation}

We consider the Lagrangian function associated with problem \eqref{eq:kldtgv_reformulate_admm},

\begin{equation}\label{eq:kldtgv_lagrangian}
    \mL(\bx,\bz,\bxi) = F_1(\bx) + F_2(\bz) + \bxi^\top\left(H\,\bx + G\,\bz\right),
\end{equation}

\noindent
where {$\bxi\in\R^{8n}$} is a vector of Lagrange multipliers, and then the augmented \linebreak \mbox{Lagrangian function }

\begin{equation}\label{eq:kldtgv_augmented_lagrangian}
    \mL_A(\bx,\bz,\bxi;\rho) = F_1(\bx) + F_2(\bz) + \bxi^\top\left(H\,\bx + G\,\bz\right) + \frac{\rho}{2}\left\|H\,\bx + G\,\bz\right\|_2^2 ,
\end{equation}

\noindent
where $\rho>0$.

Now we are ready to introduce the ADMM method for the solution of problem~\eqref{eq:kldtgv_reformulate_admm}.
Let $\bx^0\in\R^{3n}$, {$\bz^0\in\R^{8n}$, $\bxi^0\in\R^{8n}$}. 
At each step $k>0$ the ADMM method computes the new iterate $\left ( \bx^{k+1},\bz^{k+1},\bxi^{k+1} \right )$ as
follows:
\begin{equation}\label{eq:admm_method}
    \begin{split}
        \bx^{k+1} & = \displaystyle \arg\min\limits_{\bx\in\R^{3n}} \mL_A(\bx,\bz^k,\bxi^k;\rho),\\
        \bz^{k+1} & = \displaystyle \arg\min\limits_{\bz\in\R^{8n}} \mL_A(\bx^{k+1},\bz,\bxi^k;\rho),\\
        \bxi^{k+1} & = \displaystyle \bxi^k + \rho\left(H\,\bx^{k+1} + G\,\bz^{k+1}\right).
    \end{split}
\end{equation}

Note that the functions $F_1(\bx)$ and $F_2(\bz)$ in~\eqref{eq:kldtgv_reformulate_admm}
are closed, proper and convex. Moreover, the matrices $H$ and $G$ defined in~\eqref{eq:define_H_G}
are such that {$G=-I_{8n}$} and $H$ has full rank. Hence, the convergence of the method defined by~\eqref{eq:admm_method}
can be proved by applying a classical convergence result from the seminal paper by Eckstein and
Bertsekas~\cite[Theorem~8]{eckstein:1992}, which we report in a form that can be immediately applied to our
reformulation of \mbox{the problem}. 

\begin{theorem}
    Let us consider a problem of the form~\eqref{eq:kldtgv_reformulate_admm} where $F_1(\bx)$ and $F_2(\bz)$ are closed, proper
    and convex functions and $H$ has full rank. Let $\bx^0\in\R^{3n}$, $\bz^0\in{\R^{8n}}$, $\bxi^0\in{\R^{8n}}$, and $\rho>0$.
    Suppose $\{\varepsilon_k\}, \{\nu_k\} \subset \R_+$ are summable sequences such that for all $k$
    \begin{eqnarray*}
        &&\left\| \bx^{k+1} - \arg\min\limits_{\bx\in\R^{3n}} \mL_A(\bx,\bz^k,\bxi^k;\rho)\right\| \leq \varepsilon_k,\\
        &&\left\| \bz^{k+1} - \arg\min\limits_{\bz\in{\R^{8n}}} \mL_A(\bx^{k+1},\bz,\bxi^k;\rho)\right\| \leq \nu_k,\\
        &&\bxi^{k+1} = \bxi^k + \rho\left(H\,\bx^{k+1} + G\,\bz^{k+1}\right).
    \end{eqnarray*}
    
    \noindent If there exists a saddle point $(\bx^*,\bz^*,\bxi^*)$ of $ \mL(\bx,\bz,\bxi)$, then $\bx^k\rightarrow\bx^*$,
    $\bz^k\rightarrow\bz^*$ and $\bxi^k\rightarrow\bxi^*$. If such saddle point does not exist, then at least
    one of the sequences $\{\bz^k\}$ or $\{\bxi^k\}$ is unbounded.
\end{theorem}

Since we are dealing with linear constraints, we can recast \eqref{eq:admm_method} in a more convenient
form, by observing that the linear term in \eqref{eq:kldtgv_augmented_lagrangian} can be included in the
quadratic one. By introducing the vector of scaled Lagrange multipliers $ \bmu^k = \frac{1}{\rho}\bxi^k$,
the ADMM \mbox{method becomes}
\begin{eqnarray}
    \bx^{k+1} &=& \displaystyle \arg\min\limits_{\bx\in\R^{3n}} \frac{\rho}{2}\left\|H\,\bx - \bz^k + \bmu^k\right\|_2^2,\label{eq:admm_method_mod_X_problem}\\
    \bz^{k+1} &=& \displaystyle \arg\min\limits_{\bz\in{\R^{8n}}} F_2(\bz) + \frac{\rho}{2}\left\|H\,\bx^{k+1} - \bz + \bmu^k\right\|_2^2,\label{eq:admm_method_mod_Z_problem}\\
    \bmu^{k+1} &=& \displaystyle \bmu^k + H\,\bx^{k+1} + G\,\bz^{k+1}.\label{eq:admm_method_mod_mu_update}
\end{eqnarray}

In the next sections we show how the solutions to subproblems \eqref{eq:admm_method_mod_X_problem}
and \eqref{eq:admm_method_mod_Z_problem} can be computed exactly with a small computational effort.

\subsection{Solving the Subproblem in $\bx$ \label{sec:X_subproblem}}

Problem \eqref{eq:admm_method_mod_X_problem} is an overdetermined least squares problem,
since $H$ is a tall-and-skinny matrix with full rank. Hence, its solution can be computed by solving
the normal \mbox{equations system}

\begin{equation}\label{eq:X_subp_norm_eq}
    H^\top H\,\bx = H^\top\bv_x^k,
\end{equation}

\noindent
where we set $\bv_x^k = \bz^k - \bmu^k$. Starting from the definition of $H$ given in \eqref{eq:define_H_G},
we have
\begin{eqnarray*}
    H^\top H &=& \left[\begin{array}{c|c} I_n + A^\top A + \widetilde{\nabla}^\top\widetilde{\nabla}  & -\widetilde{\nabla}^\top  \\[1pt] \hline\\[-7pt]
        -\widetilde{\nabla} & I_{2n} + \widetilde{\mE}^\top\widetilde{\mE} \end{array} \right] = \\
    & =& \left[\begin{array}{c|cc} I_n + A^\top A + \widetilde{\nabla}^\top\widetilde{\nabla}  & -D_{\theta}^\top  & -D_{\theta^\bot}^\top \\[1pt] \hline\\[-7pt]
        -D_{\theta}      & I_{n} + D_{\theta}^\top D_{\theta} + \frac{1}{2}D_{\theta^\bot}^\top D_{\theta^\bot} & \frac{1}{2}D_{\theta^\bot}^\top D_{\theta}\\[5pt]
        -D_{\theta^\bot} & \frac{1}{2}D_{\theta}^\top D_{\theta^\bot} & I_{n} + \frac{1}{2}D_{\theta}^\top D_{\theta} + D_{\theta^\bot}^\top D_{\theta^\bot}
    \end{array} \right].
\end{eqnarray*}

System \eqref{eq:X_subp_norm_eq} may be quite large and expensive, also for relatively small images.
However, as pointed out in Assumption~\ref{ass:bccb_structure}, 
$A$, $D_\theta$ and $D_{\theta^\bot}$ have a BCCB structure, hence all the blocks of $H^\top H$ maintain that structure.
By recalling that BCCB matrices can be diagonalized by means of two-dimensional Discrete Fourier Transforms
(DFTs), we show how the solution to \eqref{eq:X_subp_norm_eq} can be computed expeditiously.

Let $\mF\in\C^{n\times n}$ be the matrix representing the two-dimensional DFT operator, and let $\mF^*$
denote its inverse, i.e., its adjoint. We can write $H^\top H$ as

\begin{equation}\label{eq:HtH_FFT}
    H^\top H = \left[\begin{array}{ccc}
        \mF^* & 0     & 0\\
        0     & \mF^* & 0  \\
        0     & 0     & \mF^*
    \end{array} \right] 
    \left[\begin{array}{c|cc} \Gamma & -\Delta_{\theta}^*  & -\Delta_{\theta^\bot}^* \\\hline
        -\Delta_{\theta}      & \Phi_{11} & \Phi_{12}\\
        -\Delta_{\theta^\bot} & \Phi_{21} & \Phi_{22}
    \end{array} \right] 
    \left[\begin{array}{ccc}
        \mF & 0     & 0\\
        0     & \mF & 0  \\
        0     & 0     & \mF
    \end{array} \right],
\end{equation}

\noindent
where each block of the central matrix is the diagonal complex matrix 
associated with the corresponding block in $H^\top H$, and $\Delta_\theta^*, \Delta_{\theta^\bot}^*$
denote the (diagonal) adjoint matrices of $\Delta_\theta, \Delta_{\theta^\bot}$.
By~\eqref{eq:HtH_FFT} and the definition of $\bx$, we can reformulate \eqref{eq:X_subp_norm_eq} as

\begin{equation}\label{eq:x_subp_FFTs}
    \left[\begin{array}{c|cc} \Gamma & -\Delta_{\theta}^*  & -\Delta_{\theta^\bot}^* \\ \hline
        -\Delta_{\theta}      & \Phi_{11} & \Phi_{12} \\
        -\Delta_{\theta^\bot} & \Phi_{21} & \Phi_{22}
    \end{array} \right] 
    \left[\begin{array}{c}
        \mF\,\bu\\ \mF\,\bw_1\\ \mF\,\bw_2
    \end{array} \right] =
    \left[\begin{array}{c}
        \mF\,[H^\top\bv_x^k]_1\\ \mF\,[H^\top\bv_x^k]_2\\ \mF\,[H^\top\bv_x^k]_3
    \end{array} \right],
\end{equation}

\noindent
where we split $\bw$ and $\bv_x^k$ in two and three blocks of size $n$, respectively.

Now we recall a result about the inversion of block matrices. Suppose that a square matrix $M$
is partitioned into four blocks, i.e.,

$$
M = \left[ \begin{array}{cc}
    M_{11} & M_{12} \\
    M_{21} & M_{22}
\end{array} \right];
$$

\noindent
then, if $M_{11}$ and $M_{22}$ are invertible, we have

\begin{equation}\label{eq:block_matrix_inversion}
\small
\begin{split}
    M^{-1} = & \left[\begin{array}{cc}M_{11} &M_{12} \\M_{21} &M_{22} \end{array}\right]^{-1} \\
    = & \left[\begin{array}{cc}\left(M_{11} -M_{12} M_{22} ^{-1}M_{21} \right)^{-1}&0 \\ 0 &\left(M_{22} -M_{21} M_{11} ^{-1}M_{12} \right)^{-1}\end{array}\right]
    \left[\begin{array}{cc}I &-M_{12} M_{22} ^{-1}\\-M_{21} M_{11} ^{-1}&I \end{array}\right].
\end{split}
\end{equation}

By applying \eqref{eq:block_matrix_inversion} to the matrix consisting
of the second and third block rows and columns of the matrix in \eqref{eq:x_subp_FFTs},
which we denote $\Phi$, we get
\begin{equation}\label{eq:HTH_23_block_inversion}
    \small
    \begin{split}
        \Phi^{-1} = & \left[\begin{array}{cc}
            \Phi_{11} & \Phi_{12} \\
            \Phi_{21} & \Phi_{22}
        \end{array}\right]^{-1} \\
        = &
        \left[\begin{array}{cc}
            \left(\Phi_{11} -\Phi_{12} \Phi_{22} ^{-1}\Phi_{21} \right)^{-1} & -\left(\Phi_{11} -\Phi_{12} \Phi_{22} ^{-1}\Phi_{21} \right)^{-1} \Phi_{12} \Phi_{22}^{-1} \\
            -\left(\Phi_{22} -\Phi_{21} \Phi_{11} ^{-1}\Phi_{12} \right)^{-1}\Phi_{21} \Phi_{11} ^{-1} & \left(\Phi_{22} -\Phi_{21} \Phi_{11} ^{-1}\Phi_{12} \right)^{-1}
        \end{array}\right].
    \end{split}
\end{equation}

To simplify the notation we set
\begin{equation} \label{eq:HTH_23_block_inversion_psi}
    \Psi =
    \left[\begin{array}{cc}
        \Psi_{11} & \Psi_{12} \\
        \Psi_{21} & \Psi_{22}
    \end{array}\right]
    =
    \Phi^{-1},
\end{equation}

\noindent
and observe that the matrices $\Psi_{ij} \in \C^{n \times n}$ are diagonal.
Letting $\Delta^* = \left[ \Delta_{\theta}^* \;\;\; \Delta_{\theta^\bot}^* \right]$,
applying the inversion formula \eqref{eq:block_matrix_inversion} to the whole matrix in~\eqref{eq:x_subp_FFTs}, and
using \eqref{eq:HTH_23_block_inversion} and \eqref{eq:HTH_23_block_inversion_psi}, we get

\begin{equation}\label{eq:x_subp_FFTs_inversion}
    \left[\begin{array}{cc} \Gamma & -\Delta^*\\
        -\Delta      & \Phi
    \end{array} \right]^{-1} = \left[\begin{array}{cc} \Xi^{-1} & 0 \\
        0 & \Omega^{-1}\end{array}\right]   
    \left[\begin{array}{cc}I_n & -\Delta^*\Psi \\ -\Delta\, \Gamma ^{-1} & I_{2n} \end{array}\right],
\end{equation}

\noindent
where
\begin{equation*}
    \begin{split}
        \Xi &=  \Gamma - \Delta^* \Psi \,\Delta = \Gamma -\left[ \Delta_{\theta}^* \;\;\; \Delta_{\theta^\bot}^* \right] \left[\begin{array}{cc}
            \Psi_{11} & \Psi_{12} \\
            \Psi_{21} & \Psi_{22}
        \end{array}\right] \left[\begin{array}{c}
            \Delta_{\theta} \\
            \Delta_{\theta^*}
        \end{array}\right],\\
        \Omega &= \Phi - \Delta \,\Gamma^{-1}\Delta^* = \left[\begin{array}{cc}
            \Phi_{11} & \Phi_{12} \\
            \Phi_{21} & \Phi_{22}
        \end{array}\right] -\left[\begin{array}{c}
            \Delta_{\theta} \\
            \Delta_{\theta^\bot}
        \end{array}\right] \Gamma ^{-1}\left[ \Delta_{\theta}^* \;\;\; \Delta_{\theta^\bot}^* \right].
    \end{split}
\end{equation*}

We note that $\Xi\in\C^{n \times n}$ is diagonal (and its inversion is straightforward), while $\Omega\in\C^{2n \times 2n}$ has a
$2\times 2$ block structure with blocks that are diagonal matrices belonging to $\C^{n \times n}$. Thus, we can compute $\Upsilon = \Omega^{-1}$ by applying \eqref{eq:block_matrix_inversion}:

\begin{equation}\label{eq:Omega_inversion}
    \begin{split}
        \Upsilon & = \left[ \begin{array}{cc}
            \Upsilon_{11} & \Upsilon_{12} \\
            \Upsilon_{21} & \Upsilon_{22}
        \end{array} \right] =
        \left[\begin{array}{cc}
            \Omega_{11} & \Omega_{12} \\
            \Omega_{21} & \Omega_{22}
        \end{array}\right]^{-1} \\
        & = \left[ \begin{array}{cc}
            \left(\Omega_{11} -\Omega_{12} \Omega_{22} ^{-1}\Omega_{21} \right)^{-1} & -\left(\Omega_{11} -\Omega_{12} \Omega_{22} ^{-1}\Omega_{21} \right)^{-1} \Omega_{12} \Omega_{22}^{-1} \\
            -\left(\Omega_{22} -\Omega_{21} \Omega_{11} ^{-1}\Omega_{12} \right)^{-1}\Omega_{21} \Omega_{11} ^{-1} & \left(\Omega_{22} -\Omega_{21} \Omega_{11} ^{-1}\Omega_{12} \right)^{-1}
        \end{array} \right]
    \end{split}
\end{equation}

Summing up, by \eqref{eq:x_subp_FFTs}, \eqref{eq:x_subp_FFTs_inversion}
and \eqref{eq:Omega_inversion}, the solution to \eqref{eq:X_subp_norm_eq} can be
obtained by computing

\begin{equation}\label{eq:x_subp_compute_y}
    \left[\begin{array}{c}
        \by_1\\ \by_2\\ \by_3
    \end{array} \right] = \left[\begin{array}{cc}
        \Xi^{-1} & 0  \\
        0 & \Upsilon
    \end{array} \right] \left[\begin{array}{cc}I_n & -\Delta^\top\Psi \\-\Delta\, \Gamma ^{-1} & I_{2n} \end{array}\right]
    \left[\begin{array}{c}
        \mF\,[H^\top\bv_x^k]_1\\ \mF\,[H^\top\bv_x^k]_2\\ \mF\,[H^\top\bv_x^k]_3
    \end{array} \right],
\end{equation}

\noindent
and setting

\begin{equation}\label{eq:x_subp_update_u_w}
    \bu^{k+1} = \mF^*\by_1,\qquad \bw^{k+1}_1 = \mF^*\by_2, \qquad \bw^{k+1}_2 = \mF^*\by_3.
\end{equation}

\begin{Remark}
    The only quantity in \eqref{eq:x_subp_compute_y} that varies at each iteration is $\bv_x^k$.
    Hence, the matrices $\Delta$, $\Gamma$, $\Psi$, $\Xi^{-1}$, and $\Upsilon$ can be
    computed only once before the ADMM method starts. This means that the overall cost of the
    exact solution of \eqref{eq:admm_method_mod_X_problem} at each iteration reduces
    to six two-dimensional DFTs and two matrix--vector products involving two $3\times 3$ block
    matrices with diagonal blocks of dimension~$n$.
\end{Remark}

\subsection{Solving the Subproblem in $\bz$}

By looking at the form of $F_2(\bz)$--see \eqref{eq:define_F1_F2}--and by defining the vector
$\bv_z^k = H\,\bx^{k+1} + \bmu^k$, we see that problem \eqref{eq:admm_method_mod_Z_problem}
can be split into the four problems
\begin{eqnarray}
    \bz_1^{k+1} &=& \displaystyle \arg\min\limits_{\bz_1\in{\R^{n}}} \lambda D_{KL}(\bz_1+\bgamma,\bb) + \frac{\rho}{2}\left\|\bz_1 - [\bv_z^k]_1\right\|_2^2,\label{eq:z1_subproblem}\\
    \bz_2^{k+1} &=& \displaystyle \arg\min\limits_{\bz_2\in\R^{2n}} \alpha_0 \|\bz_2\|_{2,1|\R^{2n}} + \frac{\rho}{2}\left\|\bz_2 - [\bv_z^k]_2\right\|_2^2,\label{eq:z2_subproblem}\\
    \bz_3^{k+1} &=& \displaystyle \arg\min\limits_{\bz_3\in\R^{4n}} \alpha_1 \|\bz_3\|_{2,1|\R^{4n}} + \frac{\rho}{2}\left\|\bz_3 - [\bv_z^k]_3\right\|_2^2,\label{eq:z3_subproblem}\\
    \bz_4^{k+1} &=& \displaystyle \arg\min\limits_{\bz_4\in\R^{n}} \chi_{\R_+^n}(\bz_4) + \frac{\rho}{2}\left\|\bz_4 - [\bv_z^k]_4\right\|_2^2,\label{eq:z4_subproblem}
\end{eqnarray}

\noindent
where $\bv_z^k = ( [\bv_z^k]_1, [\bv_z^k]_2, [\bv_z^k]_3, [\bv_z^k]_4 )$,
with $[\bv_z^k]_1 \in {\R^{n}}$, $[\bv_z^k]_2 \in \R^{2n}$,$[\bv_z^k]_3 \in \R^{4n}$, and \mbox{$[\bv_z^k]_4 \in \R^n$.}
Now we focus on the solution of the four subproblems.

\subsubsection{Update of $\bz_1$}

By the form of the Kullback--Leibler divergence in \eqref{eq:kl}, the minimization problem~\eqref{eq:z1_subproblem} is equivalent to

\begin{equation}\label{eq:z1_subproblem_extended}
    \min\limits_{\bz_1\in{\R^{n}}} \lambda \sum_{i=1}^{n} \left( b_i \ln \frac{b_i}{(\bz_1)_i+\gamma_i}
    + (\bz_1)_i + \gamma_i - b_i\right) + \frac{\rho}{2}\sum_{i=1}^{n} \left( (\bz_1)_i - d_i \right)^2,
\end{equation}

\noindent
where we set $\bd = [\bv_z^k]_1$ to ease the notation.
From \eqref{eq:z1_subproblem_extended} it is clear that the problem in $\bz_1$ can be split into {$n$} problems of the form

\begin{equation}\label{eq:z1_1D_subproblem}
    \min\limits_{z\in\R} \lambda \left(-b \ln(z+\gamma) + z\right) + \frac{\rho}{2}\left(z - d\right)^2.
\end{equation}

Since the objective function of this problem is strictly convex, its solution can be determined
by setting the gradient equal to zero, i.e., by solving

$$
\lambda \left( -\frac{b}{z+\gamma} + 1\right) + \rho\left(z - d\right) = 0,
$$

\noindent
which leads to the quadratic equation

\begin{equation}\label{eq:z1_subproblem_2nd_order_eq}
    z^2 + \left(\frac{\lambda}{\rho} + \gamma - d\right)z - \frac{\lambda}{\rho} \left(\frac{\rho}{\lambda}\gamma d - \gamma + b\right) = 0.
\end{equation}

{Since, by looking at the domain of the objective function in \eqref{eq:z1_1D_subproblem}, $z+\gamma$ has to be strictly positive, we set each entry of $\bz_1^{k+1}$ as the largest solution of the corresponding quadratic
    Equation \eqref{eq:z1_subproblem_2nd_order_eq}.}

\subsubsection{Update of $\bz_2$ and $\bz_3$}

{The minimization problems \eqref{eq:z2_subproblem} and \eqref{eq:z3_subproblem} correspond
    to the computation of the proximal operators of the functions $f(\bz_2) = \frac{\alpha_0}{\rho}\|\bz_2\|_{2,1|\R^{2n}}$
    and $g(\bz_3) = \frac{\alpha_1}{\rho}\|\bz_3\|_{2,1|\R^{4n}}$, respectively.}

By the definitions given in~\eqref{eq:def_21norm_z2} and~\eqref{eq:def_21norm_z3}, we see that
the two (2,1)-norms correspond to the sum of two-norms of vectors in $\R^2$ and $\R^4$, respectively.
This means that the computation of both the proximal operators can be split into the computation of
$n$ proximal operators of functions that are scaled two-norms in either $\R^2$ or $\R^4$.

The proximal operator of the function $f(\by) = c\|\by\|$, $c > 0$, at a vector $\bd$ is

$$
\prox_{c\|\cdot\|}(\bd) = \arg\min\limits_{\by} c\|\by\| + \frac{1}{2}\|\by-\bd\|^2.
$$

It can be shown (see, e.g., \cite{beck:2017book} [Chapter~6]) that

\begin{equation}\label{eq:proxi_euclidean_norm}
    \prox_{c\|\cdot\|}(\bd) = \left(1-\frac{c}{\max\left\lbrace\|\bd\|,\,c \right\rbrace} \right)\bd = \max\left\lbrace\frac{\|\bd\|-c}{\|\bd\|},\,0 \right\rbrace\bd.
\end{equation}

Hence, for the update of $\bz_2$ we proceed as follows. By setting $\bd = [\bv_z^k]_2$ and
$c=\frac{\alpha_0}{\rho}$, for each $i = 1,\ldots,n $ we have
$$
\left( (\bz_2^{k+1})_i,\, (\bz_2^{k+1})_{n+i} \right) = \prox_{c\|\cdot\|}\left( (d_i, \, d_{n+i}) \right).
$$

To update $\bz_3$, we set $\bd = [\bv_z^k]_3$ and $c=\frac{\alpha_1}{\rho}$ and compute

$$
\left( (\bz_3^{k+1})_i, \, (\bz_3^{k+1})_{n+i}, \, (\bz_3^{k+1})_{2n+i}, \, (\bz_3^{k+1})_{3n+i} \right)
= \prox_{c\|\cdot\|}\left( (d_i, \, d_{n+i}, \, d_{2n+i}, \, d_{3n+i}) \right).
$$

\subsubsection{Update of $\bz_4$}

It is straightforward to verify that the update of $\bz_4$ in~\eqref{eq:z4_subproblem} can be obtained as

$$
\bz_4^{k+1} = \Pi_{\R_+^n}\left( [\bv_z^k]_4\right),
$$

\noindent
where $\Pi_{\R_+^n}$ is the Euclidean projection onto the nonnegative orthant in $\R^n$.

\subsection{Summary of the ADMM Method}

For the sake of clarity, in Algorithm~\ref{alg:admm_kldtgv} we sketch the ADMM version for solving problem~\eqref{eq:kldtgv_reformulate}.

In many image restoration applications, a reasonably good starting guess for $\bu$
is often available. For example, if $A$ represents a blur operator,
a common choice is to set $\bu^0$ equal to the the noisy and blurry image. We make
this choice for $\bu^0$.
By numerical experiments we also verified that once $\bx = (\bu, \bw)$ has been initialized,
it is convenient to set $\bu^1=\bu^0$, $\bw_1^1=\bw_1^0$ and $\bw_2^1=\bw_2^0$ and to shift
the order of the updates in the ADMM scheme \eqref{eq:admm_method_mod_X_problem}--\eqref{eq:admm_method_mod_mu_update},
so that a ``more effective'' initialization of $\bz$ and $\bmu$ is
performed. We see from line~\ref{alg:admm_kldtgv_stopping_criterion} of Algorithm~\ref{alg:admm_kldtgv}
that the algorithm stops when the relative change in the restored image $\bu$ goes below
a certain threshold $tol \in(0,\,1)$ or a maximum number of iterations $k_{\max}$ is reached.
Finally, we note that for the case of the KL-TGV$^2$ model, corresponding to $\theta= 0$ and $a=1$, we have that  $D_\theta = D_H$ and $D_{\theta^\bot} = D_V$; hence, we use the initialization  $\bw_1^0 = D_H\bu^0$ and $\bw_2^0 = D_V\bu^0$.

\begin{algorithm}[H]
    \caption{ADMM for problem \eqref{eq:kldtgv_reformulate}.\label{alg:admm_kldtgv}}
    {\small
        \begin{algorithmic}[1]
            \STATE Let $\bu^0 \in\R^n$, $\bw_1^0 = D_{\theta}\bu^0$, $\bw_2^0 = D_{\theta^\bot}\bu^0$, $\bmu^0 = 0$, $\bz^0 = 0$,  $\lambda,\rho\in(0,\,+\infty)$, $\alpha_0,\alpha_1\in(0,\,1)$
            \STATE Compute matrices $\Delta$, $\Gamma$, $\Psi$, $\Xi^{-1}$, and $\Upsilon$ as specified in Section~\ref{sec:X_subproblem}
            \STATE Let $k=0$, $\bu^1=\bu^0$, $\bw^1=\bw^0$, $stop = false$, $tol\in(0,\,1)$, $k_{\max}\in\N$
            \WHILE{$\mbox{not} \; stop \; \mbox{and} \; k \leq k_{\max}$}
            \STATE Compute $\bz^{k+1}$ by solving the four subproblems \eqref{eq:z1_subproblem}--\eqref{eq:z4_subproblem}
            \STATE Compute $\bmu^{k+1}$ as in \eqref{eq:admm_method_mod_mu_update}
            \STATE $k = k+1$
            \STATE Compute $\bu^{k+1}$, $\bw_1^{k+1}$ and $\bw_2^{k+1}$ by \eqref{eq:x_subp_compute_y} and \eqref{eq:x_subp_update_u_w}
            \STATE Set $stop = \left(\|\bu^{k+1}-\bu^{k}\|< tol\,\|\bu^{k}\|\right)$ \label{alg:admm_kldtgv_stopping_criterion}
            \ENDWHILE
        \end{algorithmic}
    }
\end{algorithm}


\section{Numerical Results\label{sec:results}}
All the experiments were carried out using MATLAB R2018a on a 3.50 GHz Intel Xeon E3
with 16 GB of RAM and Windows operating system. In this section, we first illustrate
the effectiveness of Algorithm~\ref{alg:direction_estimation} for the estimation of the image direction
by comparing it with the one given in~\cite{kongskov:2017} and by analysing its sensitivity
to the degradation in the image to be restored. Then, we present numerical experiments that
demonstrate the improvement of the KL-DTGV$^2$ model upon the KL-TGV$^2$ model
for the restoration of directional images corrupted by Poisson noise.

\begin{figure}[htbp]
\begin{subfigure}{0.33\textwidth}
    \begin{center}
        \includegraphics[width=5cm]{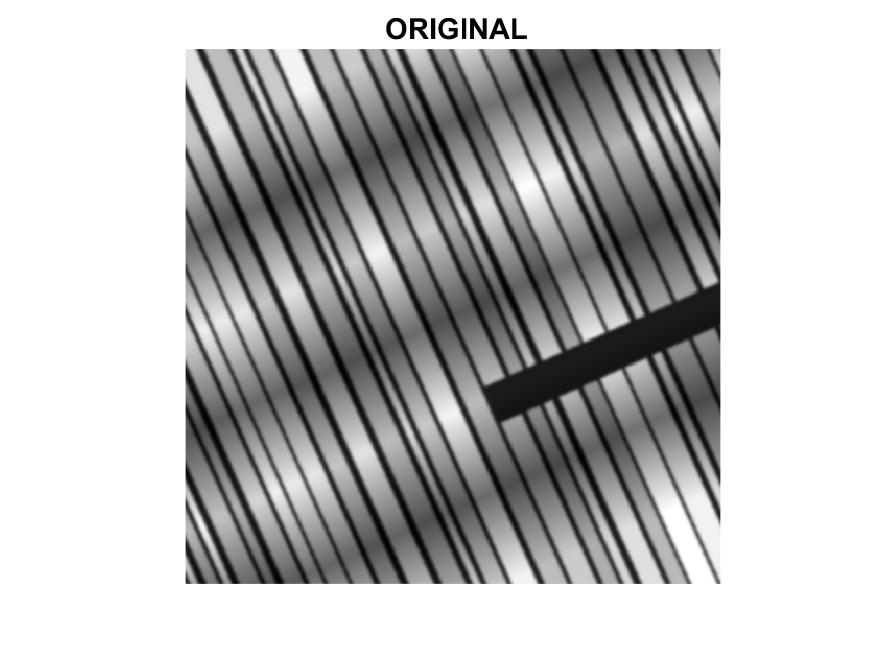}
    \end{center}
\end{subfigure}
\begin{subfigure}{0.66\textwidth}
    \begin{center}
        \includegraphics[width=5cm]{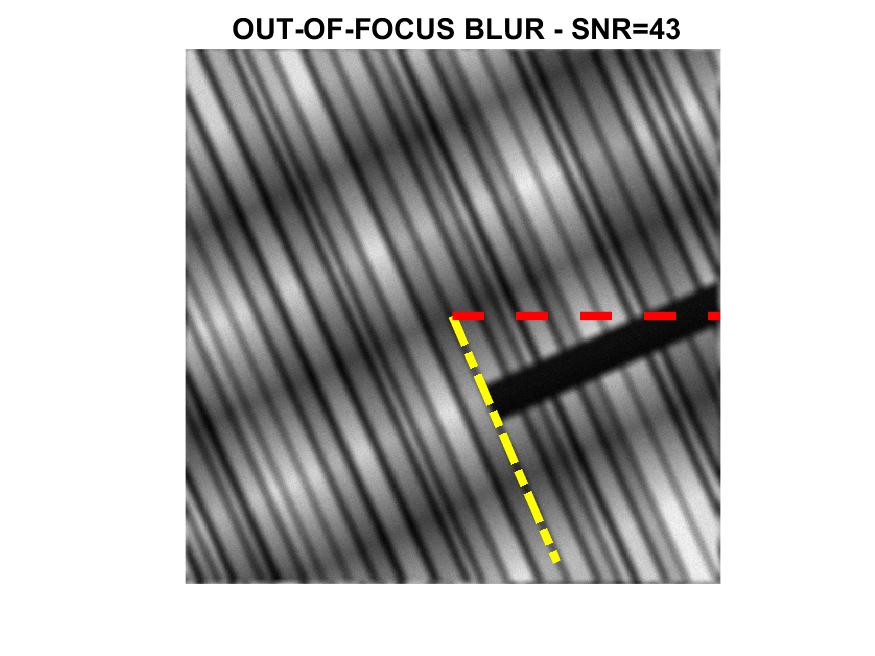}
        \includegraphics[width=5cm]{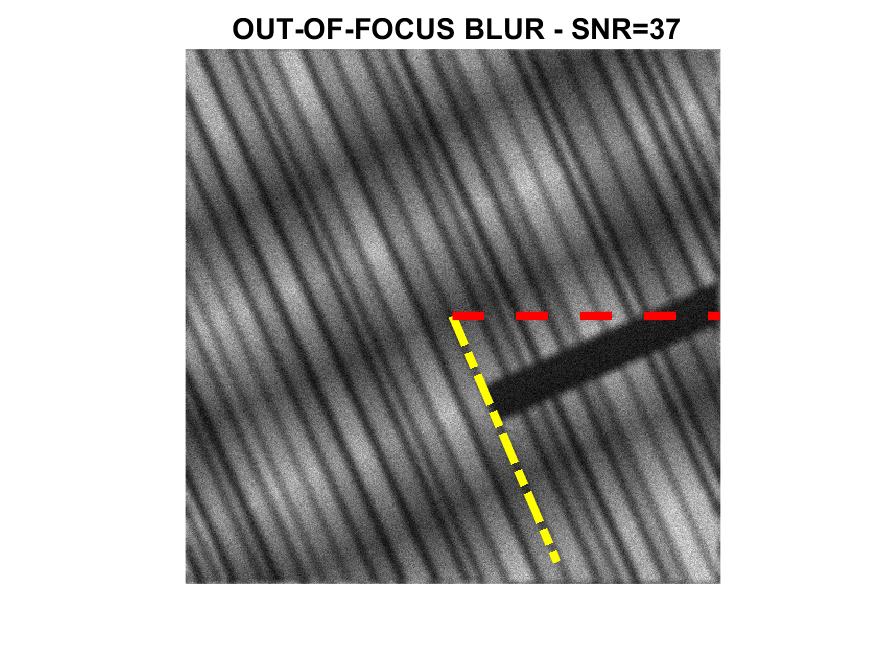} \\[-3mm]
        \includegraphics[width=5cm]{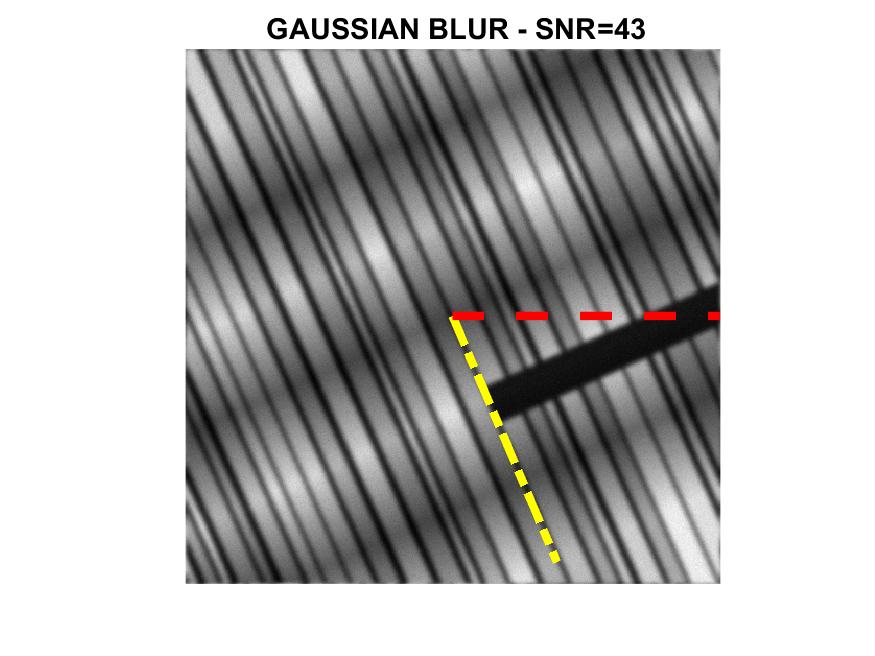}
        \includegraphics[width=5cm]{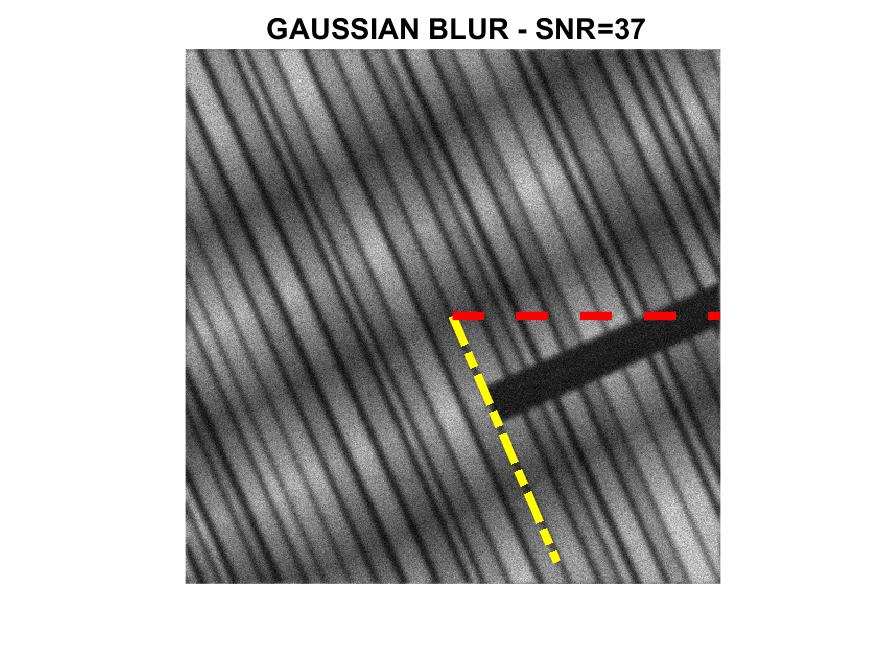}
    \end{center}
\end{subfigure}
\vspace{-9pt}
\caption{Test problem \texttt{phantom}: original and corrupted images.
    The yellow dash-dotted line indicates the direction estimated by Algorithm \ref{alg:direction_estimation}
    and the red dashed line the direction estimated by the method in \cite{kongskov:2017}.\label{fig:phantom}}
\end{figure}

\begin{figure}[htbp]
\begin{subfigure}{0.33\textwidth}
    \hspace*{6mm}\includegraphics[width=5cm]{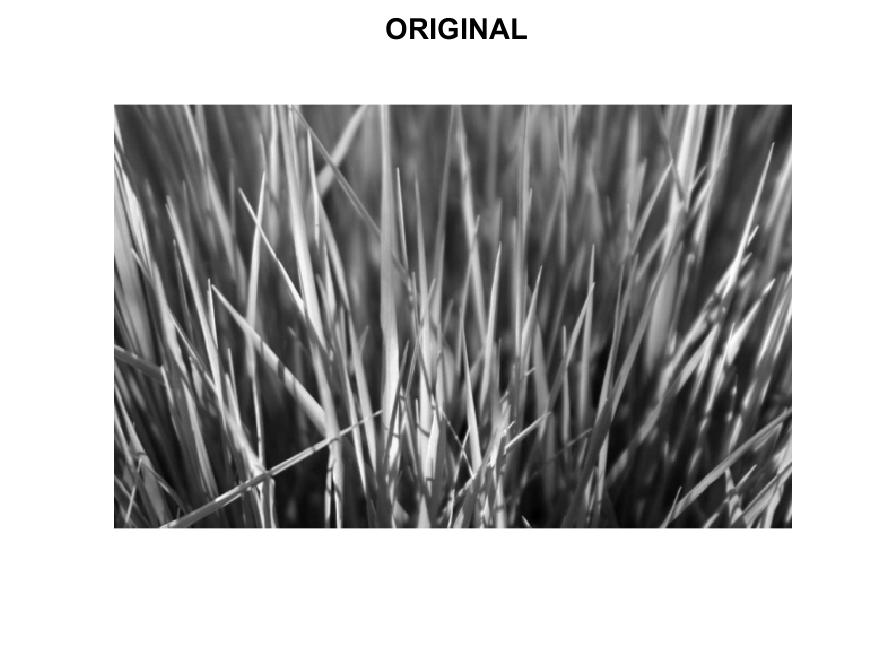}
\end{subfigure}
\begin{subfigure}{0.66\textwidth}
    \includegraphics[width=5cm]{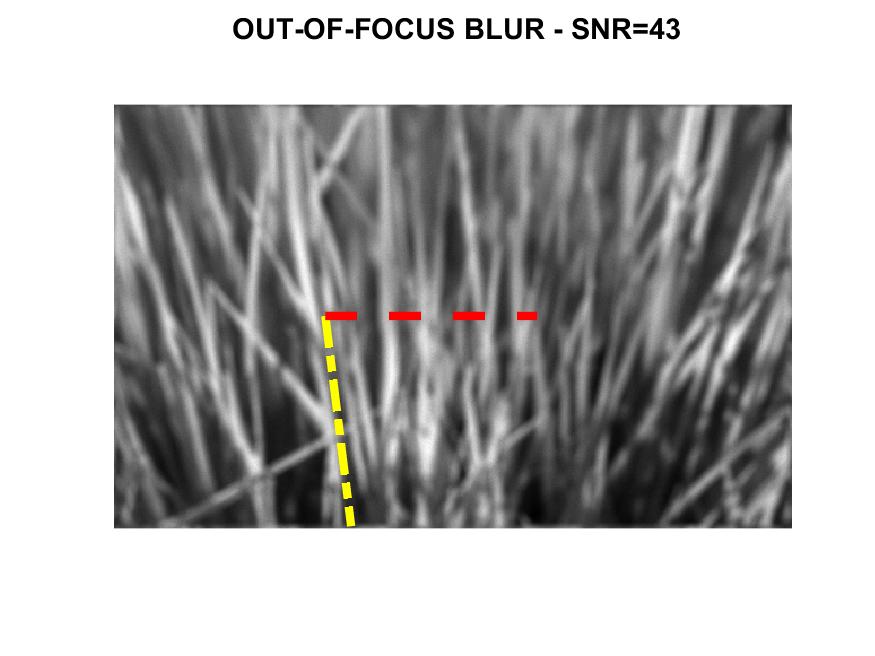}
    \includegraphics[width=5cm]{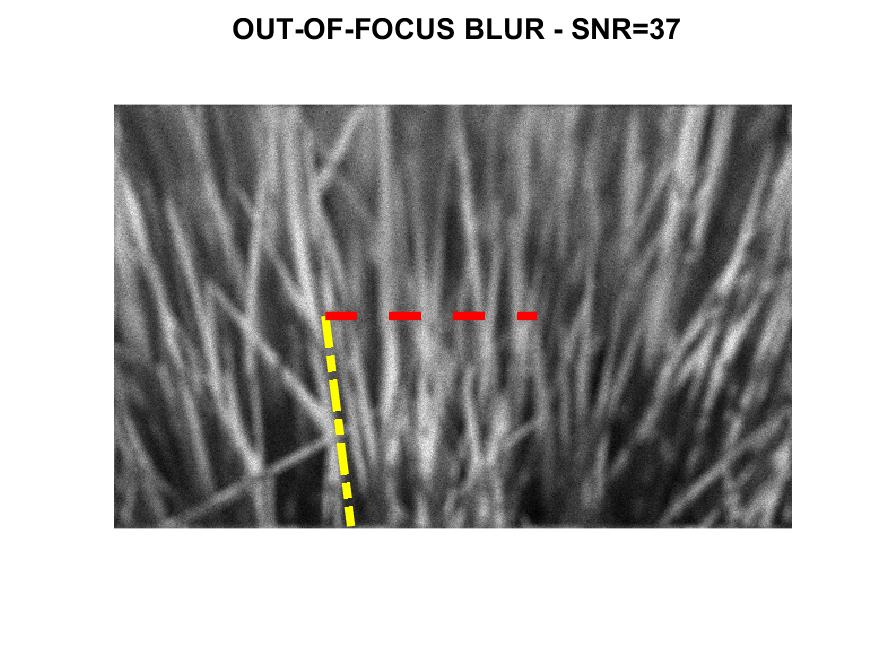} \\[-3mm]
    \includegraphics[width=5cm]{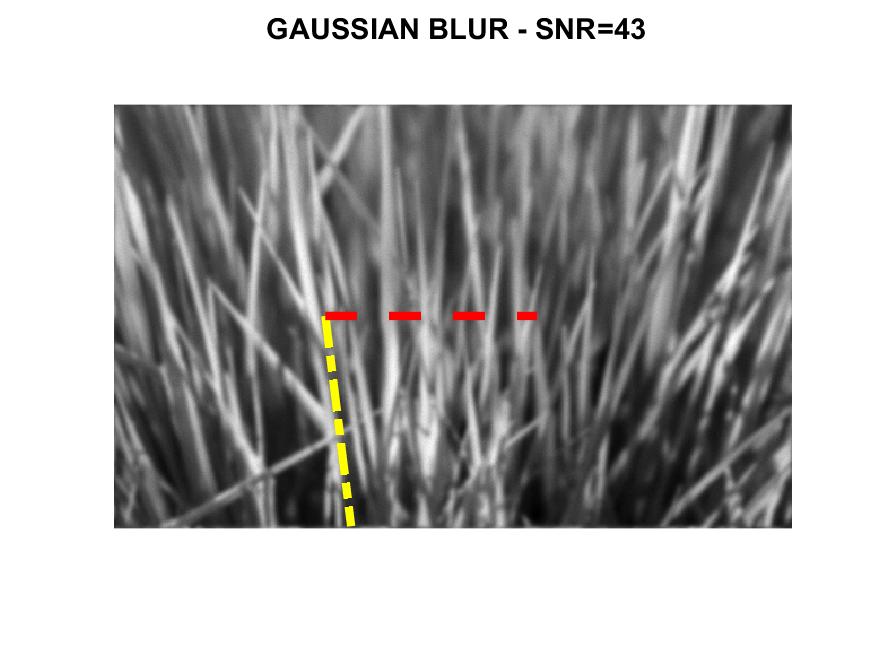}
    \includegraphics[width=5cm]{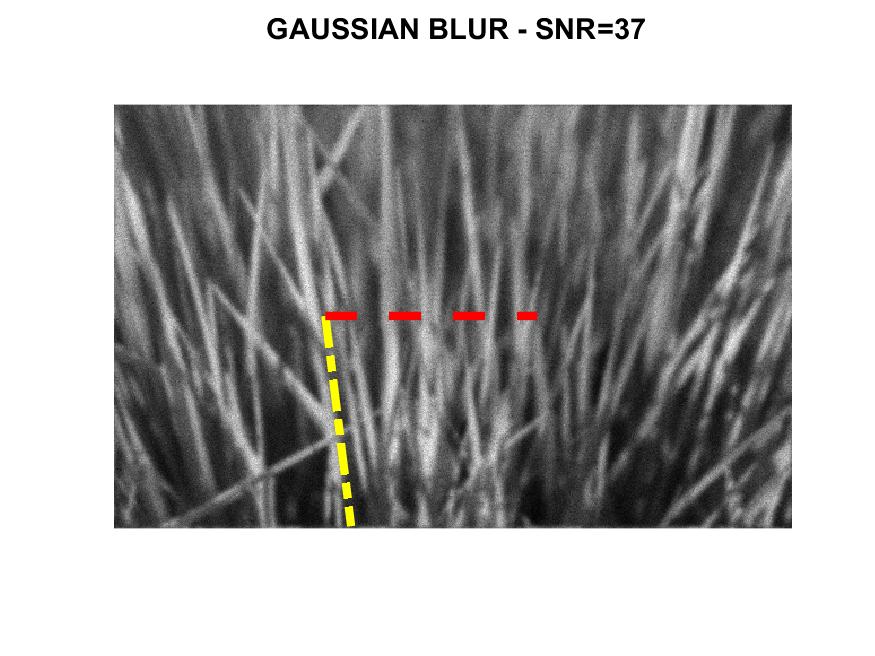}
\end{subfigure}
\vspace{-9pt}
\caption{Test problem \texttt{grass}: original and corrupted images.
    The yellow dash-dotted line indicates the direction estimated by Algorithm \ref{alg:direction_estimation}
    and the red dashed line the direction estimated by the method in \cite{kongskov:2017}.\label{fig:grass}}
\end{figure}

\begin{figure}[htbp]
\begin{subfigure}{0.33\textwidth}
    \begin{center}
        \includegraphics[width=5cm]{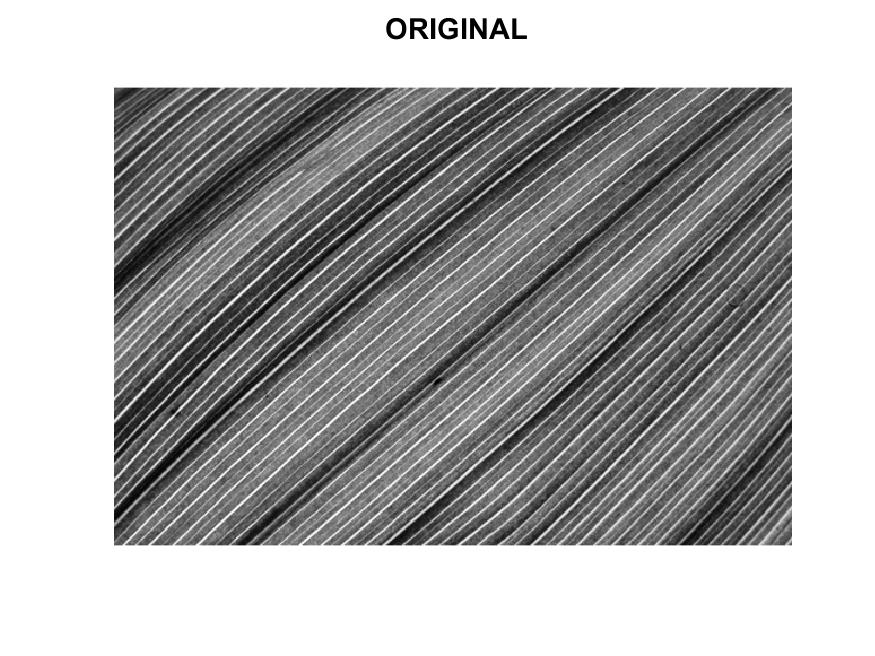}
    \end{center}
\end{subfigure}
\begin{subfigure}{0.66\textwidth}
    \begin{center}
        \includegraphics[width=5cm]{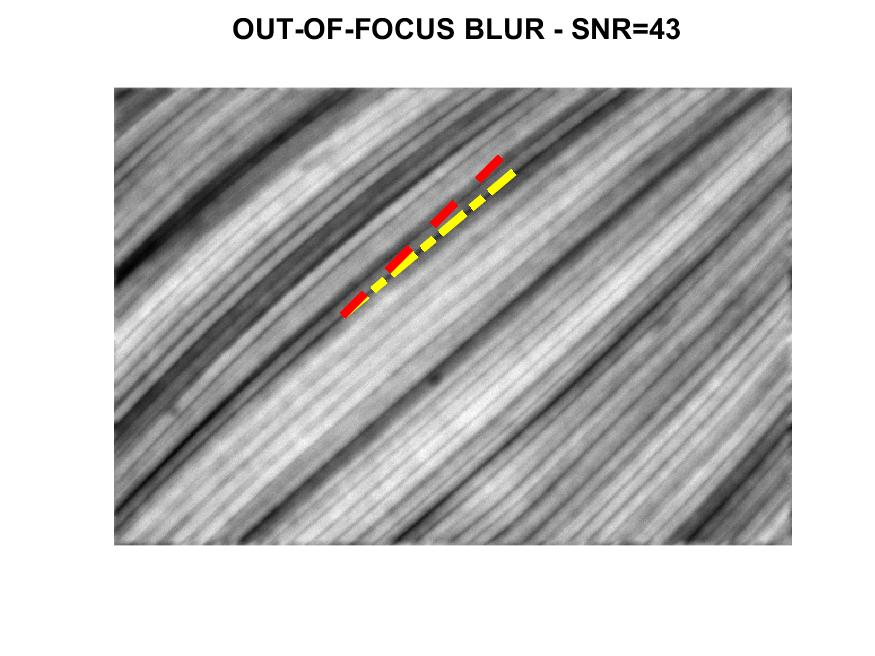}
        \includegraphics[width=5cm]{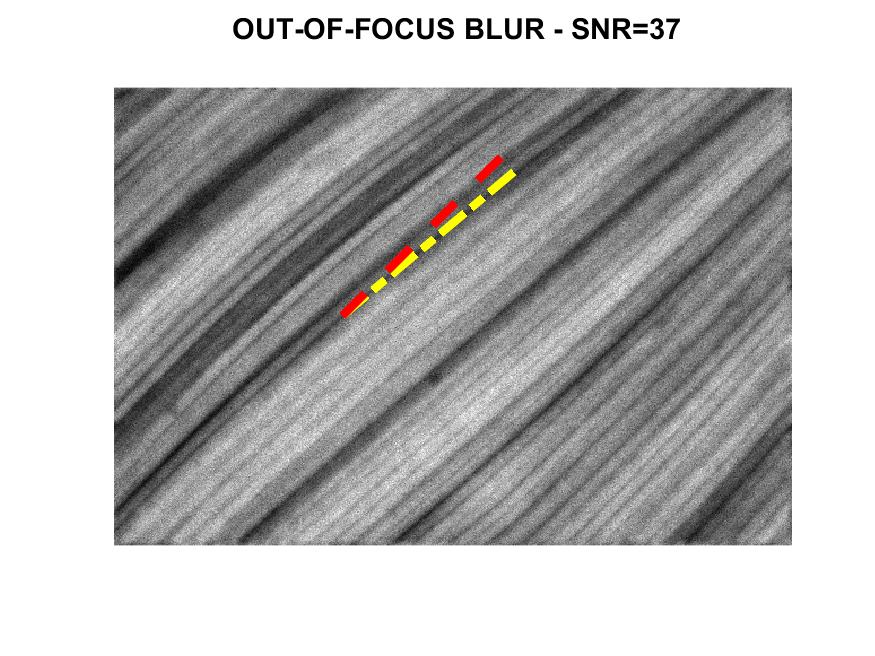} \\[-3mm]
        \includegraphics[width=5cm]{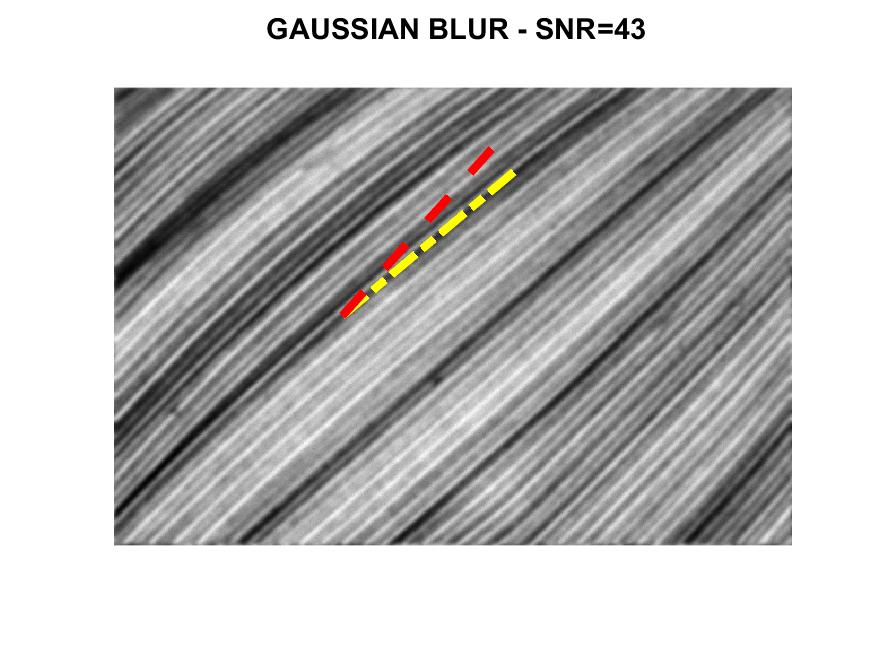}
        \includegraphics[width=5cm]{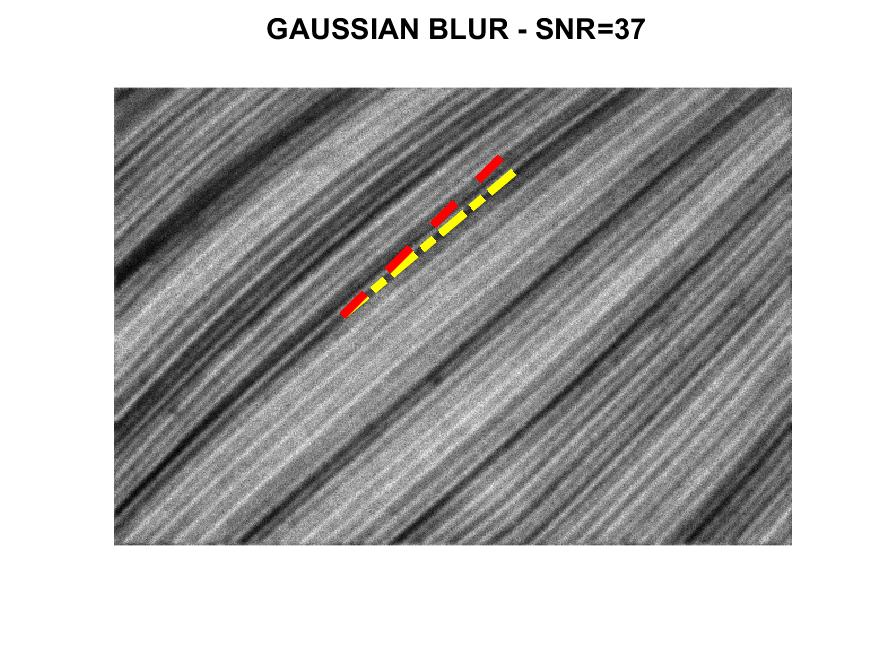}
    \end{center}
\end{subfigure}
\vspace{-9pt}
\caption{Test problem \texttt{leaves}: original and corrupted images.
    The yellow dash-dotted line indicates the direction estimated by Algorithm \ref{alg:direction_estimation}
    and the red dashed line the direction estimated by the method in \cite{kongskov:2017}.\label{fig:leaves}}
\end{figure}

\begin{figure}[htbp]
\begin{subfigure}{0.33\textwidth}
    \begin{center}
        \includegraphics[width=5cm]{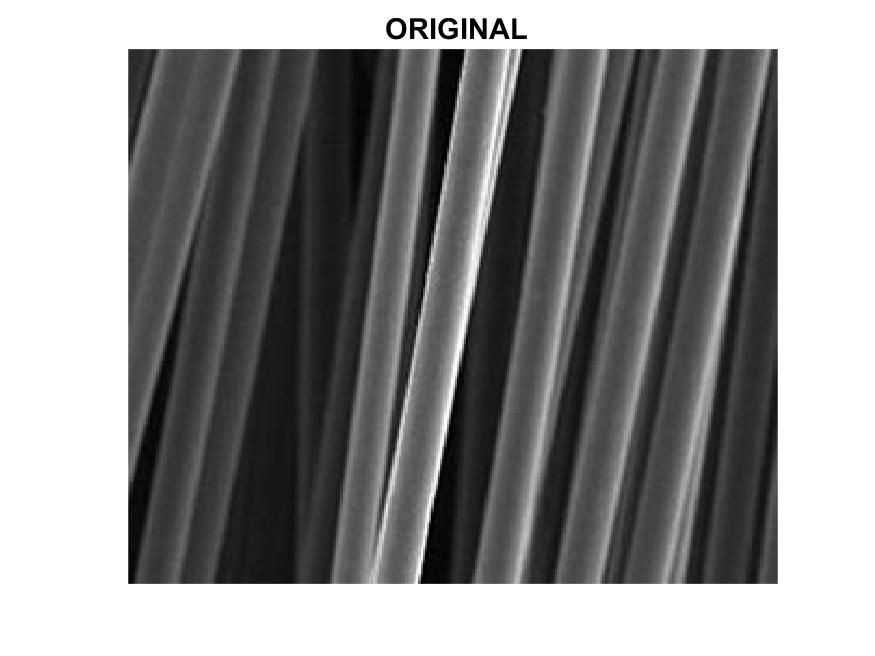}
    \end{center}
\end{subfigure}
\begin{subfigure}{0.66\textwidth}
    \begin{center}
        \includegraphics[width=5cm]{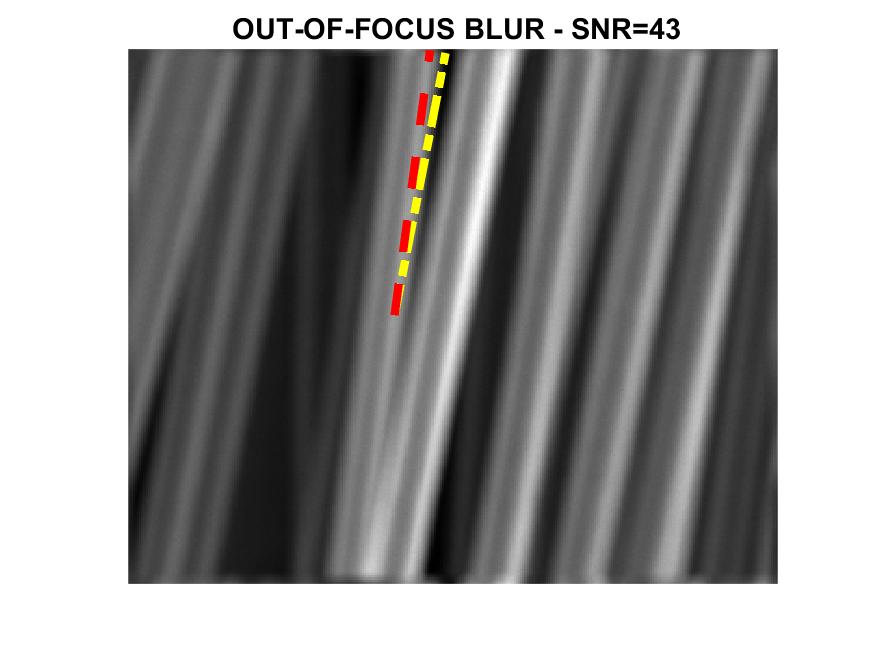}
        \includegraphics[width=5cm]{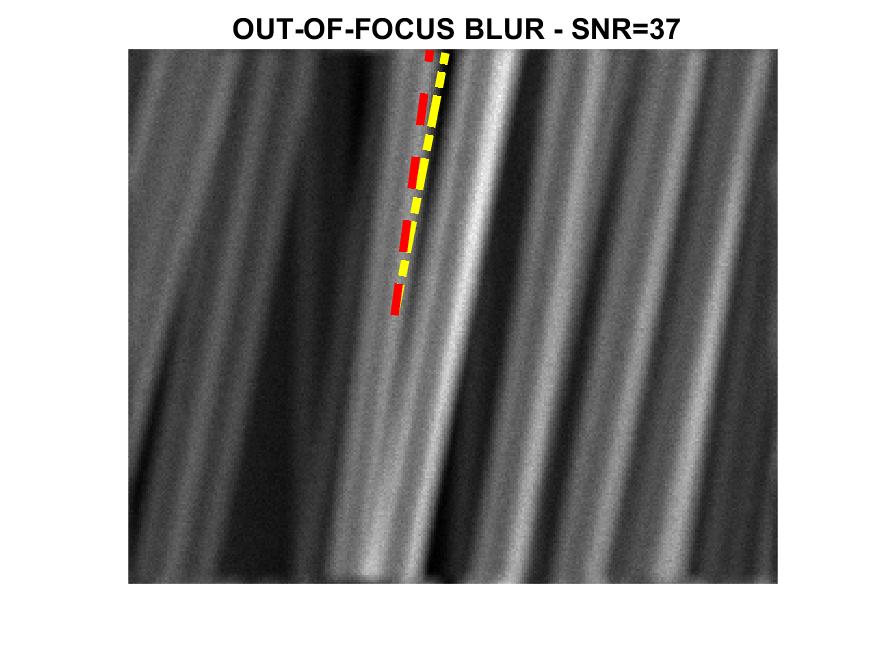} \\[-3mm]
        \includegraphics[width=5cm]{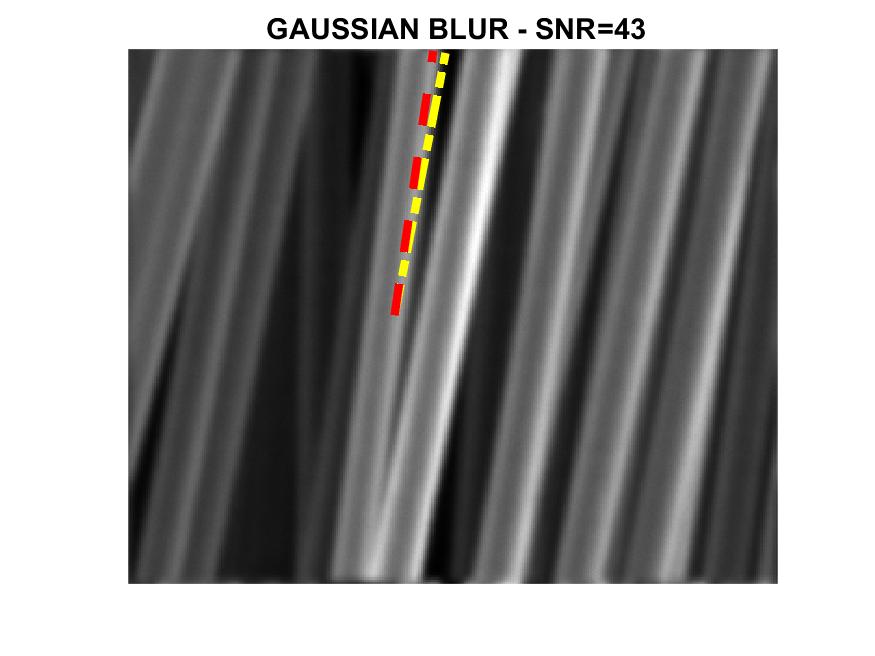}
        \includegraphics[width=5cm]{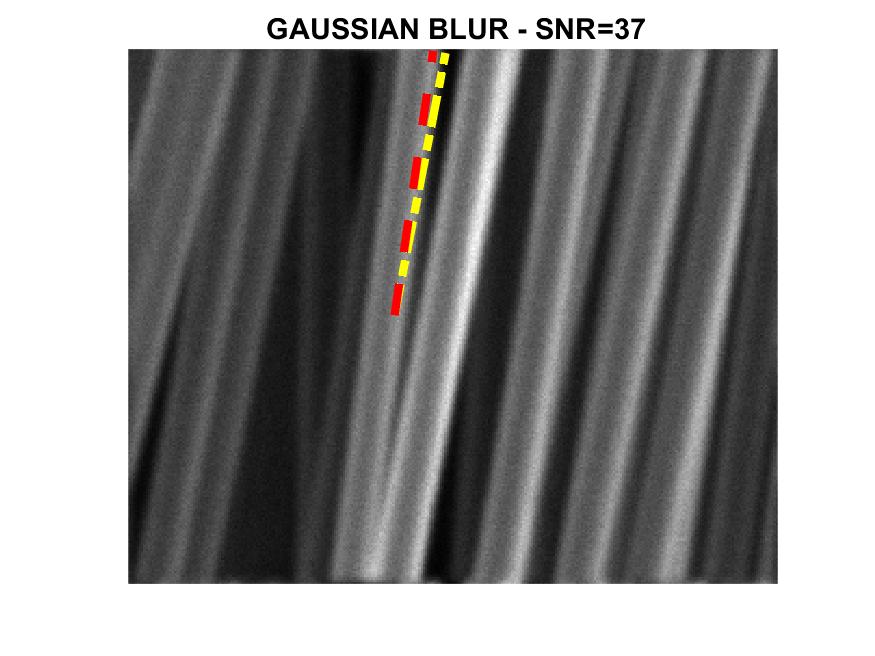}
    \end{center}
\end{subfigure}
\vspace{-9pt}
\caption{Test problem \texttt{carbon}: original and corrupted images.
    The yellow dash-dotted line indicates the direction estimated by Algorithm \ref{alg:direction_estimation}
    and the red dashed line the direction estimated by the method in \cite{kongskov:2017}.\label{fig:carbon}}
\end{figure}

Four directional images named \texttt{phantom} ($512\times512$), \texttt{grass} ($375\times600$),
\texttt{leaves} ($203\times300$)  and \texttt{carbon} ($247\times300$) were used in the experiments.
The first image is a piecewise affine fibre phantom image obtained with the
\texttt{fibre\_phantom\_pa} MATLAB function available from \url{http://www2.compute.dtu.dk/\~pcha/HDtomo/} (downloaded: 20 September 2020). 
The second and third images represent grass and veins of leaves, respectively, which naturally exhibit
a directional structure. The last image is a Scanning Electron Microscope (SEM) image of carbon fibres.
The images are shown in {Figures}~\ref{fig:phantom}--\ref{fig:carbon}.

To simulate experimental data, each reference image was convolved with two PSFs, one corresponding
to a Gaussian blur with variance 2, generated by the \texttt{psfGauss} function from~\cite{nagy:2004},
and the other corresponding to an out-of-focus blur with radius 5, obtained with the
function \texttt{fspecial} from the MATLAB Image Processing Toolbox.
{To take into account the existence of some background emission, a constant term  $\gamma$ equal to $10^{-10}$ was added to all pixels of the blurry image.}
The resulting image was corrupted by Poisson noise, using the MATLAB
function \texttt{imnoise}. The intensities of the original images
were pre-scaled to get noisy and blurry images with Signal to Noise Ratio (SNR) equal to 43 and 37 dB.
We recall that in the case of Poisson noise, which affects the photon counting process,
the SNR is estimated as \cite{bertero:2018}
\begin{equation*}
    \text{SNR} = 10\log_{10}\left(\frac{N_\text{exact}}{\sqrt{N_\text{exact} + N_\text{background}}} \right),
\end{equation*}
where $N_\text{exact}$ and $N_\text{background}$ are the total number of photons in
the image to be recovered and in the background term, respectively. Finally, the corrupted images
were scaled to have their maximum intensity values equal to 1. For each test problem, the noisy and blurry
images are shown in Figures~\ref{fig:phantom}--\ref{fig:carbon}.
%
%
%

%
%

\subsection{Direction Estimation\label{sec:estimation_results}}
\begin{figure}[htbp]
\begin{center}
    \includegraphics[width=3.5cm]{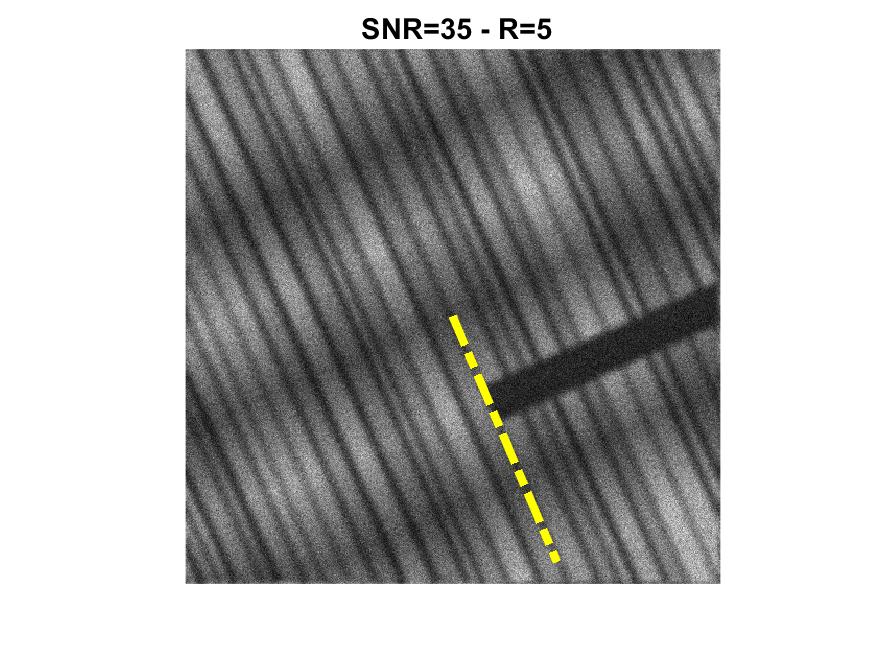}\hspace{-.6cm}
    \includegraphics[width=3.5cm]{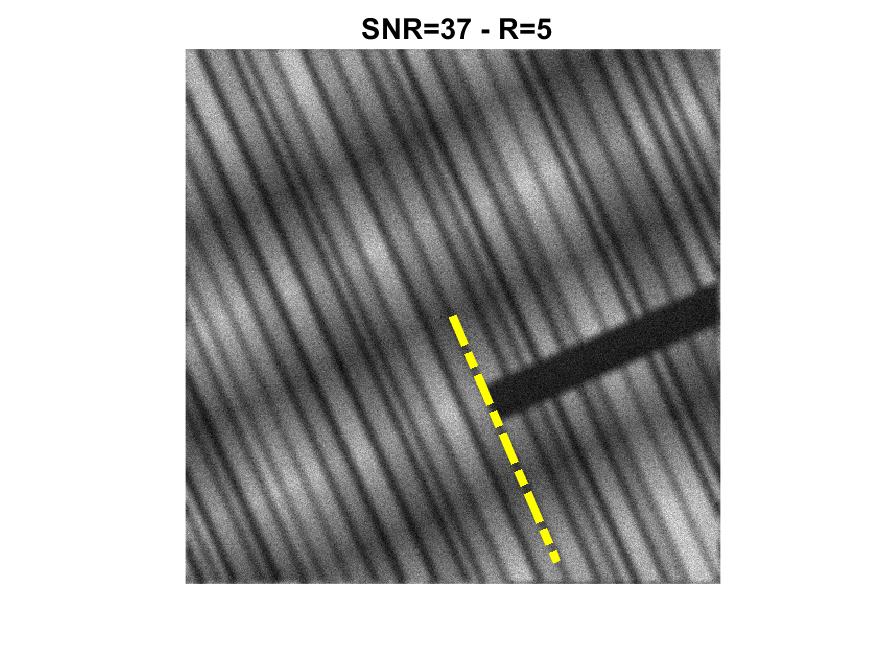}\hspace{-.6cm}
    \includegraphics[width=3.5cm]{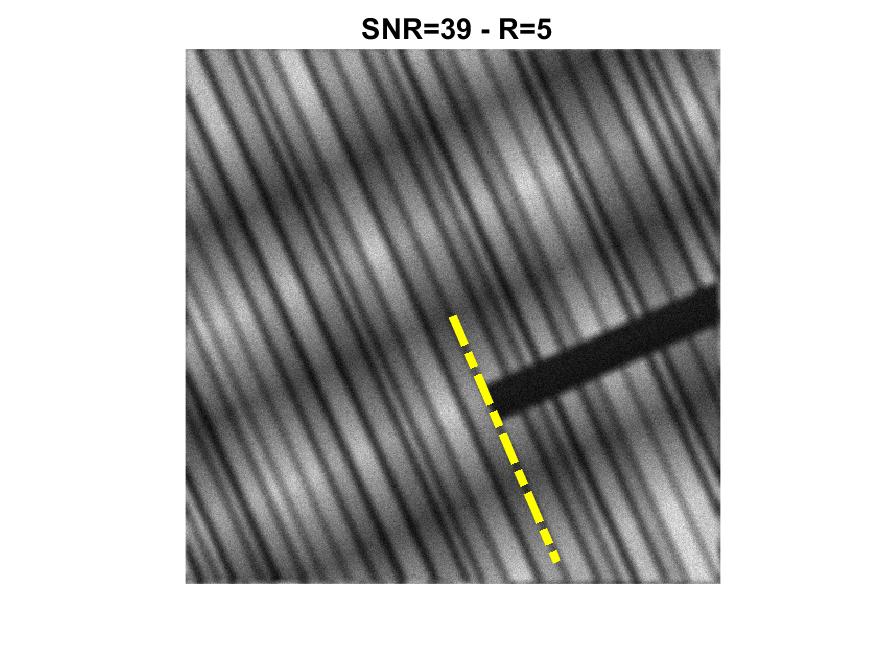}\hspace{-.6cm}
    \includegraphics[width=3.5cm]{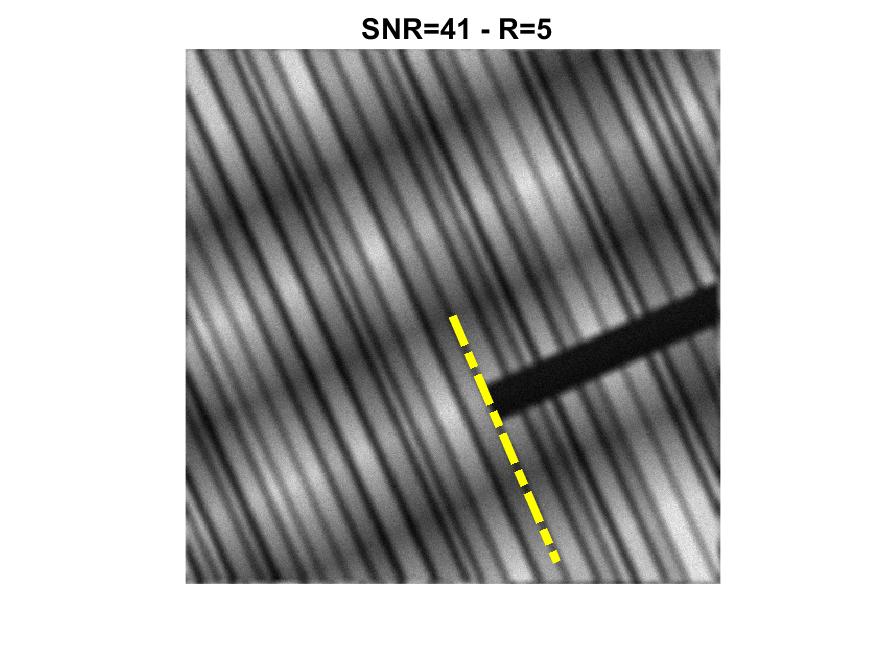}\hspace{-.6cm}
    \includegraphics[width=3.5cm]{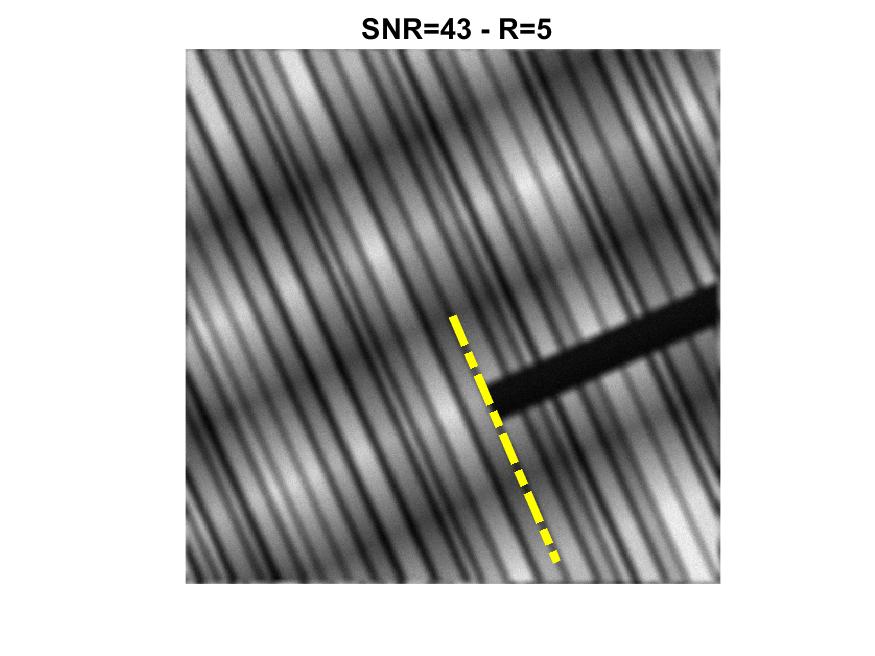}
    \\[-3mm]
    \includegraphics[width=3.5cm]{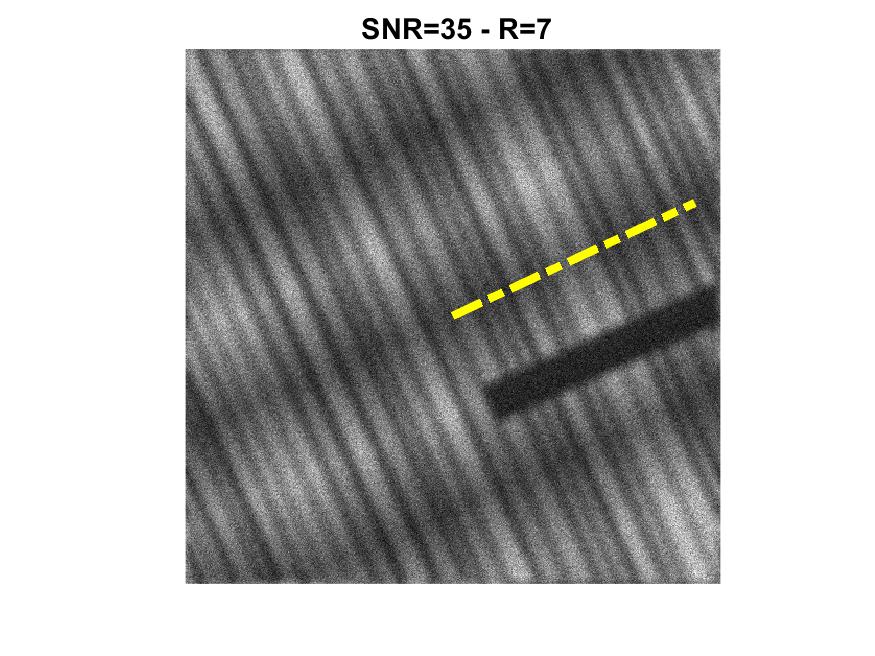}\hspace{-.6cm}
    \includegraphics[width=3.5cm]{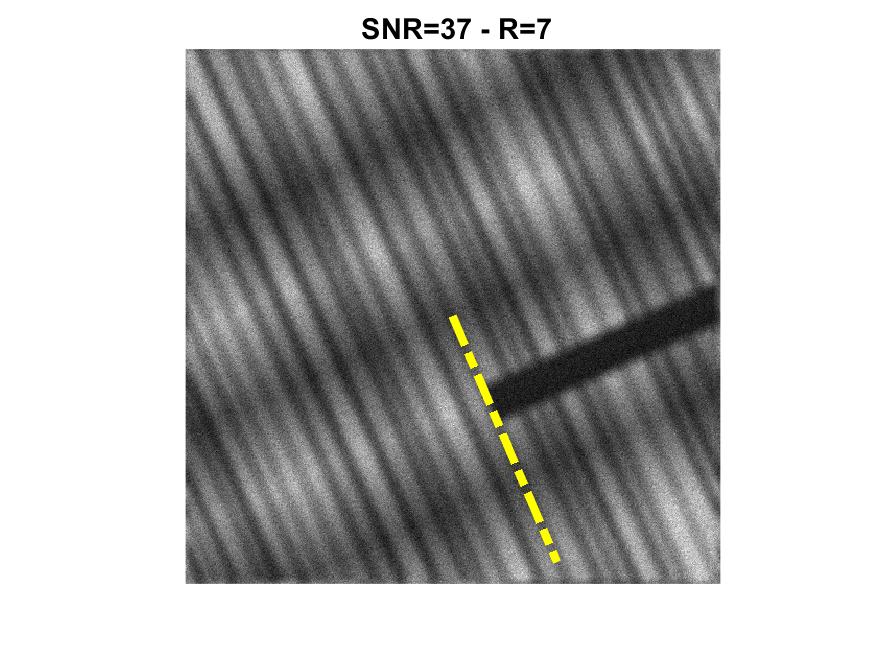}\hspace{-.6cm}
    \includegraphics[width=3.5cm]{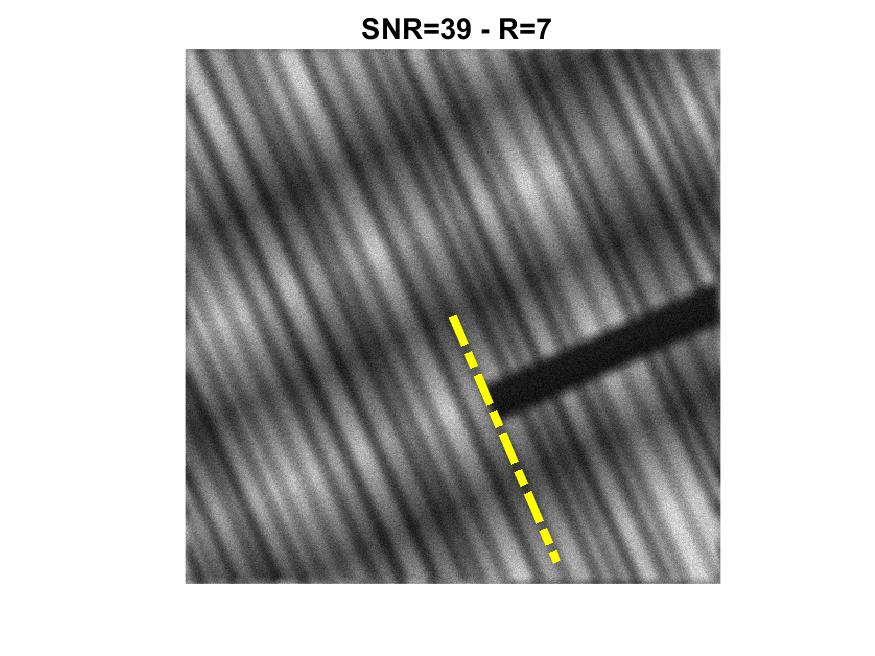}\hspace{-.6cm}
    \includegraphics[width=3.5cm]{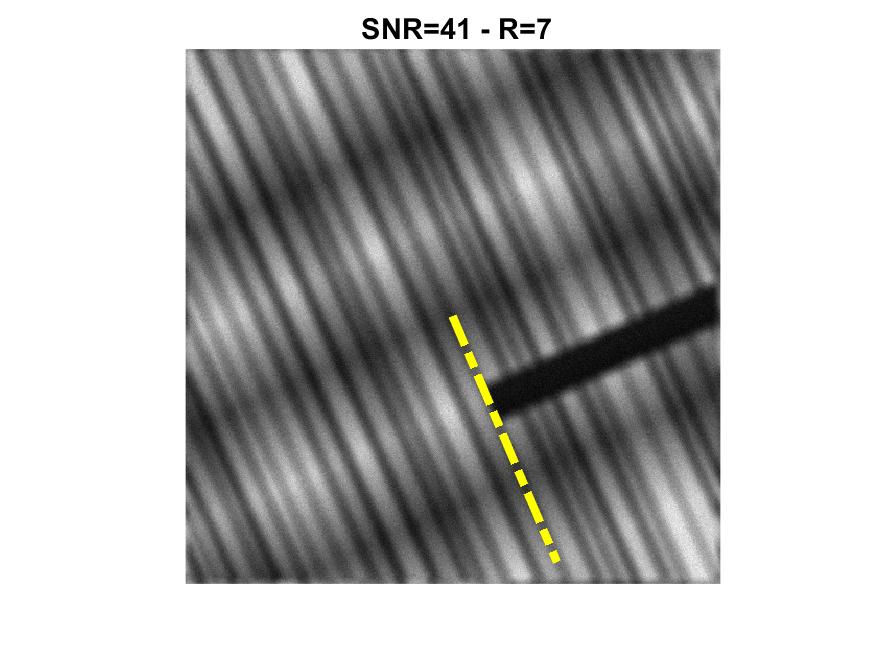}\hspace{-.6cm}
    \includegraphics[width=3.5cm]{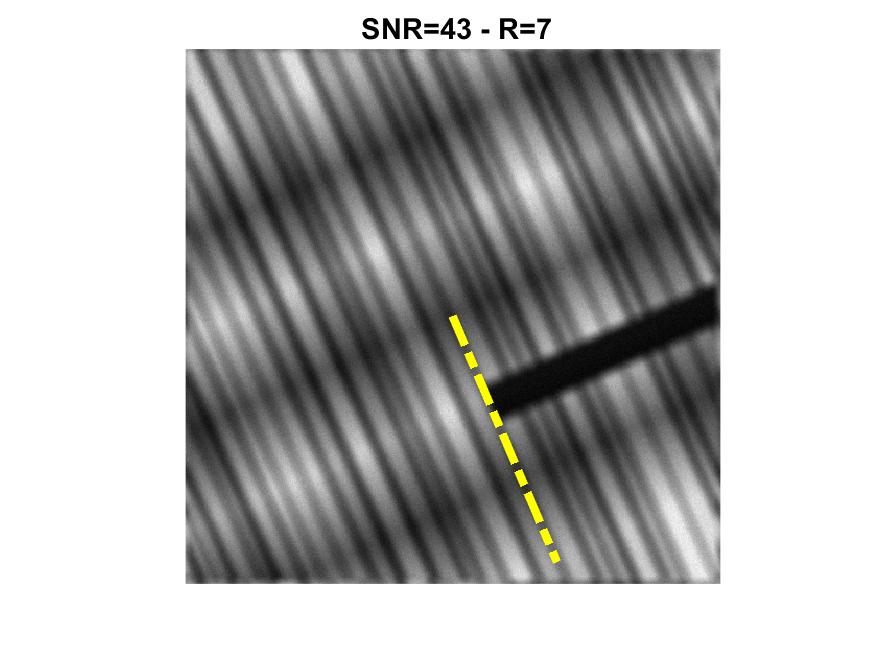}
    \\[-3mm]
    \includegraphics[width=3.5cm]{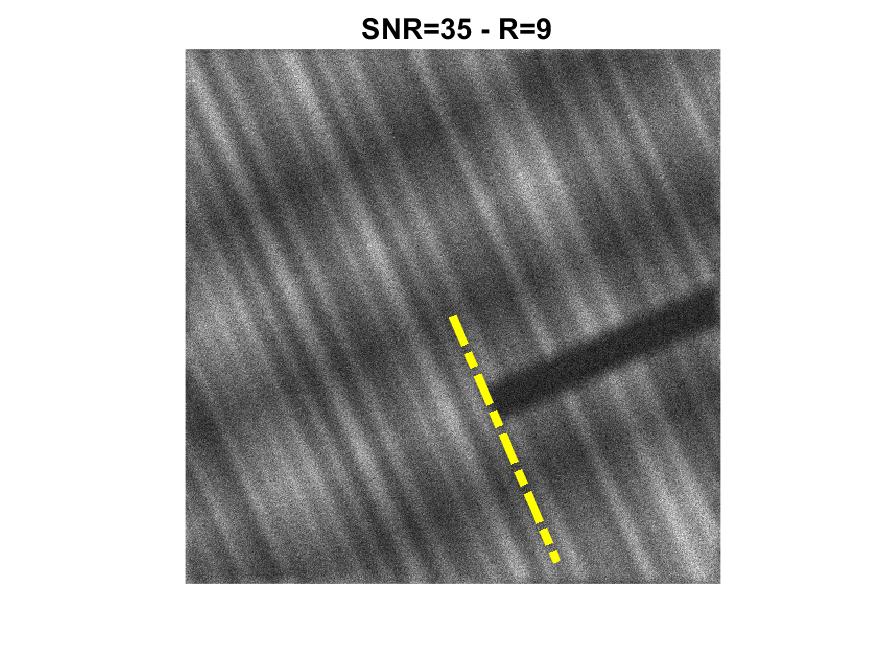}\hspace{-.6cm}
    \includegraphics[width=3.5cm]{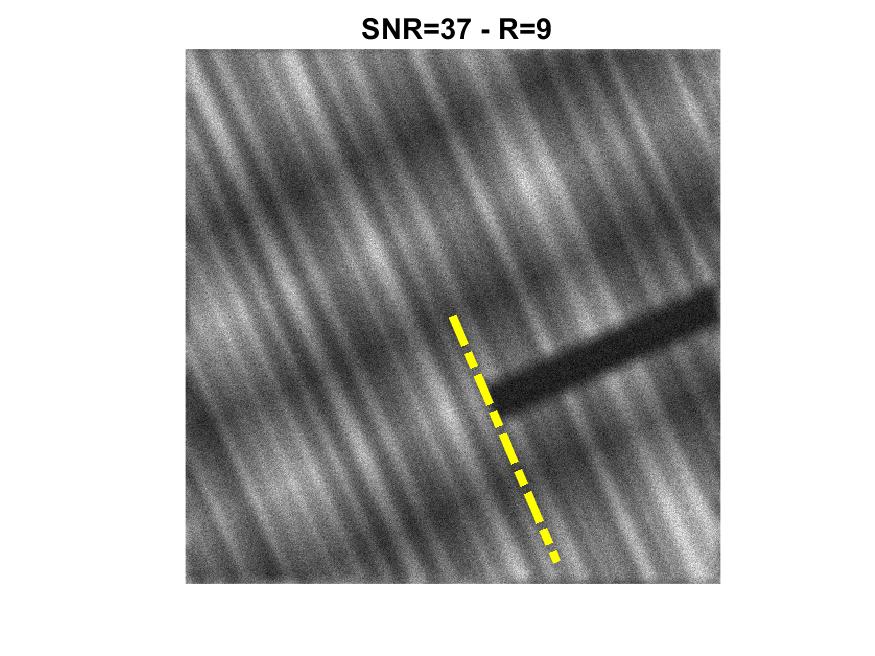}\hspace{-.6cm}
    \includegraphics[width=3.5cm]{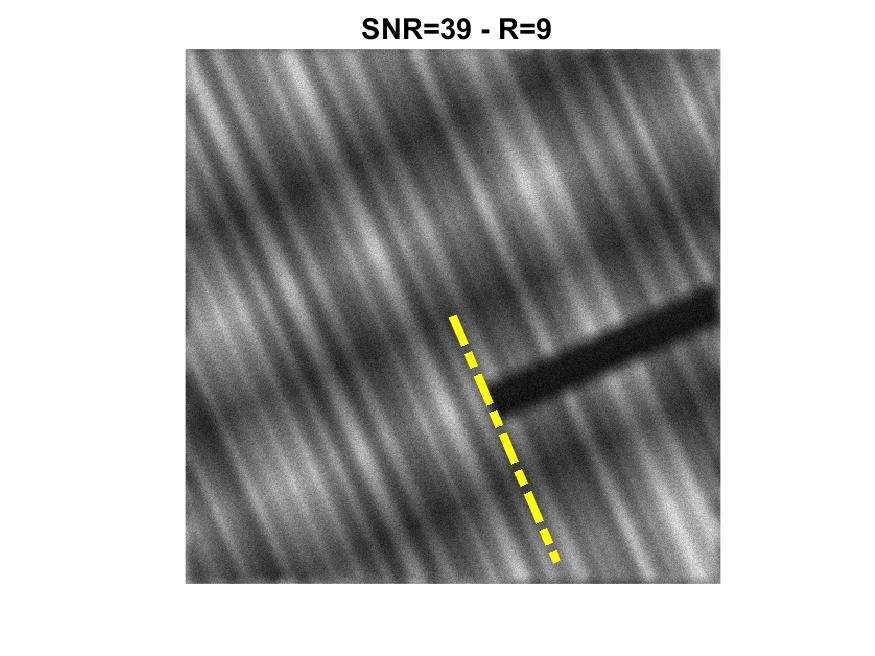}\hspace{-.6cm}
    \includegraphics[width=3.5cm]{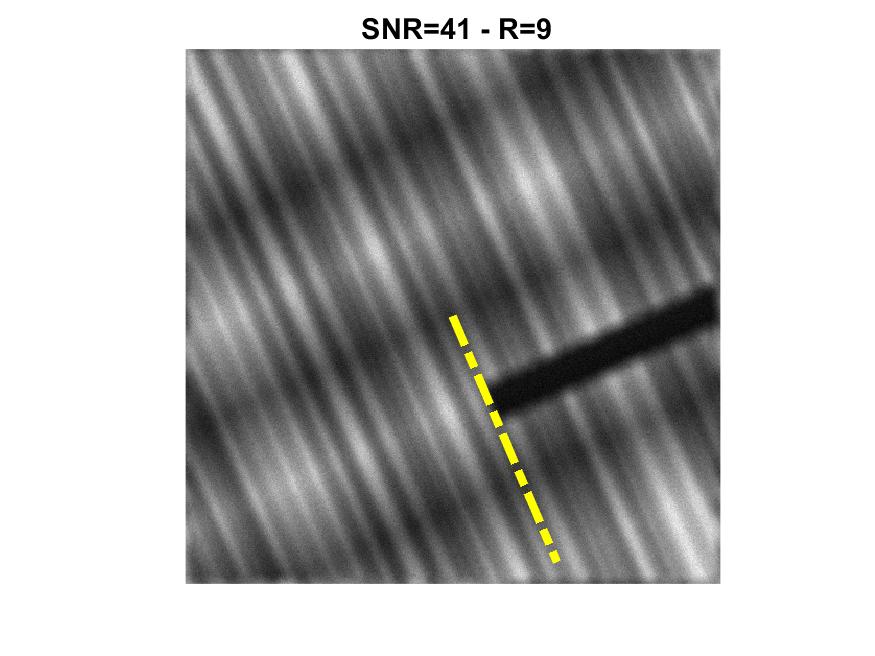}\hspace{-.6cm}
    \includegraphics[width=3.5cm]{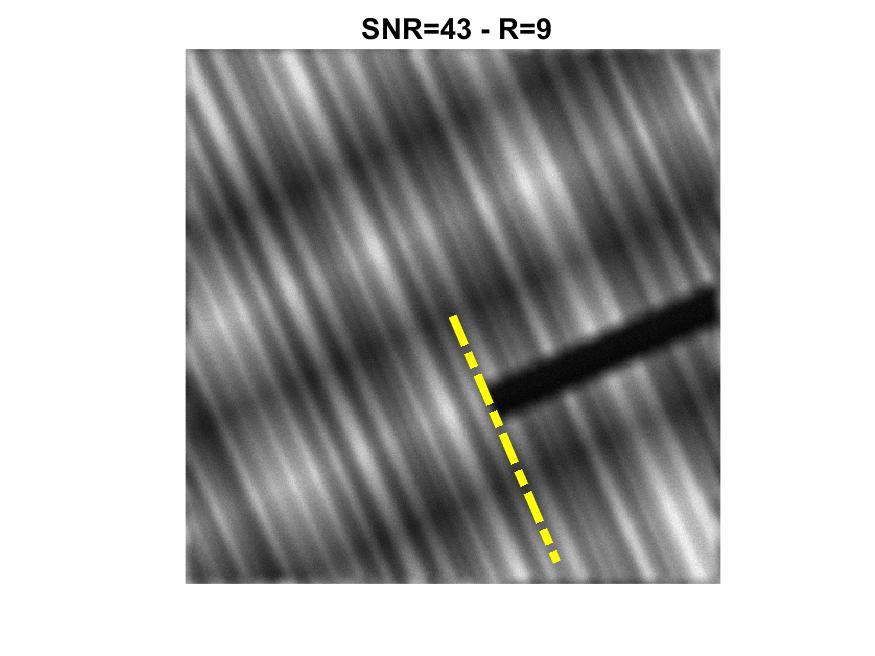}
\end{center}    
\vspace{-9pt}
\caption{Direction estimation for \texttt{phantom} with SNR $=35,37,39,41,43$ dB (from left to right) and out-of-focus blur with radius R $=5, 7,9$ (from top to bottom).\label{fig:phantomdir}}
\end{figure}

\begin{figure}[htbp]
\begin{center}
    \includegraphics[width=4.5cm]{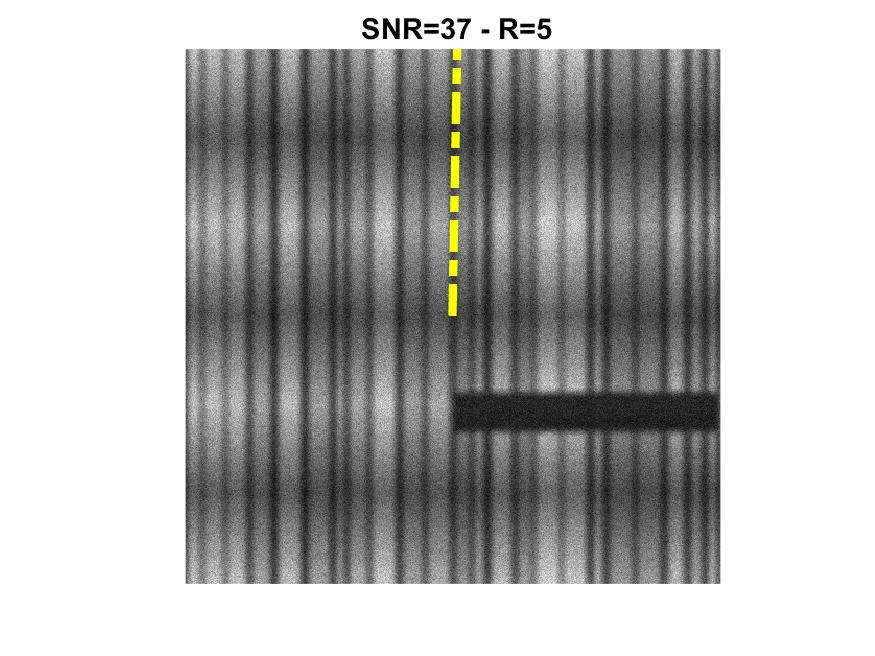}\hspace{-.6cm}
    \includegraphics[width=4.5cm]{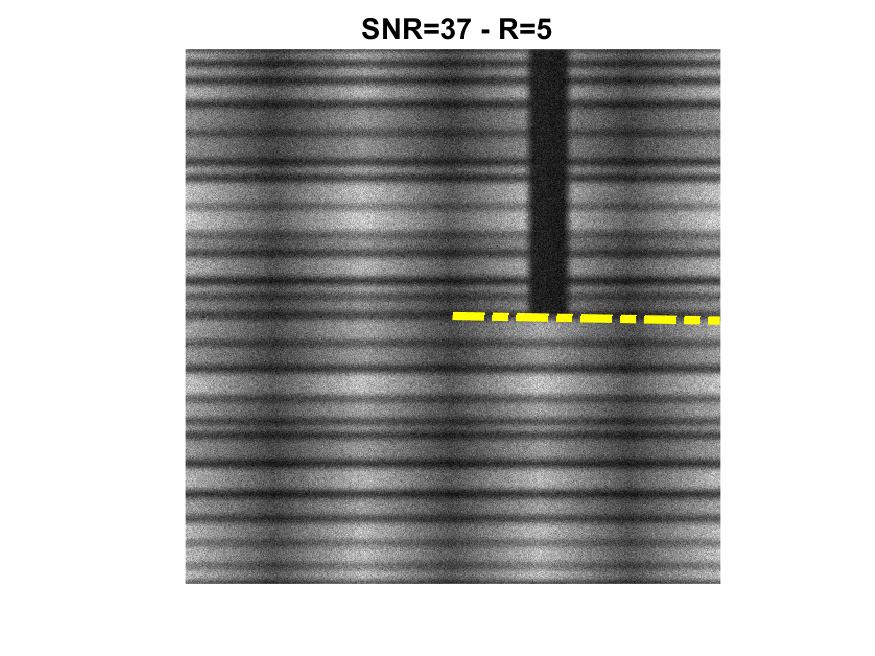}\hspace{-.6cm}
    \includegraphics[width=4.5cm]{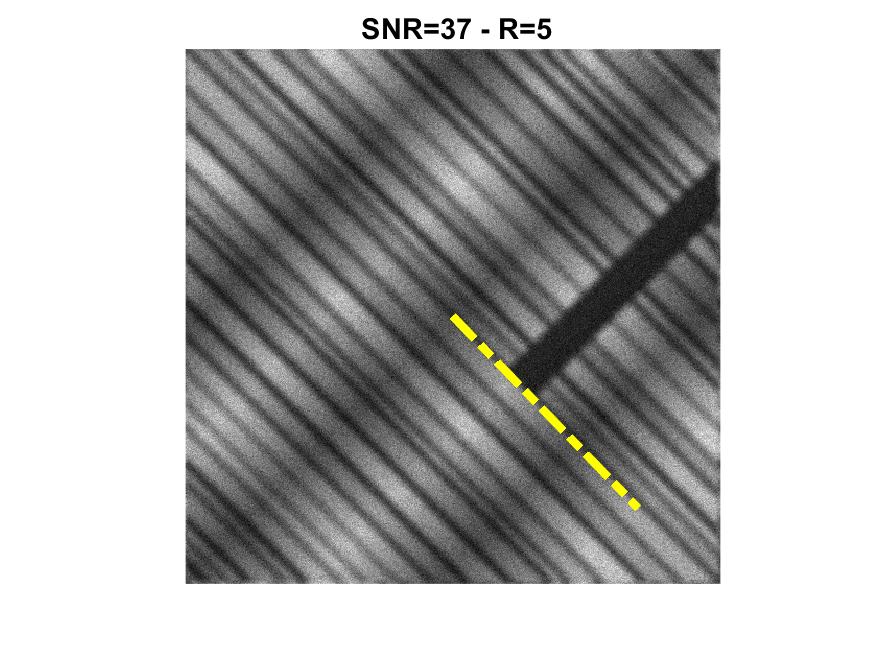} \\[-3mm]
    \includegraphics[width=4.5cm]{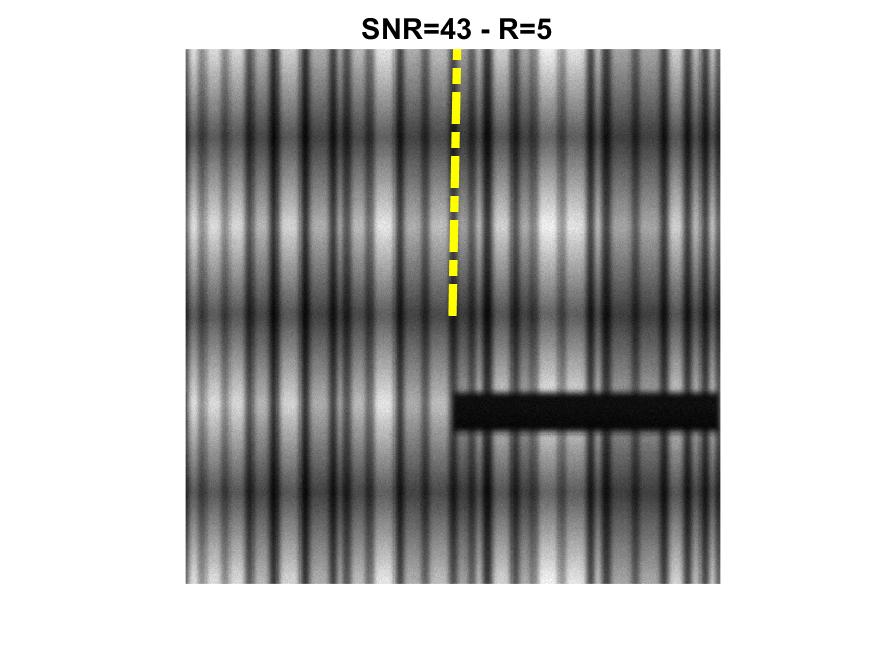}\hspace{-.6cm}
    \includegraphics[width=4.5cm]{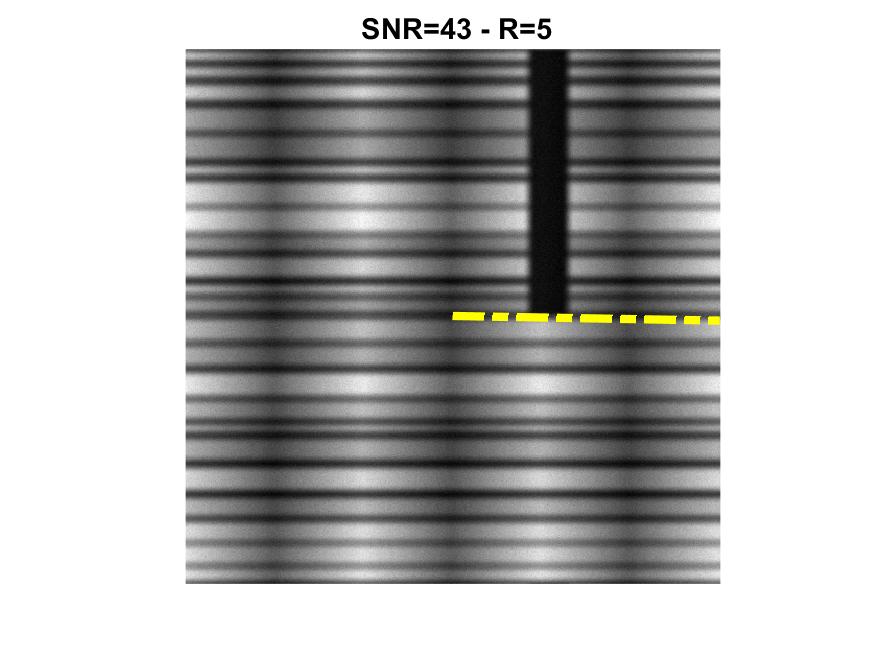}\hspace{-.6cm}
    \includegraphics[width=4.5cm]{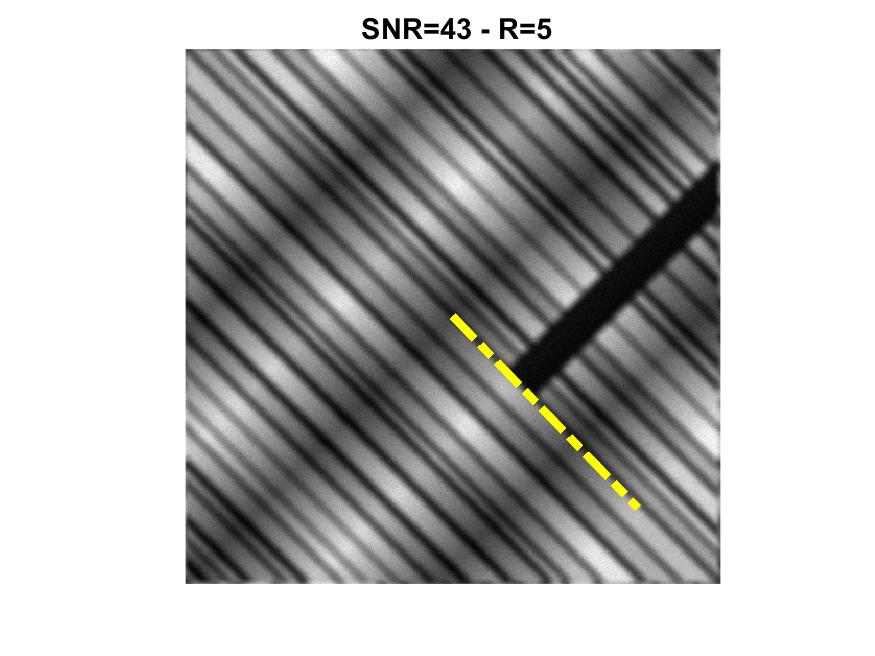}
\end{center}
\vspace{-9pt}
\caption{Direction estimation for \texttt{phantom} with SNR $=37, 43$ (top, bottom) and out-of-focus blur with radius R $=5$.\label{fig:phantomdir2}}
\end{figure}

In Figures~\ref{fig:phantom}--\ref{fig:carbon} we compare Algorithm~\ref{alg:direction_estimation}
with the algorithm proposed in~\cite{kongskov:2017}, showing that Algorithm~\ref{alg:direction_estimation} always  
correctly estimates the main direction of the four test images. We also test the robustness of our algorithm with respect to noise and blur.
In Figure~\ref{fig:phantomdir} we show the estimated main direction of the \texttt{phantom} image corrupted
by Poisson noise with SNR $=35,37,39,41,43$ dB and out-of-focus blurs with radius R $= 5, 7, 9$.
In only one case (SNR $=35$, R $=7$) Algorithm~\ref{alg:direction_estimation} fails, returning as estimate
the orthogonal direction, i.e., the direction corresponding to the large black line and the background color gradient.
Finally, we test Algorithm~\ref{alg:direction_estimation} on a phantom image with vertical, horizontal
and diagonal main directions corresponding to $\theta = 0,90,45$. The results, in Figure~\ref{fig:phantomdir2}, show
that our algorithm is not sensitive to the specific directional structure of the image.
%
%

\subsection{Image Deblurring}
We compare the quality of the restorations obtained by using the DTGV$^2$ and TGV$^2$ regularizers and ADMM for the solution of both models. 
In all the tests, the value of the penalty parameter was set as $\rho=10$ and the value of the stopping threshold as $tol=10^{-4}$. A maximum number of $k_\text{max}=500$ iterations was allowed. By following \cite{kongskov:2017,kongskov:2019}, the weight parameters of DTGV were chosen as $\alpha_0=\beta$ and $\alpha_1=(1-\beta)$ with $\beta=2/3$.
For each test problem, the value of the regularization parameter $\lambda$ was tuned by a trial-and-error strategy. This strategy consisted in running ADMM with initial guess $\bu^0 = \bb$ several times on each test image, varying the value of $\lambda$ at each execution. For all the runs the stopping criterion for ADMM and the values of $\alpha_0$, $\alpha_1$ and $\rho$ were the same as described above. The value of $\lambda$ yielding the smallest Root Mean Square Error (RMSE) at the last iteration was chosen as the ``optimal'' value.

\begin{figure}[b!]
\begin{center}
\vspace*{-12pt}
\includegraphics[width=0.25\columnwidth]{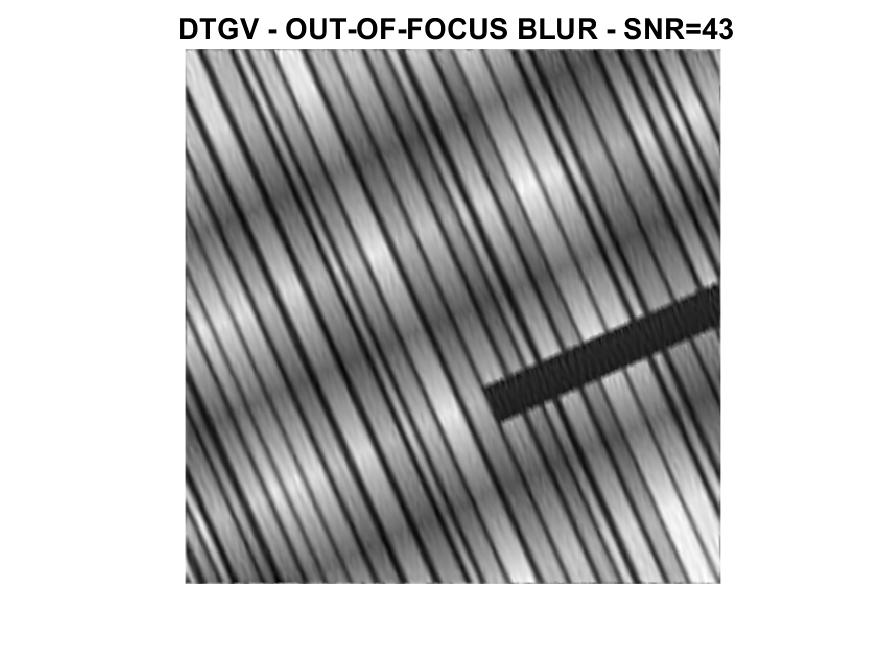}\hspace{-.6cm}
\includegraphics[width=0.25\columnwidth]{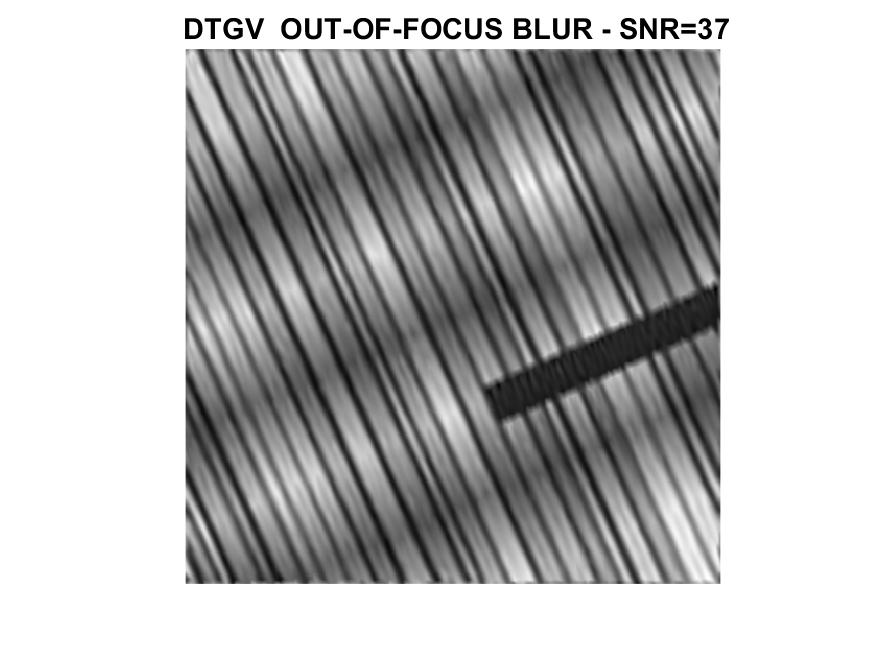}\hspace{-.6cm}
\includegraphics[width=0.25\columnwidth]{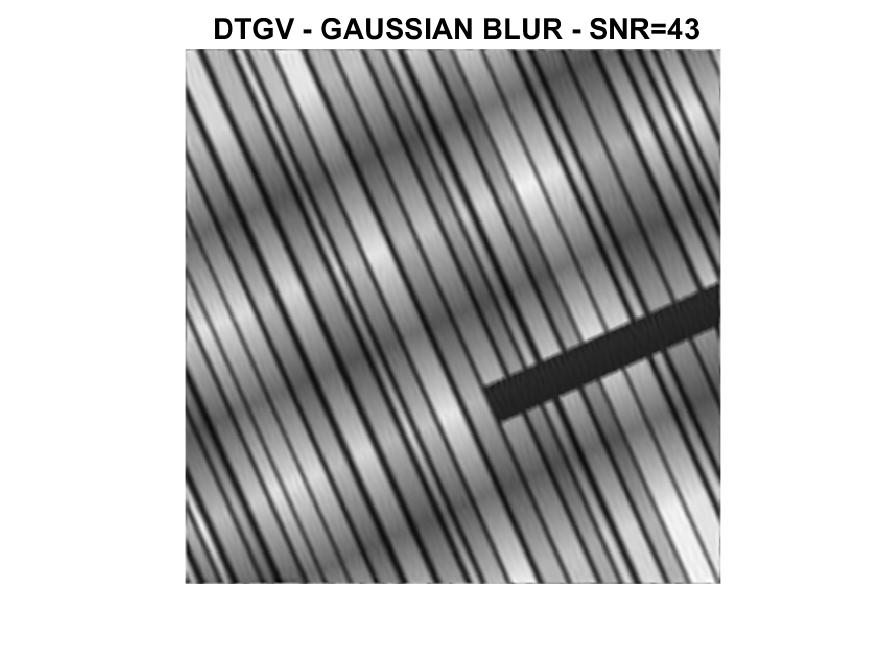}\hspace{-.6cm}
\includegraphics[width=0.25\columnwidth]{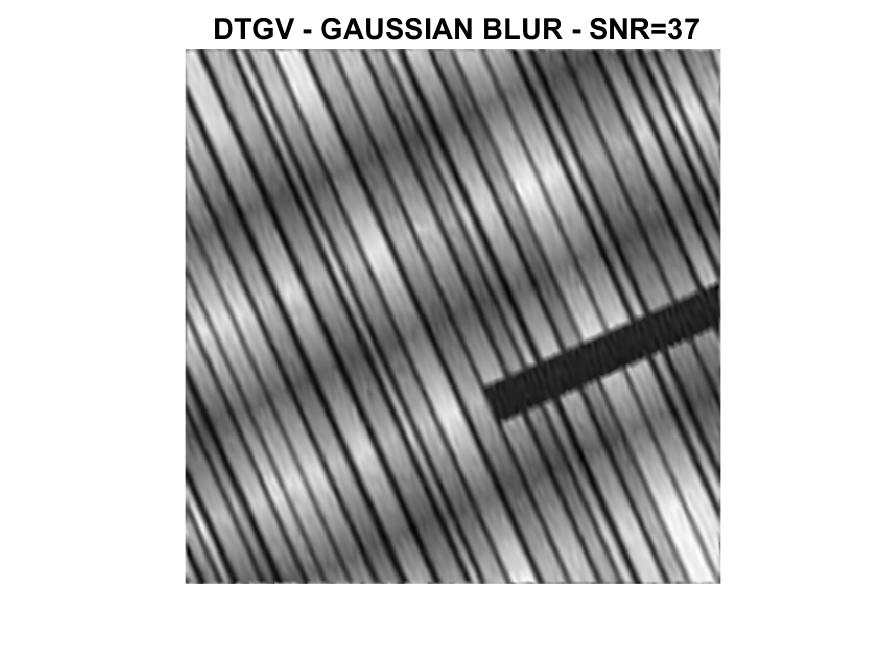} \\[-3mm]
\includegraphics[width=0.25\columnwidth]{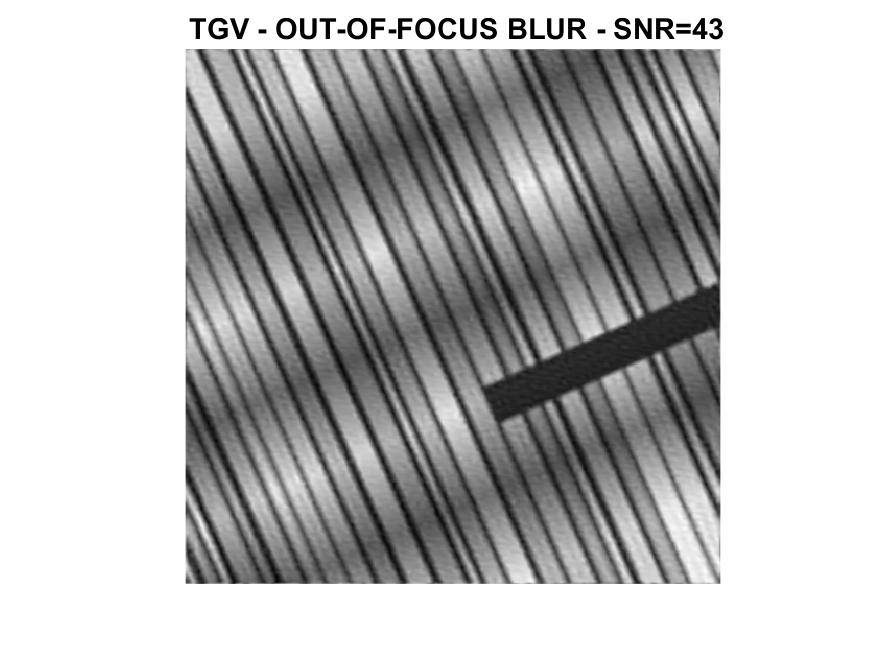}\hspace{-.6cm}
\includegraphics[width=0.25\columnwidth]{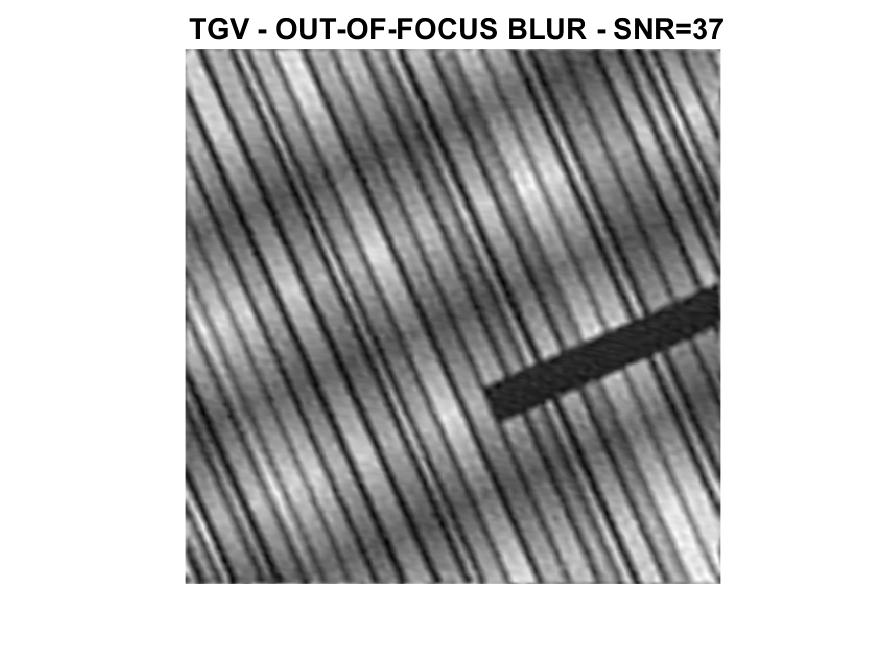}\hspace{-.6cm}
\includegraphics[width=0.25\columnwidth]{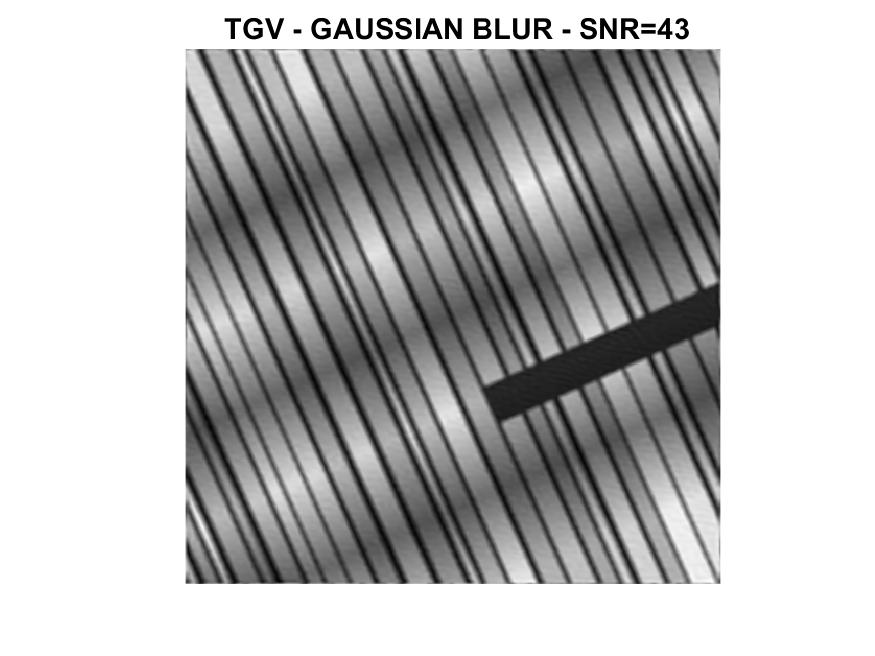}\hspace{-.6cm}
\includegraphics[width=0.25\columnwidth]{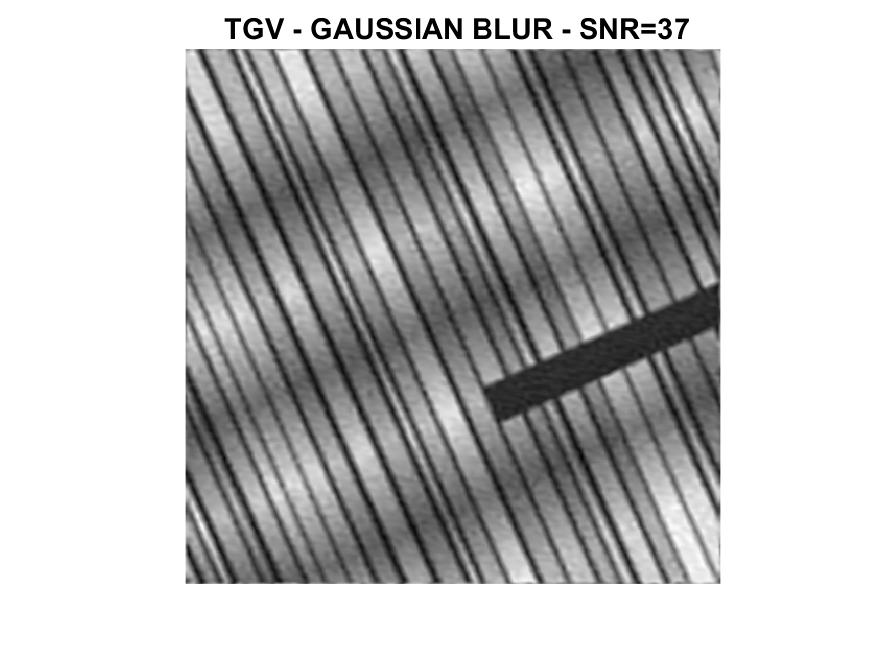} \\[-3mm]
\caption{Test problem \texttt{phantom}: images restored with DTGV$^2$ (\textbf{top}) and TGV$^2$ (\textbf{bottom}).\label{fig:phantomrec}}
\end{center}
\end{figure}

\begin{figure}[htbp]
\begin{center}
\includegraphics[width=0.25\columnwidth]{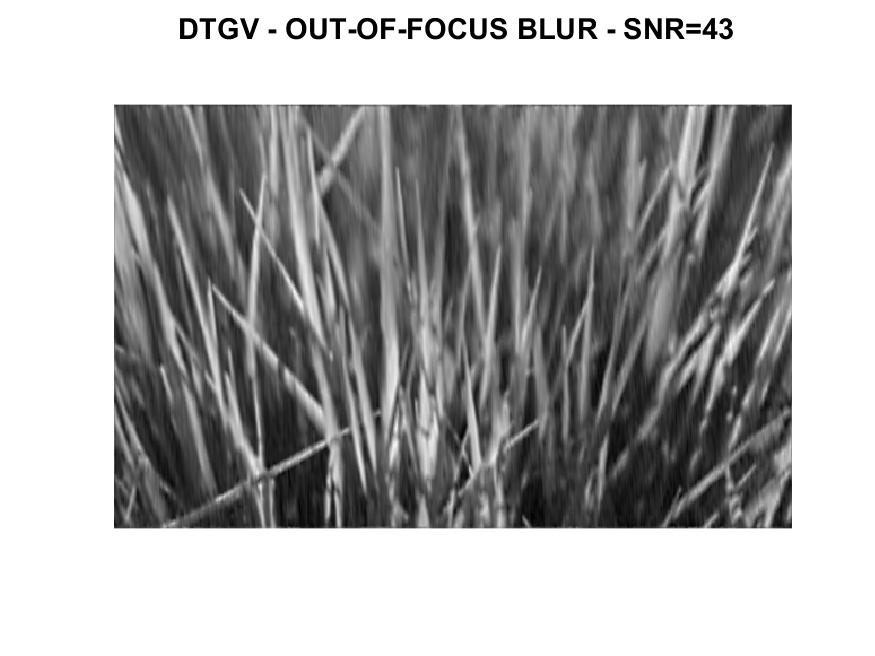}\hspace{-.6cm}
\includegraphics[width=0.25\columnwidth]{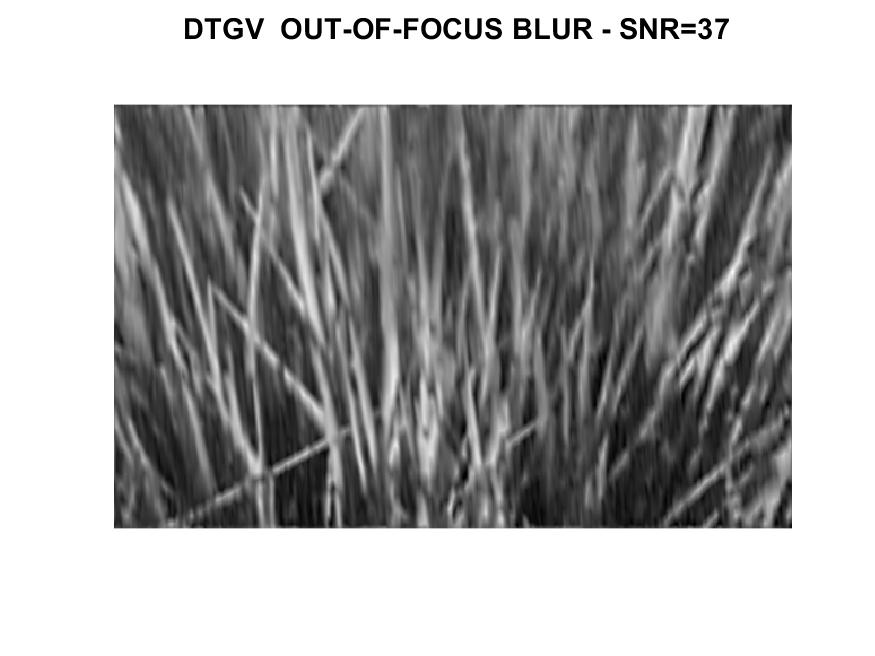}\hspace{-.6cm}
\includegraphics[width=0.25\columnwidth]{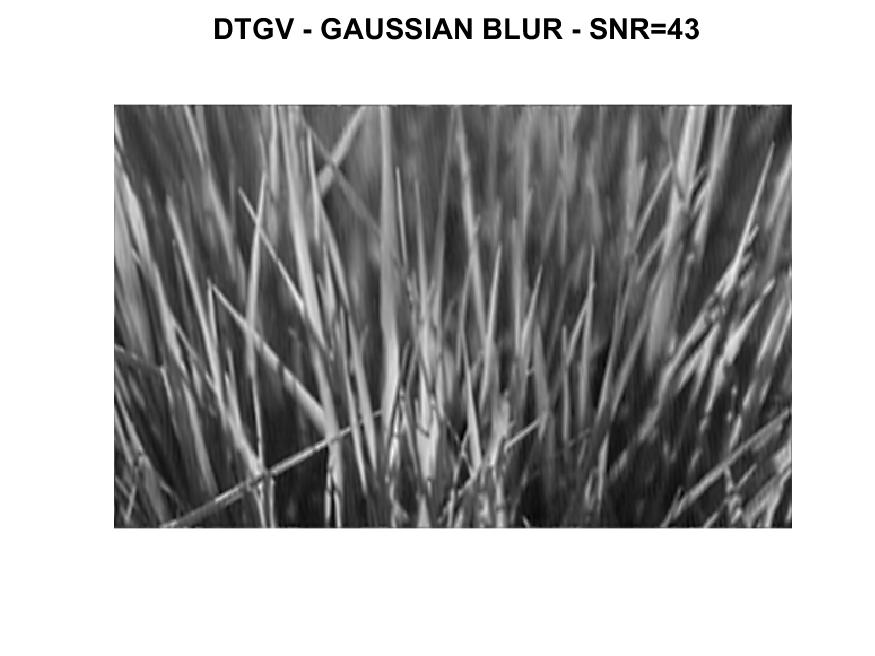}\hspace{-.6cm}
\includegraphics[width=0.25\columnwidth]{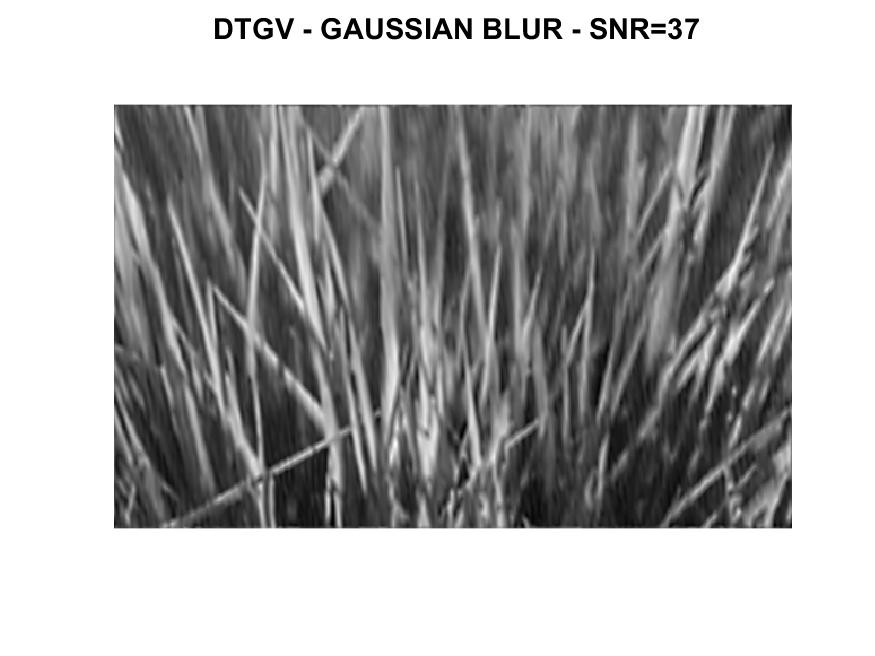} \\[-3mm]
\includegraphics[width=0.25\columnwidth]{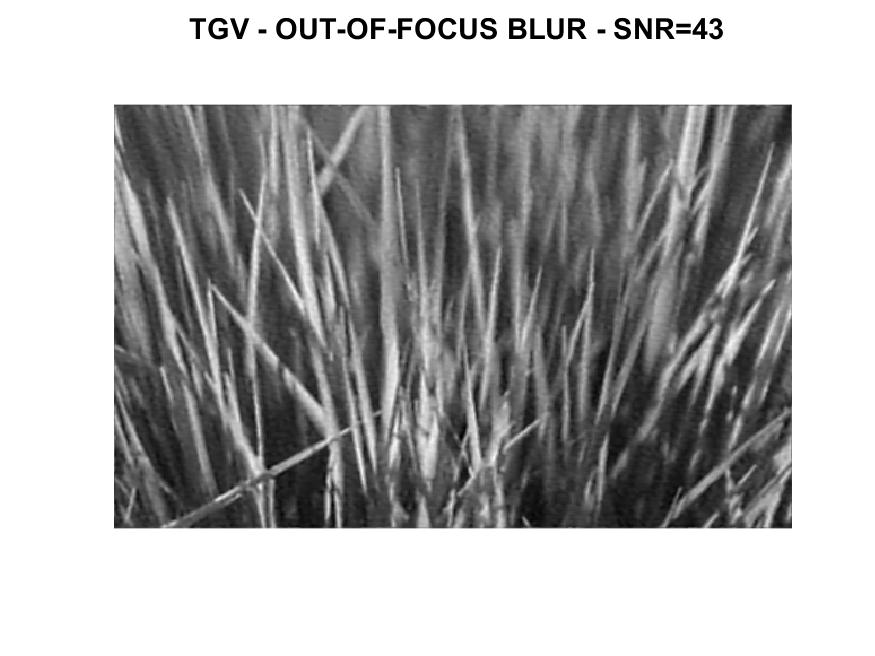}\hspace{-.6cm}
\includegraphics[width=0.25\columnwidth]{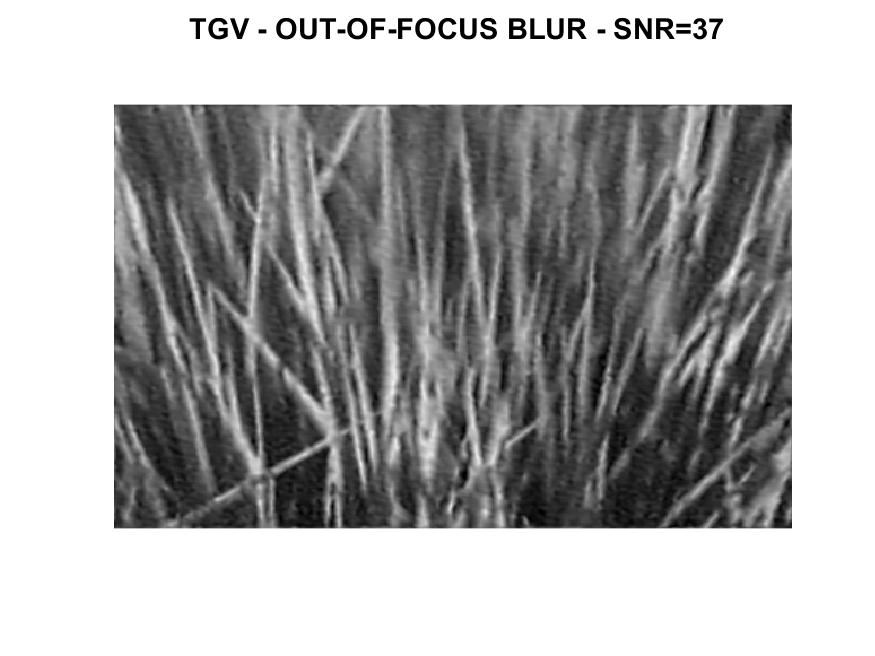}\hspace{-.6cm}
\includegraphics[width=0.25\columnwidth]{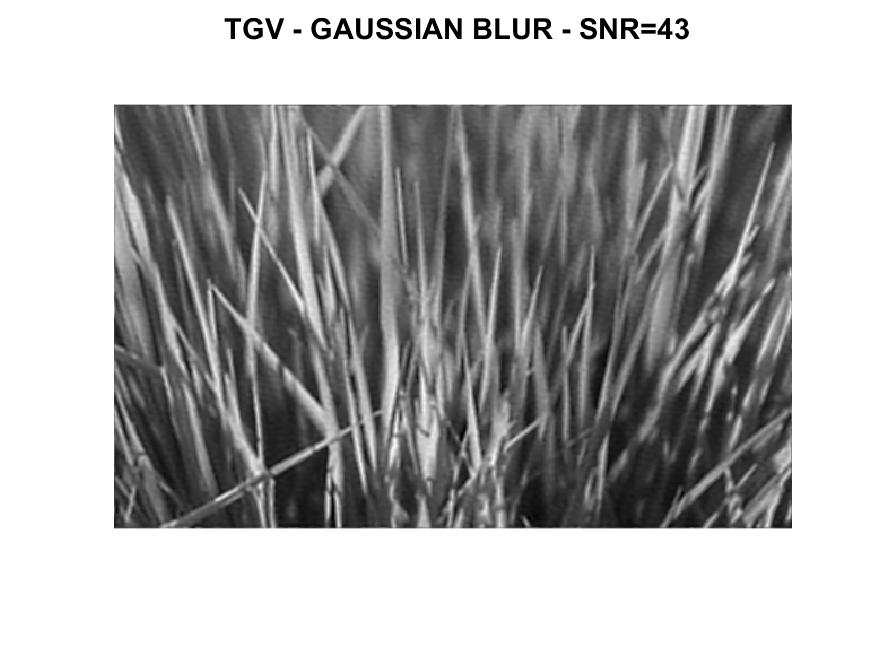}\hspace{-.6cm}
\includegraphics[width=0.25\columnwidth]{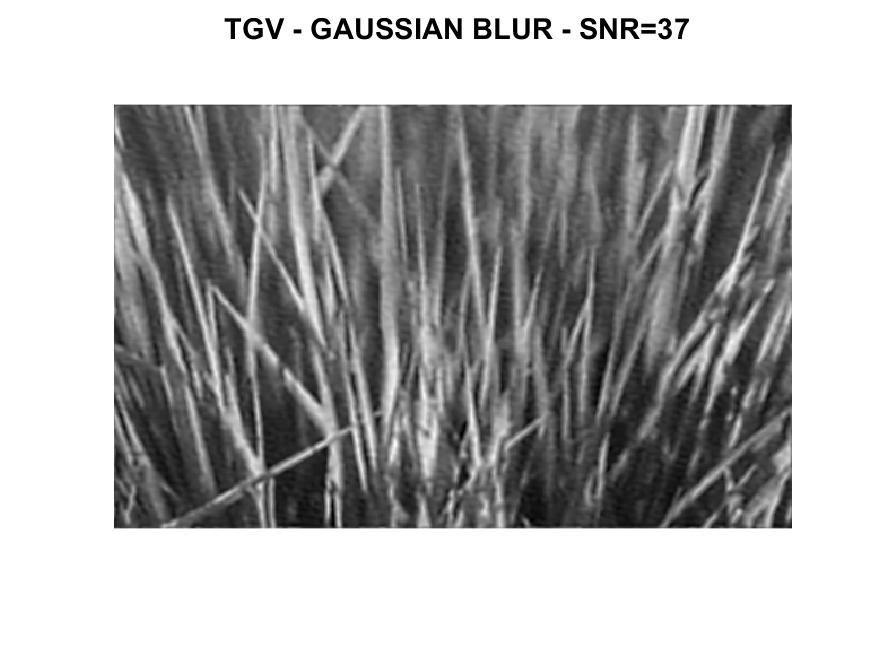} \\[-3mm]
\vspace{-9pt}
\caption{Test problem \texttt{grass}: images restored with DTGV$^2$ (\textbf{top}) and TGV$^2$ (\textbf{bottom}).\label{fig:grassrec}}
\end{center}
\end{figure}

\begin{figure}[htbp]
\begin{center}
    \includegraphics[width=0.25\columnwidth]{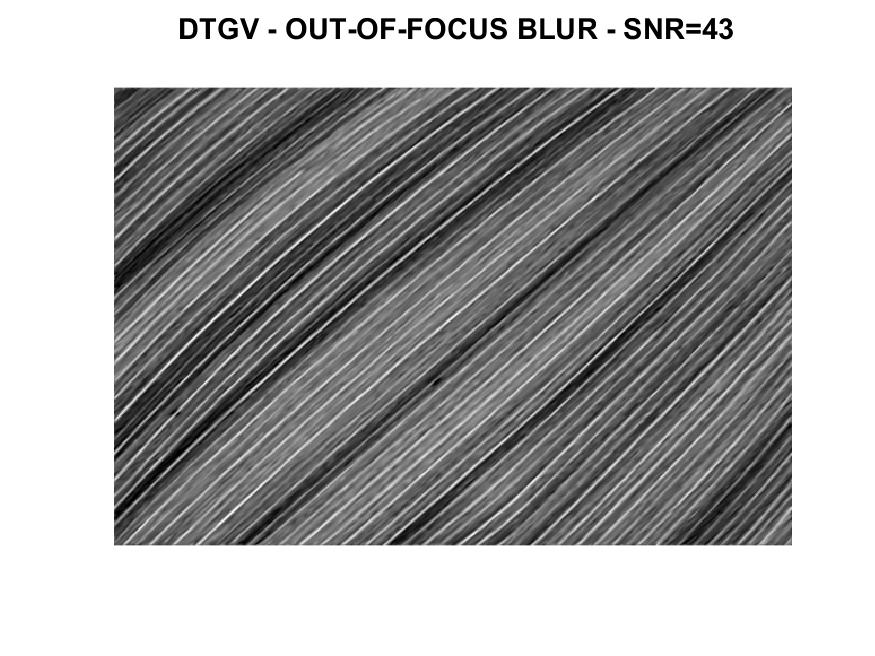}\hspace{-.6cm}
    \includegraphics[width=0.25\columnwidth]{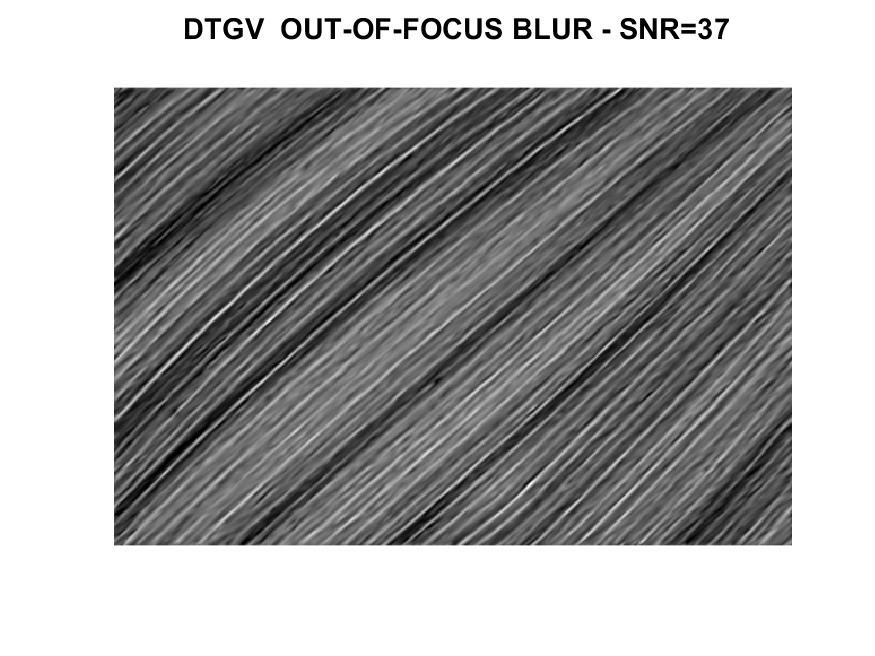}\hspace{-.6cm}
    \includegraphics[width=0.25\columnwidth]{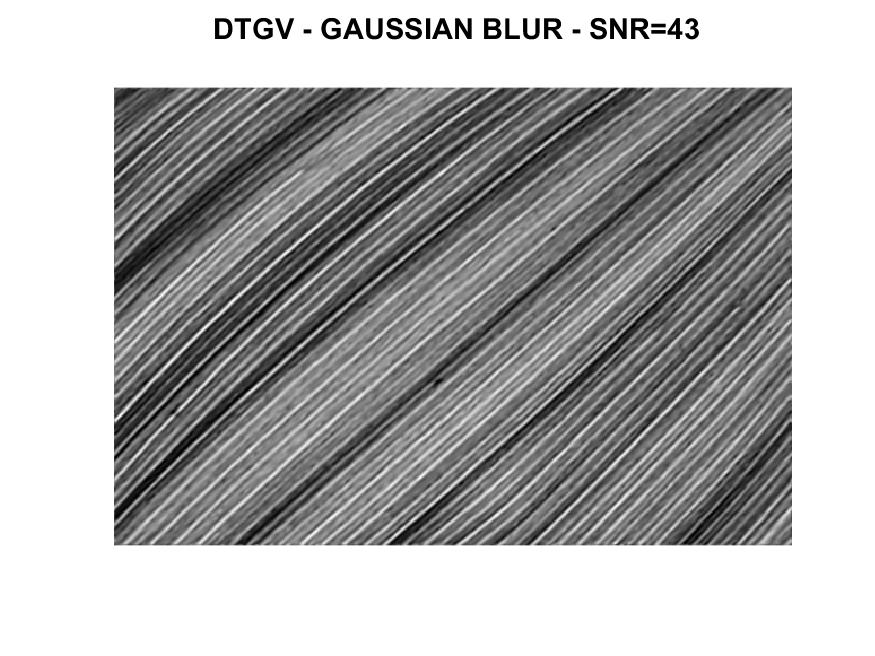}\hspace{-.6cm}
    \includegraphics[width=0.25\columnwidth]{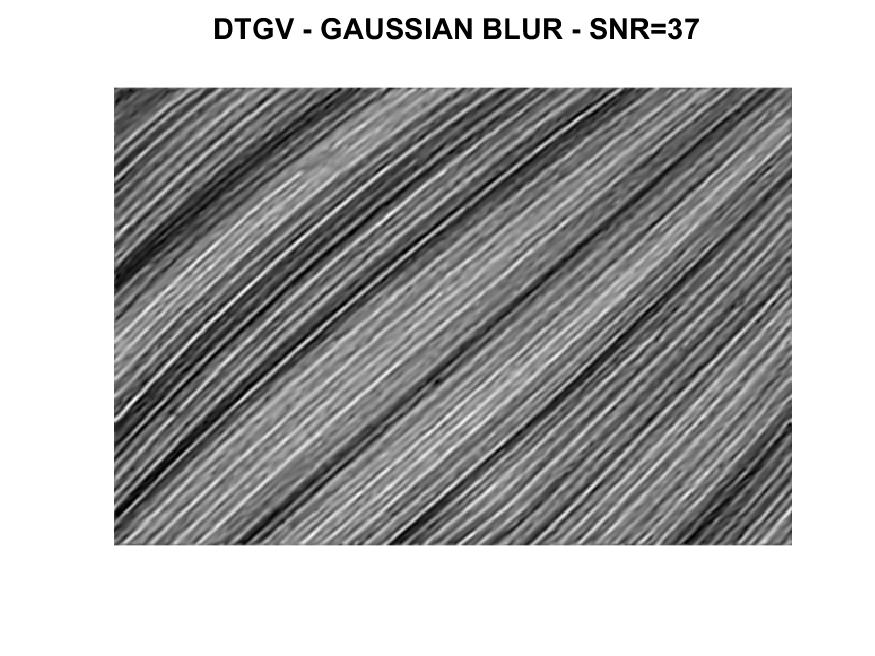} \\[-3mm]
    \includegraphics[width=0.25\columnwidth]{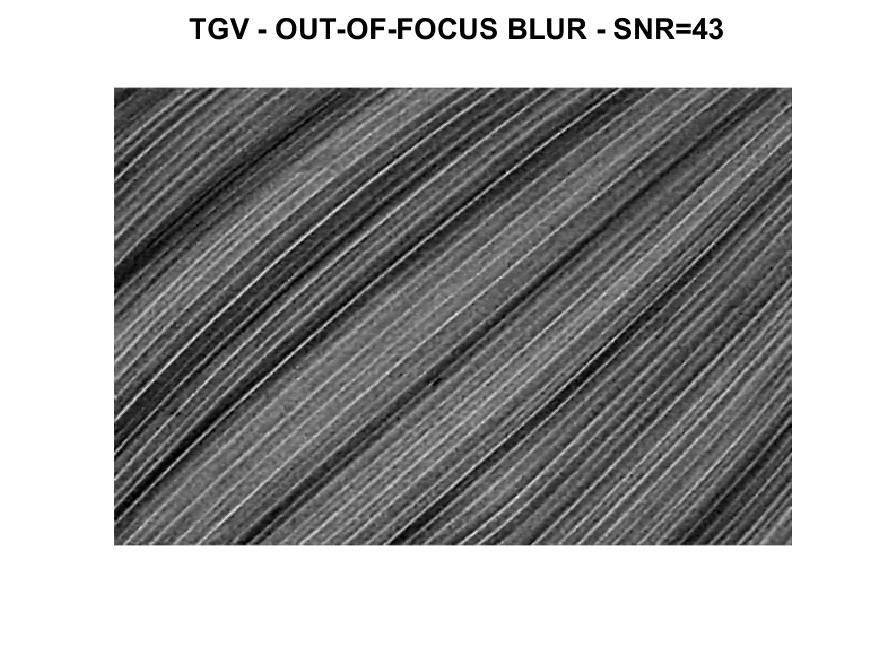}\hspace{-.6cm}
    \includegraphics[width=0.25\columnwidth]{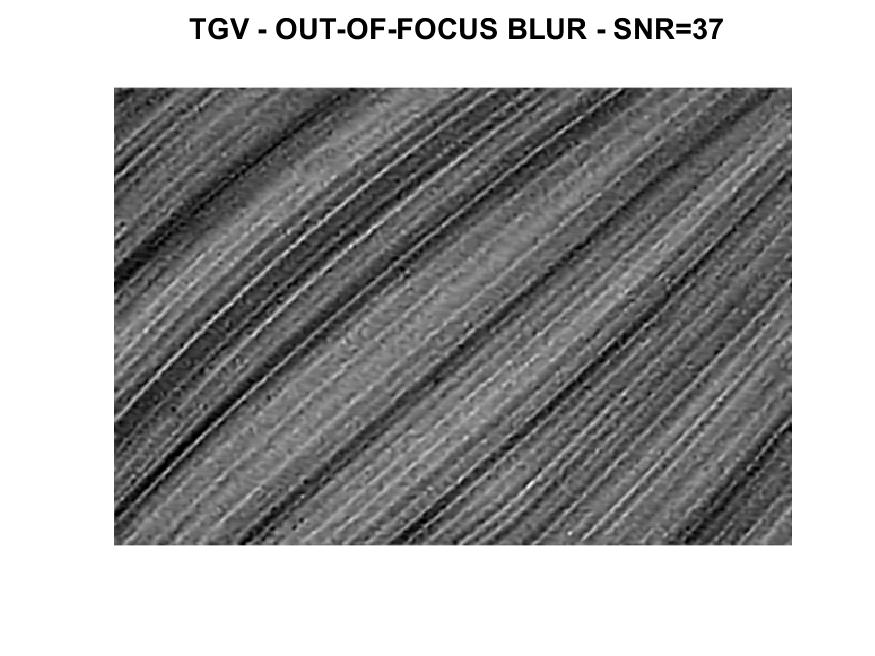}\hspace{-.6cm}
    \includegraphics[width=0.25\columnwidth]{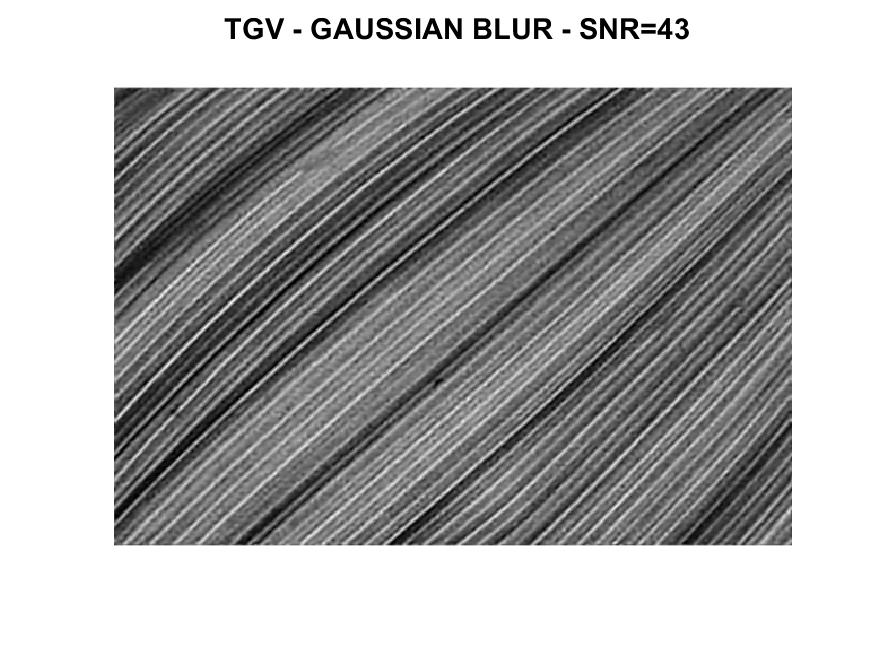}\hspace{-.6cm}
    \includegraphics[width=0.25\columnwidth]{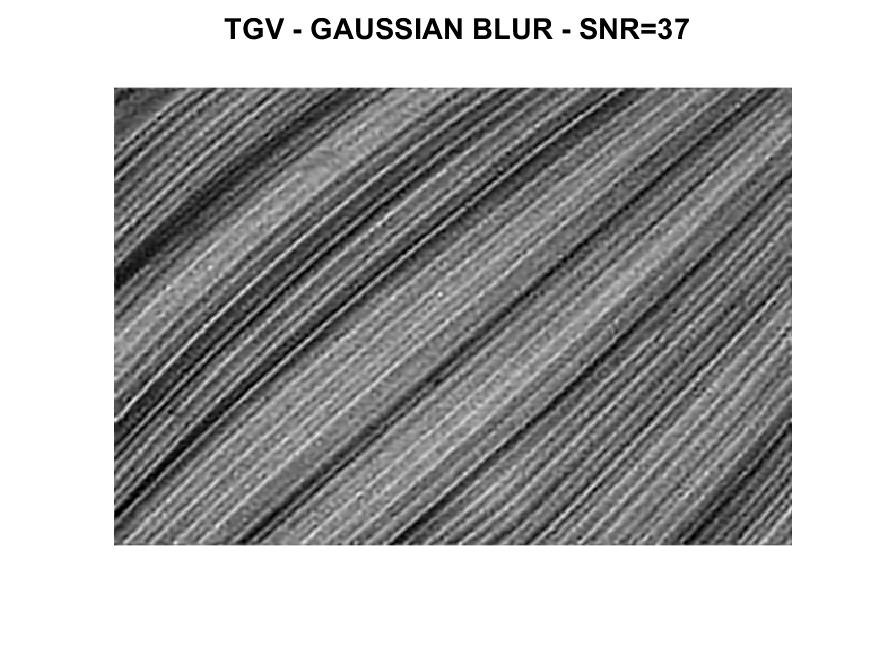} \\[-6mm]
\end{center}
\vspace{-9pt}
\caption{Test problem \texttt{leaves}: images restored with DTGV$^2$ (\textbf{top}) and TGV$^2$ (\textbf{bottom}).\label{fig:leavesrec}}
\end{figure}

\begin{figure}[htbp]
\begin{center}
    \includegraphics[width=0.25\columnwidth]{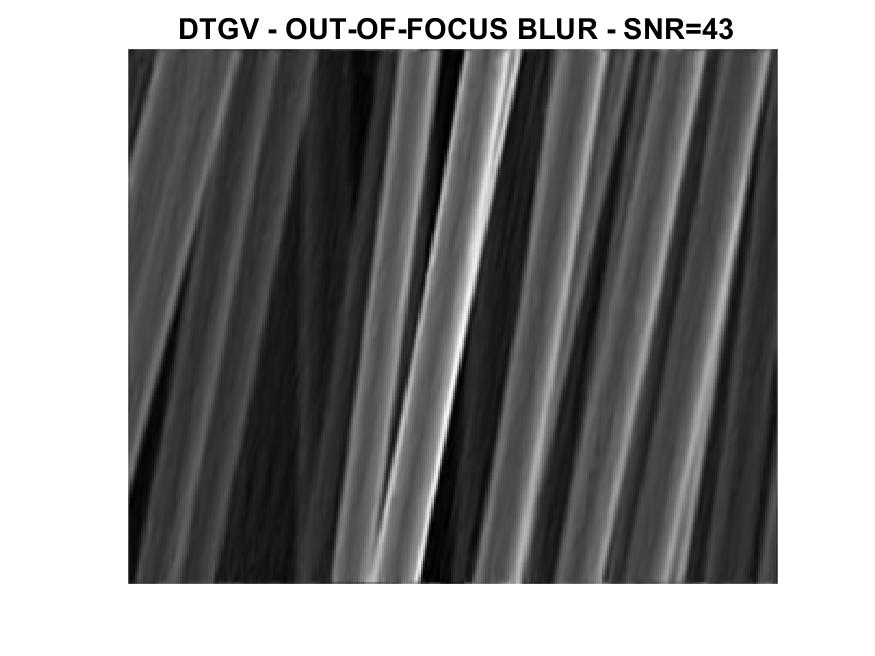}\hspace{-.6cm}
    \includegraphics[width=0.25\columnwidth]{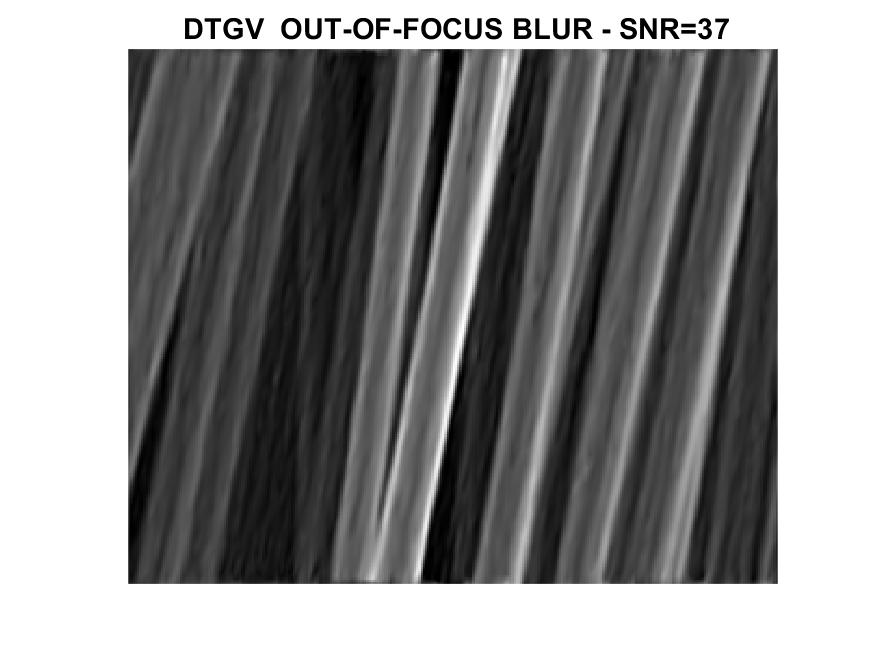}\hspace{-.6cm}
    \includegraphics[width=0.25\columnwidth]{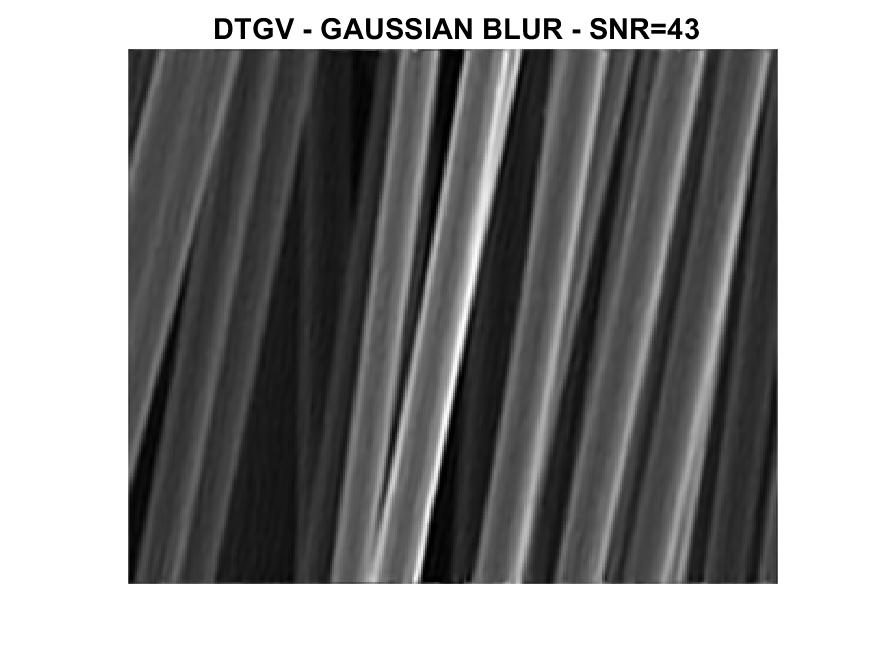}\hspace{-.6cm}
    \includegraphics[width=0.25\columnwidth]{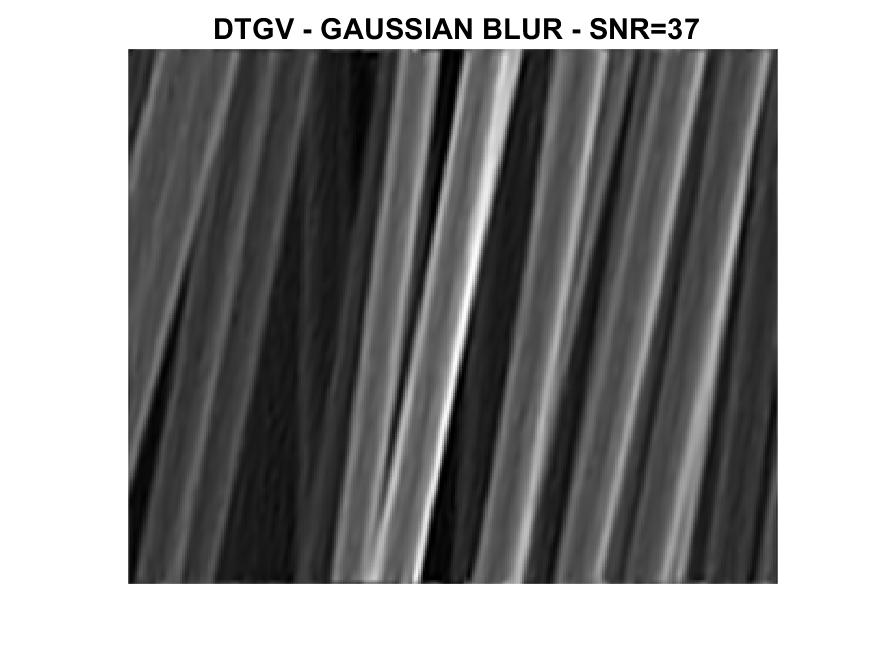} \\[-3mm]
    \includegraphics[width=0.25\columnwidth]{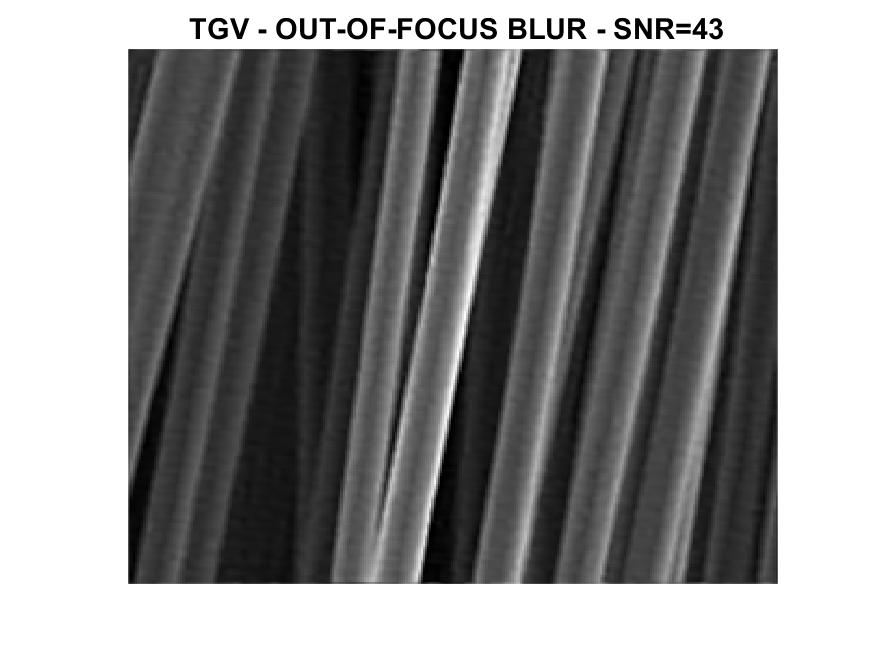}\hspace{-.6cm}
    \includegraphics[width=0.25\columnwidth]{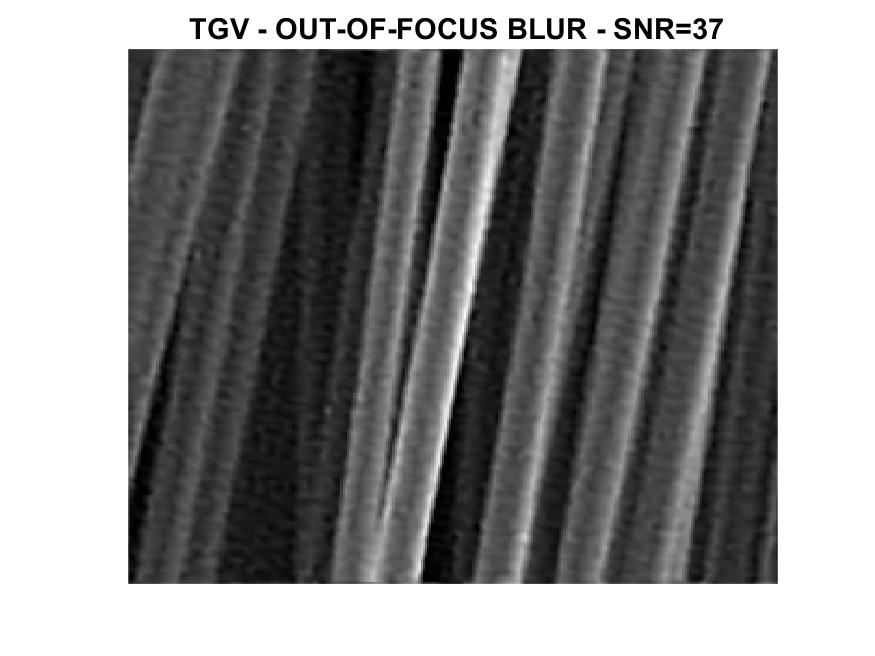}\hspace{-.6cm}
    \includegraphics[width=0.25\columnwidth]{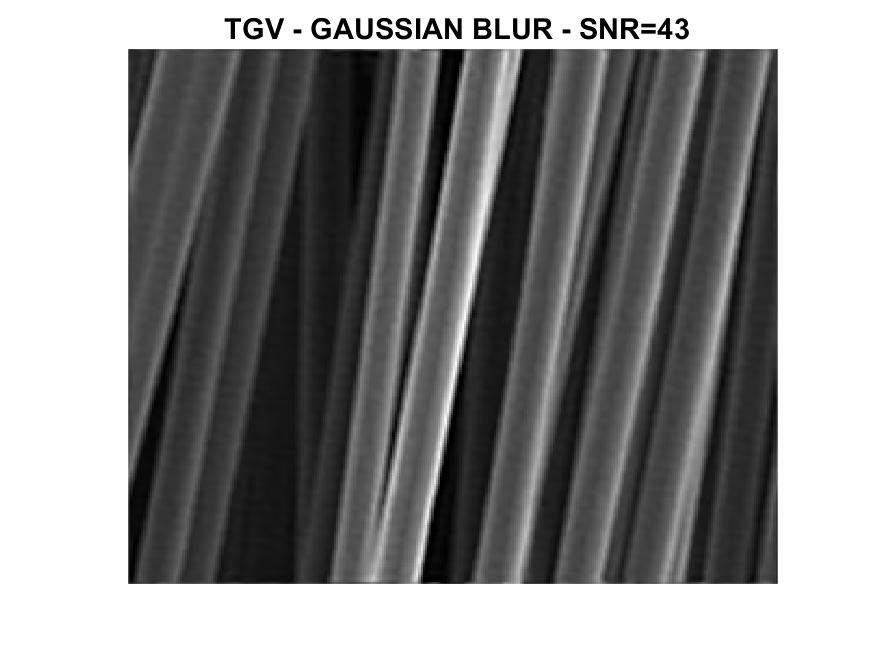}\hspace{-.6cm}
    \includegraphics[width=0.25\columnwidth]{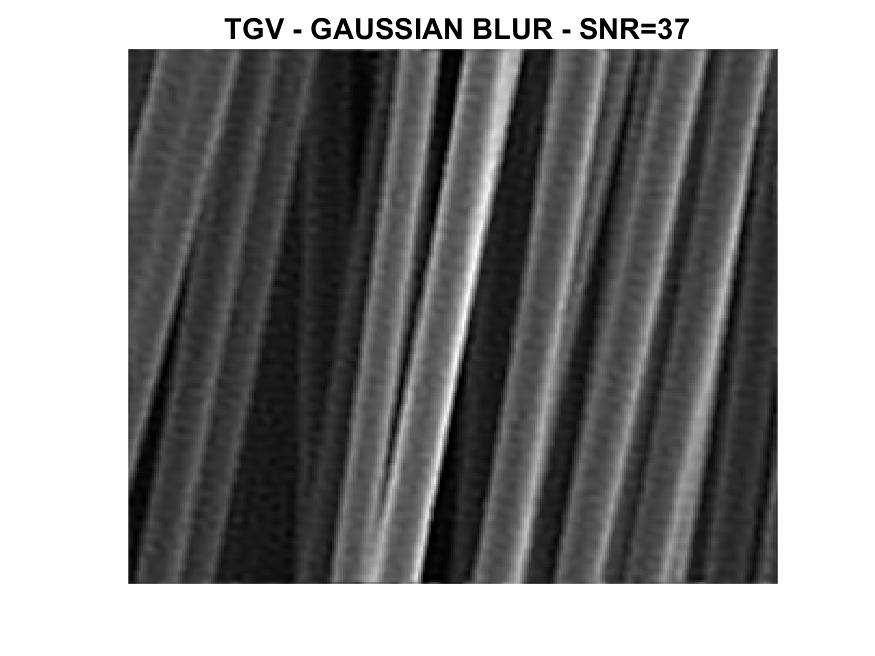} \\[-6mm]
\end{center}
\vspace{-9pt}
\caption{Test problem \texttt{carbon}: images restored with DTGV$^2$ (\textbf{top}) and TGV$^2$ (\textbf{bottom}).\label{fig:carbonrec}}
\end{figure}

\begin{figure}[htbp]
\begin{center}
    \includegraphics[width=0.25\columnwidth]{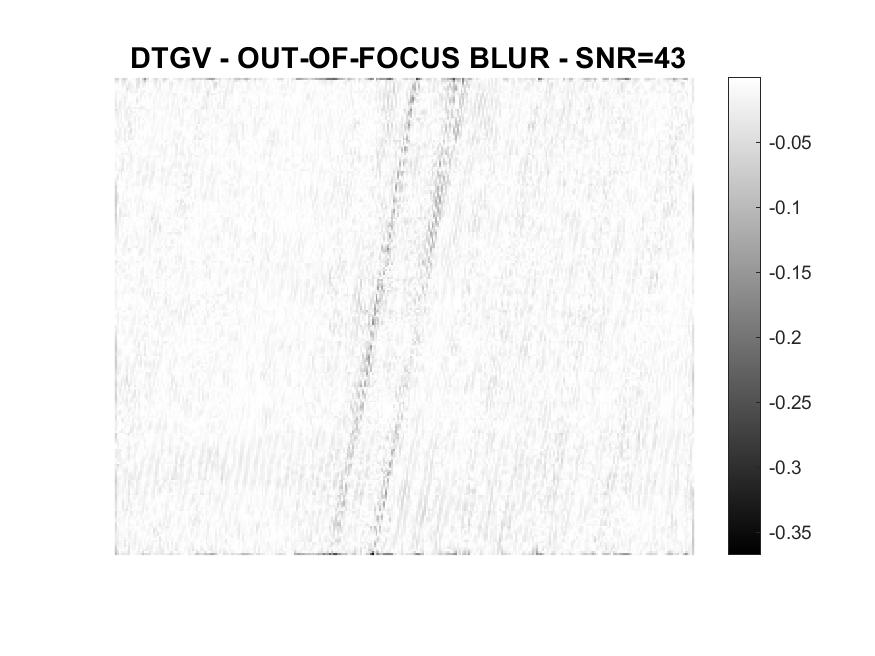}\hspace{-.4cm}
    \includegraphics[width=0.25\columnwidth]{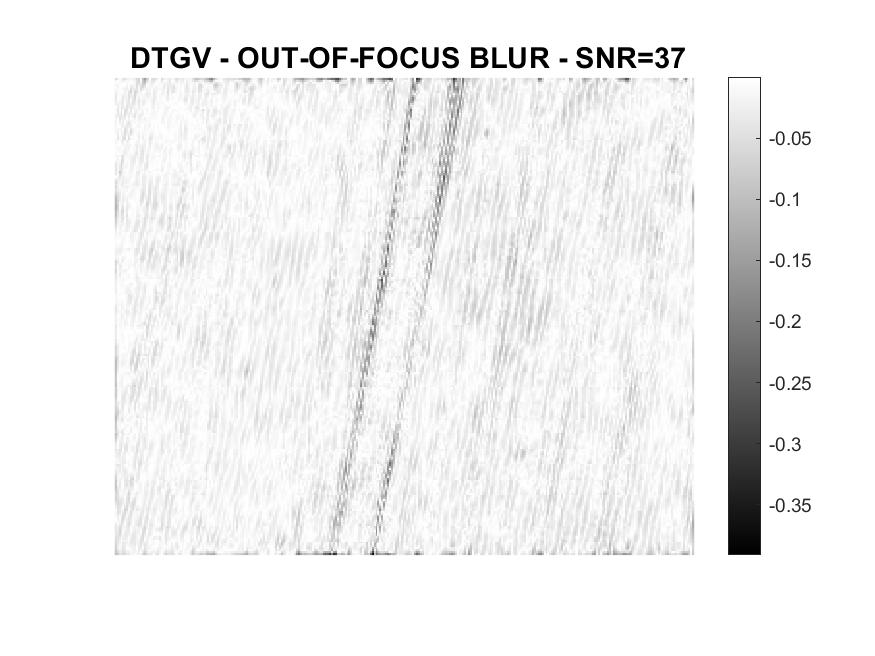}\hspace{-.4cm}
    \includegraphics[width=0.25\columnwidth]{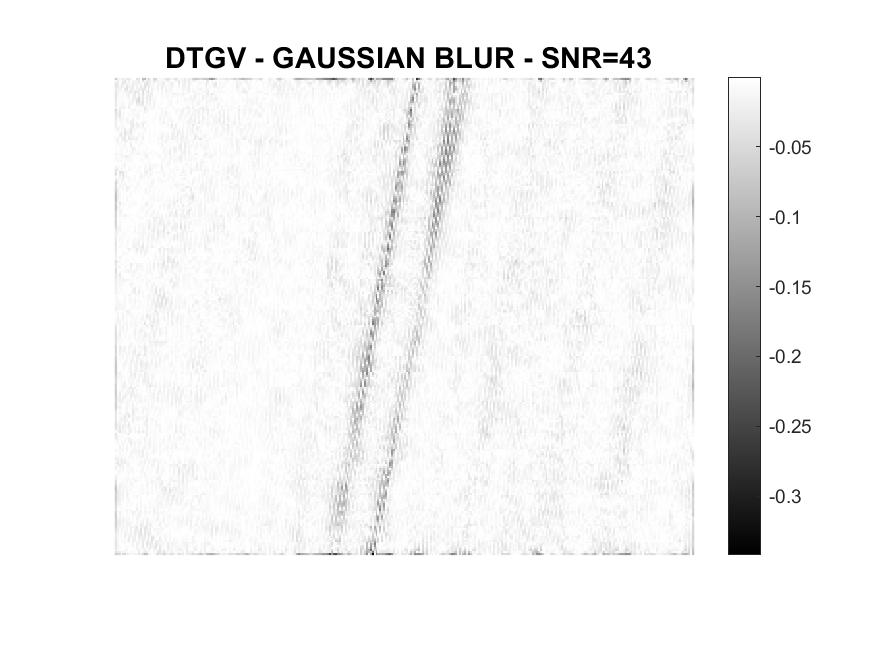}\hspace{-.4cm}
    \includegraphics[width=0.25\columnwidth]{DTGVimage4blur3highdiff} \\[-3mm]
    \includegraphics[width=0.25\columnwidth]{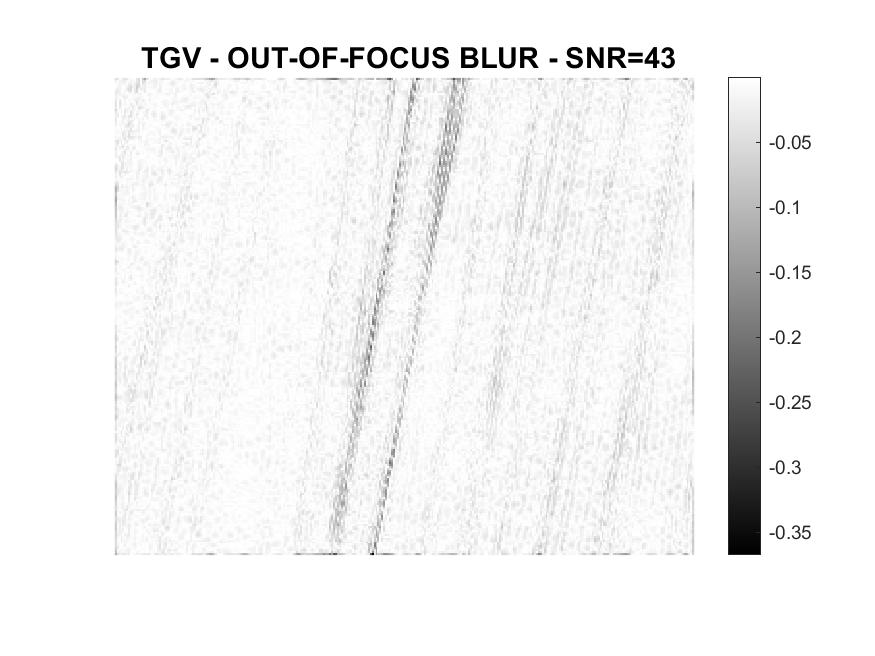}\hspace{-.4cm}
    \includegraphics[width=0.25\columnwidth]{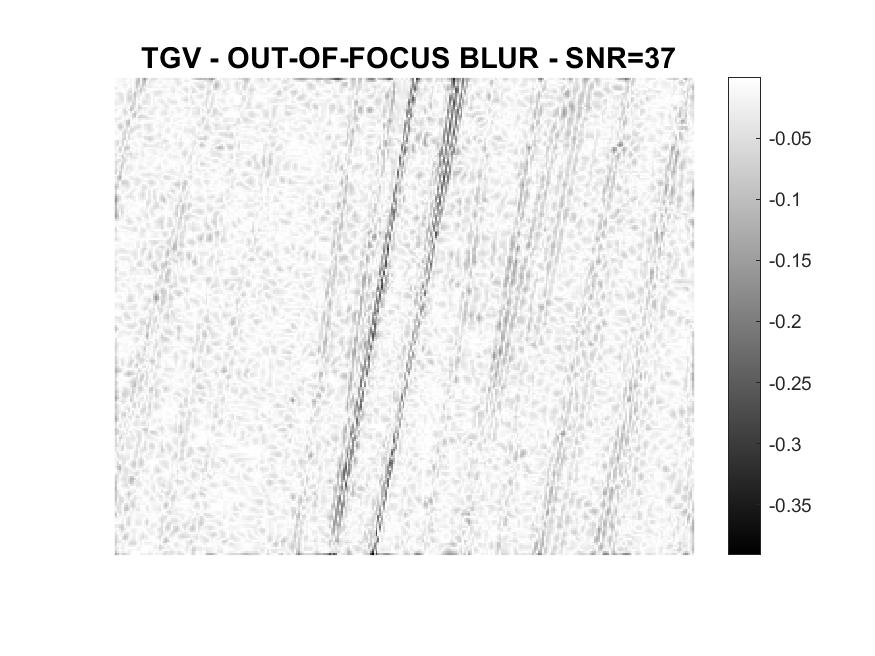}\hspace{-.4cm}
    \includegraphics[width=0.25\columnwidth]{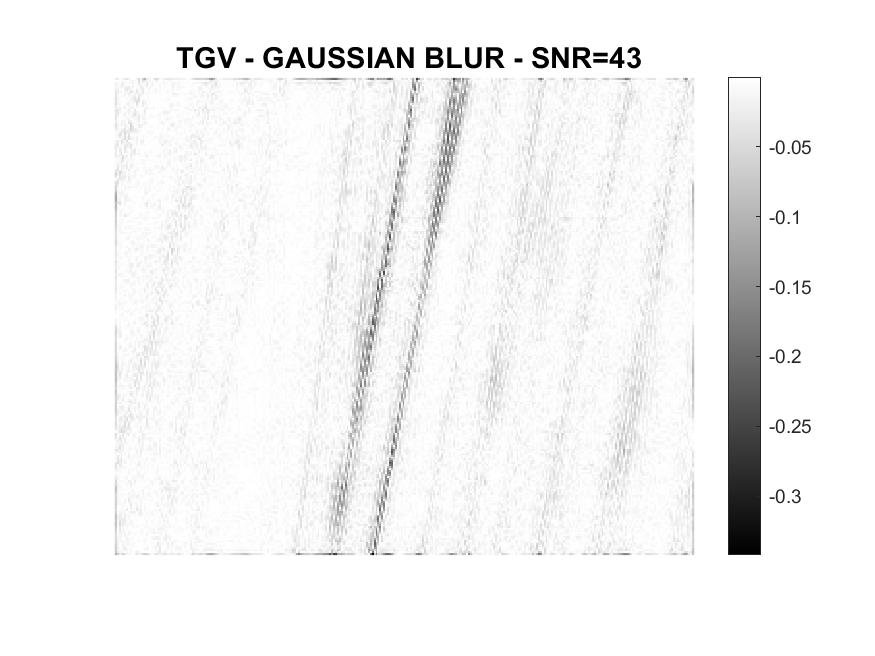}\hspace{-.4cm}
    \includegraphics[width=0.25\columnwidth]{TGVimage4blur3highdiff} \\[-3mm]
\end{center}
\vspace{-9pt}
\caption{Test problem \texttt{carbon}: difference images with DTGV$^2$ (\textbf{top}) and TGV$^2$ (\textbf{bottom}).\label{fig:carbondiff}} 
\end{figure}

The numerical results are summarized in Table~\ref{tab:1}, where the RMSE, the Improved Signal to Noise Ratio (ISNR) \cite{banham:1997}, and the structural similarity (SSIM) index \cite{wang:2004} are used to give a quantitative evaluation of the quality of the restorations. As a measure of the computational cost, the number of iterations and the time in seconds are reported. \mbox{Table~\ref{tab:1}} also shows, for each test problem, the values of the regularization parameter $\lambda$. The restored images are shown in  Figures~\ref{fig:phantomrec}--\ref{fig:carbonrec}. For the \texttt{carbon} test problem, Figure~\ref{fig:carbondiff} shows the error images, i.e., the images obtained as the absolute difference between the original image and the restored one. The values of the pixels of the error images have been scaled in the range $[m,M]$ where $m$ and $M$ are the minimum and maximum pixel value of the DTGV$^2$ and TGV$^2$ error images. 

From the results, it is evident that the DTGV$^2$ model outperforms the TGV$^2$ one in terms of quality of the restoration. A visual inspection of {the figures} shows that the DTGV$^2$ regularization is very effective in removing the noise, while for high noise levels the TGV$^2$ reconstructions still exhibit noise artifacts. Finally, by observing the ``Iters'' column of the table, we can conclude that, on average, the TGV$^2$ regularization requires less ADMM iterations to achieve a relative change in the restoration that is below the fixed threshold. However, the computational time per iteration is very small and also ADMM for the KL-DGTV$^2$ regularization is efficient.

\begin{figure}[htbp]
\begin{center}
    \includegraphics[width=0.24\textwidth]{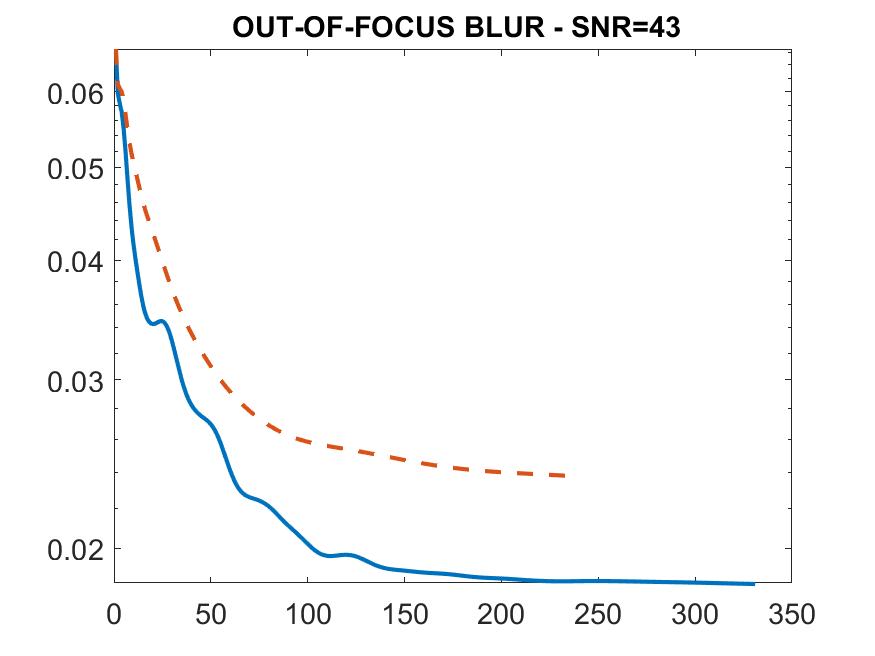}
    \includegraphics[width=0.24\textwidth]{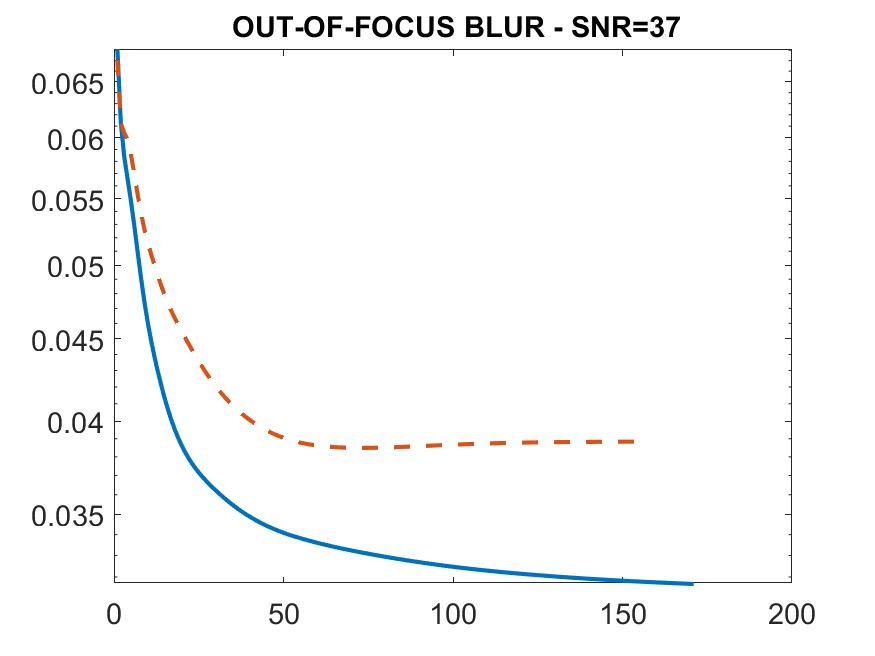}
    \includegraphics[width=0.24\textwidth]{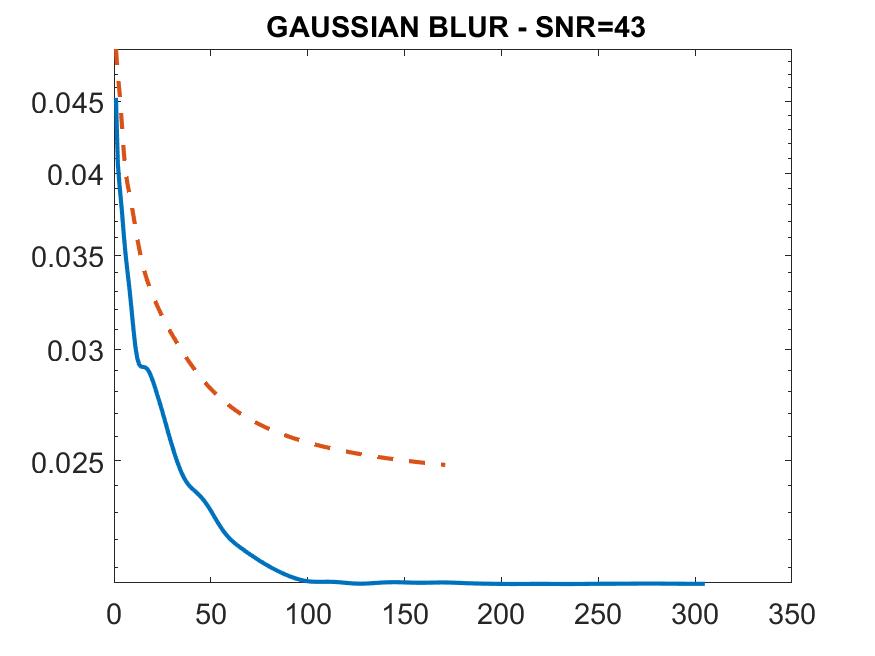}
    \includegraphics[width=0.24\textwidth]{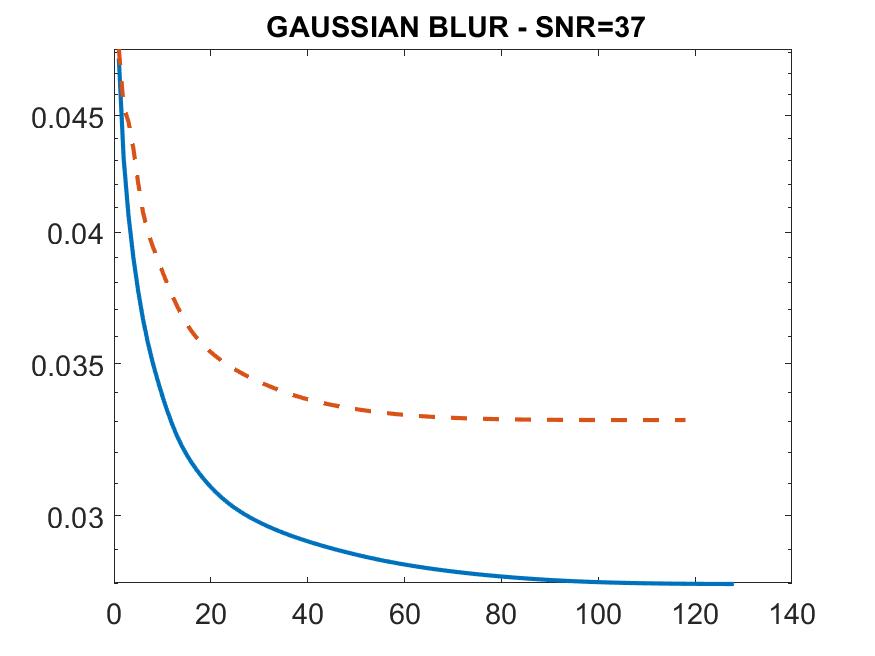}\\[-3mm]
\end{center}
\vspace{-9pt}
\caption{Test problem \texttt{carbon}: RMSE history for the KL-DGTV$^2$ (continuous line) and KL-TGV$^2$ (dashed line) models.\label{fig:phantomerr}}
\end{figure}

Finally, to illustrate the behaviour of ADMM, in Figure \ref{fig:phantomerr} we plot the RMSE history for the \texttt{carbon} test problem. A similar RMSE behaviour has been observed in all the numerical experiments.

%
%
\begin{table}[t!]
\begin{center}
    \caption{Numerical results for the test problems.\label{tab:1}}
    {\small
        \begin{tabular}{cclcccccc} 
            \toprule
            \textbf{Blur} & \textbf{SNR}  &  \textbf{Model}     & \boldmath{$\lambda$} & \textbf{RMSE}  &  \textbf{ISNR} &  \textbf{MSSIM} &  \textbf{Iters}  & \textbf{Time       } \\
            \midrule
            \multicolumn{9}{c}{\texttt{phantom}} \\
            \midrule
            \multirow{4}{*}{Out-of-focus}
            & \multirow{2}{*}{43}  & DTGV & 57.5  & 2.2558$\times 10^{-2}$ & 9.5472 & 9.3007$\times 10^{-1}$ & 86 & 10.95 \\ 
            &                      & TGV  & 275   & 2.8043$\times 10^{-2}$ & 7.6568 & 8.9887$\times 10^{-1}$ & 89 & 11.33 \\
            \cline{2-9}
            & \multirow{2}{*}{37} & DTGV & 3.25  & 3.7573$\times 10^{-2}$ & 7.4431 & 8.5823$\times 10^{-1}$ & 122 & 15.45 \\
            &                     & TGV  & 22.5  & 4.1719$\times 10^{-2}$ & 6.5339 & 8.4061$\times 10^{-1}$ & 52  & 6.64 \\
            \midrule
            \multirow{4}{*}{Gaussian}
            & \multirow{2}{*}{43} & DTGV & 25 & 1.5530$\times 10^{-2}$ & 9.1966 & 9.7829$\times 10^{-1}$ & 56 & 7.17 \\
            &                     & TGV  & 100& 1.8100$\times 10^{-2}$ & 7.8667 & 9.7200$\times 10^{-1}$ & 45 & 5.76 \\
            \cline{2-9}
            & \multirow{2}{*}{37} & DTGV & 3    & 2.5498$\times 10^{-2}$ & 9.0841 & 9.2994$\times 10^{-1}$ & 90 & 11.41 \\
            &                     & TGV  & 17.5 & 3.0674$\times 10^{-2}$ & 7.4788 & 9.0199$\times 10^{-1}$ & 53 & 6.76 \\
            \midrule
            \multicolumn{9}{c}{\texttt{grass}}\\
            \midrule
            \multirow{4}{*}{Out-of-focus}
            & \multirow{2}{*}{43} & DTGV & 60   & 3.6313$\times 10^{-2}$ & 7.7364 & 8.7262$\times 10^{-1}$ & 136 & 15.55 \\
            &                     & TGV  & 550  & 3.6575$\times 10^{-2}$ & 7.6738 & 8.7188$\times 10^{-1}$ & 179 & 20.39 \\
            \cline{2-9}
            & \multirow{2}{*}{37} & DTGV & 50  & 5.6164$\times 10^{-2}$ & 4.7390 & 7.6165$\times 10^{-1}$ & 160 & 18.56 \\
            &                     & TGV  & 55  & 5.7604$\times 10^{-2}$ & 4.5191 & 7.4566$\times 10^{-1}$ & 72  & 8.31 \\
            \midrule
            \multirow{4}{*}{Gaussian}
            & \multirow{2}{*}{43} & DTGV & 65  & 2.9883$\times 10^{-2}$ & 6.3343 & 9.2764$\times 10^{-1}$ & 106 & 12.08 \\
            &                     & TGV  & 650 & 3.0814$\times 10^{-2}$ & 6.0676 & 9.2523$\times 10^{-1}$ & 136 & 15.48 \\
            \cline{2-9}
            & \multirow{2}{*}{37} & DTGV & 5.5 & 4.2274$\times 10^{-2}$ & 4.7973 & 8.5615$\times 10^{-1}$ & 98 & 11.13 \\
            &                     & TGV  & 35  & 4.3936$\times 10^{-2}$ & 4.4624 & 8.4795$\times 10^{-1}$ & 54 & 6.18 \\
            \midrule
            \multicolumn{9}{c}{\texttt{leaves}}\\
            \midrule
            \multirow{4}{*}{Out-of-focus}
            & \multirow{2}{*}{43} & DTGV & 125 & 6.2767$\times 10^{-2}$ & 7.4978 & 8.2099$\times 10^{-1}$ & 251 & 31.18 \\
            &                     & TGV  & 1100& 8.2397$\times 10^{-2}$ & 5.1342 & 7.1557$\times 10^{-1}$ & 435 & 53.74 \\
            \cline{2-9}
            & \multirow{2}{*}{37} & DTGV & 12.5 & 9.5597$\times 10^{-2}$ & 4.1497 & 6.3065$\times 10^{-1}$ & 257 & 31.87 \\
            &                     & TGV  & 90   & 1.1874$\times 10^{-1}$ & 2.2665 & 4.3294$\times 10^{-1}$ & 113 & 14.03 \\
            \midrule
            \multirow{4}{*}{Gaussian}
            & \multirow{2}{*}{43} & DTGV & 150  & 7.3332$\times 10^{-2}$ & 4.8675 & 7.7456$\times 10^{-1}$ & 236 & 29.13 \\
            &                     & TGV  & 1750 & 8.0857$\times 10^{-2}$ & 4.0190 & 7.3001$\times 10^{-1}$ & 380 & 46.77 \\
            \cline{2-9}
            & \multirow{2}{*}{37} & DTGV & 12.5 & 9.0999$\times 10^{-2}$ & 3.3907 & 6.6469$\times 10^{-1}$ & 148 & 18.36 \\
            &                     & TGV  & 100  & 1.0308$\times 10^{-1}$ & 2.3081 & 5.6534$\times 10^{-1}$ & 103 & 12.85 \\
            \midrule
            \multicolumn{9}{c}{\texttt{carbon}}\\
            \midrule
            \multirow{4}{*}{Out-of-focus}
            & \multirow{2}{*}{43} & DTGV & 150 & 1.8360$\times 10^{-2}$ & 1.2830$\times 10^{1}$ & 9.4734$\times 10^{-1}$ & 331 & 13.78 \\
            &                     & TGV  & 850 & 2.3825$\times 10^{-2}$ & 1.0567$\times 10^{1}$ & 9.3671$\times 10^{-1}$ & 233 & 9.73 \\
            \cline{2-9}
            & \multirow{2}{*}{37} & DTGV & 20  & 3.1682$\times 10^{-2}$ & 8.2416 & 8.6294$\times 10^{-1}$ & 171 & 7.07\\
            &                     & TGV  & 150 & 3.8840$\times 10^{-2}$ & 6.4723 & 8.2237$\times 10^{-1}$ & 155 & 6.55\\
            \midrule
            \multirow{4}{*}{Gaussian}
            & \multirow{2}{*}{43} & DTGV & 250 & 2.0453$\times 10^{-2}$ & 8.6178 & 9.5974$\times 10^{-1}$ & 305 & 12.53 \\
            &                     & TGV  & 950 & 2.4839$\times 10^{-2}$ & 6.9302 & 9.5698$\times 10^{-1}$ & 171 & 7.12 \\
            \cline{2-9}
            & \multirow{2}{*}{37} & DTGV & 15  & 2.7995$\times 10^{-2}$ & 6.2017 & 9.3007$\times 10^{-1}$ & 128 & 5.36 \\
            &                     & TGV  & 150 & 3.3061$\times 10^{-2}$ & 4.7572 & 8.9690$\times 10^{-1}$ & 118 & 4.73 \\
            \bottomrule
        \end{tabular}
    }
\end{center}
\end{table}

%
%

\section{Conclusions\label{sec:conclusions}}

We dealt with the use of the Directional TGV regularization in the case of directional images corrupted by Poisson noise. We presented the KL-DTGV$^2$ model and introduced a two-block ADMM version for its minimization. Finally, we proposed an effective strategy for the estimation of the main direction of the image. Our numerical experiments show that for Poisson noise the DTGV$^2$ regularization provides superior restoration performance compared with the standard TGV$^2$ regularization, thus remarking the importance of taking into account the texture structure of the image. A crucial ingredient for the success of the model was the proposed direction estimation strategy, which proved to be more reliable than those proposed in the literature.

Possible future work includes the use of space-variant regularization terms and the analysis of automatic strategies for the selection of the regularization parameters.

\end{document}